\documentclass[12pt]{article}
\usepackage[]{amsmath,amssymb}
\usepackage{amscd}
\usepackage{latexsym}
\usepackage{cite}

\newtheorem{definition}{Definition}[section]
\newtheorem{theorem}[definition]{Theorem}
\newtheorem{lemma}[definition]{Lemma}
\newtheorem{corollary}[definition]{Corollary}

\newtheorem{example}[definition]{Example}

\newtheorem{note}[definition]{Note}

\newtheorem{proposition}[definition]{Proposition}

\def\N{\mathbb N}

\def\Z{\mathbb Z}
\def\F{\mathbb F}

\begin{document}
\title{\bf  
Evaluation modules for the \\
$q$-tetrahedron algebra
}
\author{
Tatsuro Ito, Hjalmar Rosengren, and
Paul Terwilliger}
\date{}

\maketitle
\begin{abstract}
Let $\F$ denote an algebraically closed field, and fix a nonzero
$q \in \F$ that is not a root of unity.
We consider the $q$-tetrahedron algebra 
$\boxtimes_q$ over $\F$.
It is known that each finite-dimensional irreducible 
$\boxtimes_q$-module of type 1 is a tensor product of
evaluation modules. This paper contains a comprehensive
description of the evaluation modules for $\boxtimes_q$.
This description includes the following topics.
Given an evaluation module $V$ for $\boxtimes_q$,
we display 24
bases for $V$ that we find
attractive. For each basis
we give the matrices
that represent the $\boxtimes_q$-generators. We give
the transition matrices between certain pairs of bases
among the 24. It is known that the    
cyclic group $\Z_4$ acts on $\boxtimes_q$ as a group of
automorphisms. We describe what happens when $V$ is
twisted via an element of $\Z_4$. We discuss how
evaluation modules for $\boxtimes_q$ are related to
Leonard pairs of $q$-Racah type.

 \end{abstract}

\section{Introduction}
The $q$-tetrahedron algebra $\boxtimes_q$ was introduced in
\cite{qtet}. This algebra is associative, 
noncommutative, and infinite-dimensional.
It is defined by
 generators and relations.
There are eight generators, and it is natural to
identify each of these
with an orientation on an edge in
a tetrahedron.
From this point of view the generating set looks
as follows.
In the tetrahedron,
a pair of opposite edges are each oriented in both directions.
The four remaining edges are each oriented in one direction,
to create a directed 4-cycle.
Thus the cyclic group $\Z_4$ acts transitively
on the vertex set of the tetrahedron, in a manner that
preserves the set of edge-orientations.
The relations in $\boxtimes_q$ are described as follows.
For each doubly oriented edge of the tetrahedron,
the product of the two 
edge-orientations is 1.
For each pair of edge-orientations that create a directed 2-path
involving three distinct vertices,
 these
edge-orientations 
satisfy a $q$-Weyl relation. 
For each pair of edges in the tetrahedron that are opposite
and singly oriented,
the two associated 
edge-orientations 
satisfy the cubic $q$-Serre relations.
By construction, the $\Z_4$ action on the tetrahedron
induces 
a $\Z_4$ action on $\boxtimes_q$ as a group
of automorphisms.  

\medskip
\noindent
We will be discussing the quantum enveloping
algebra 
$U_q(\mathfrak{sl}_2)$,
the loop algebra
$U_q(L(\mathfrak{sl}_2))$
\cite[Section~8]{qtet},
and an
algebra $\mathcal A_q$ called 
the positive part of
 $U_q(\widehat{\mathfrak{sl}}_2)$
\cite[Definition~9.1]{qtet}.
These algebras are related to 
 $\boxtimes_q$  
in the following way.
Each face of the tetrahedron is surrounded by
three edges, of which two are singly oriented and one is doubly oriented.
The resulting four edge-orientations 
generate a subalgebra of $\boxtimes_q$ that is
isomorphic to 
$U_q(\mathfrak{sl}_2)$
\cite[Proposition~7.4]{qtet},
\cite[Proposition~4.3]{miki}.
Upon removing one doubly oriented edge from the tetrahedron,
the remaining
six edge-orientations generate a subalgebra of
$\boxtimes_q$ that is isomorphic to 
$U_q(L(\mathfrak{sl}_2))$
\cite[Proposition~4.3]{miki}.
For each pair of edges in the tetrahedron that are opposite
and singly oriented, the two associated edge-orientations
generate a subalgebra of $\boxtimes_q$ that is isomorphic to
$\mathcal A_q$
  \cite[Proposition~4.1]{miki}.

\medskip
\noindent The above containments reveal
a close relationship between
the representation theories of
$\boxtimes_q$, 
$U_q(L(\mathfrak{sl}_2))$, and $\mathcal A_q$.
Before discussing the details, we comment on
$U_q(L(\mathfrak{sl}_2))$.
In \cite{chari}, Chari and Pressley
classify up to isomorphism the finite-dimensional irreducible  
$U_q(L(\mathfrak{sl}_2))$-modules.
This classification involves a
bijection 
between the following two sets:
(i) the isomorphism classes of finite-dimensional irreducible
$U_q(L(\mathfrak{sl}_2))$-modules of type 1;
(ii) the polynomials in  
 one variable that have 
 constant coefficient 1.
The polynomial is called the Drinfel'd polynomial.

\medskip \noindent
The representation theories for $\boxtimes_q$ and
$U_q(L(\mathfrak{sl}_2))$ are related as follows.
Let $V$ denote a 
$\boxtimes_q$-module. 
Earlier we mentioned
a subalgebra of $\boxtimes_q$ that is isomorphic to 
$U_q(L(\mathfrak{sl}_2))$. Upon restricting the 
$\boxtimes_q$ action on $V$  to
this subalgebra, $V$ becomes
a $U_q(L(\mathfrak{sl}_2))$-module.
The restriction procedure yields a map 
from the set of 
$\boxtimes_q$-modules to the
set of $U_q(L(\mathfrak{sl}_2))$-modules.
By 
\cite[Remark~1.8]{nonnil} and
\cite[Remark~10.5]{qtet},
 this map induces a bijection between
the following two sets: (i) the isomorphism 
classes of 
finite-dimensional irreducible 
$\boxtimes_q$-modules of type 1; (ii) the
isomorphism classes of finite-dimensional
irreducible
$U_q(L(\mathfrak{sl}_2))$-modules of type 1
whose associated Drinfel'd polynomial does not vanish at $1$.
(We follow the normalization conventions from \cite{miki}).
In \cite{miki}, Miki extends the above bijective correspondence 
to include finite-dimensional modules that 
are not necessarily irreducible.
\medskip

\noindent
The representation theories for $\boxtimes_q$ and
$\mathcal A_q$ are related as follows. 
A finite-dimensional
$\mathcal A_q$-module
is called {\it NonNil} whenever
the two $\mathcal A_q$-generators 
are not nilpotent on the module
\cite[Definition~1.3]{nonnil}.
Let $V$ denote a $\boxtimes_q$-module.
Earlier we mentioned a subalgebra of
$\boxtimes_q$ that is isomorphic to
$\mathcal A_q$.
Upon restricting the $\boxtimes_q$ action
on $V$ to this subalgebra, $V$ becomes
a $\mathcal A_q$-module.
This yields a map
from the set of 
$\boxtimes_q$-modules to the
set of 
$\mathcal A_q$-modules.
By 
\cite[Remark~10.5]{qtet} this map induces a bijection between
the following two sets: (i) the isomorphism 
classes of 
finite-dimensional irreducible 
$\boxtimes_q$-modules of type 1; (ii) the
isomorphism classes of NonNil finite-dimensional
irreducible
$\mathcal A_q$-modules of type (1,1).

\medskip
\noindent 
We just related the representation theories
of 
$\boxtimes_q$ and 
$\mathcal A_q$.
To illuminate this relationship
we bring in the concept of a Leonard pair
\cite{terLS,
ter24,
ter:LPQR,
ter:LBUB,
terPA,
tersplit,
vidter}
and tridiagonal pair
\cite{TD00,
tdanduqsl2hat,
class,
TDpairKraw}.
 Roughly speaking, a Leonard pair
consists of two diagonalizable 
linear transformations of a finite-dimensional
vector space, each of which  acts in an
irreducible tridiagonal fashion
on an eigenbasis for the other one
\cite[Definition~1.1]{terLS}. The Leonard pairs 
are classified 
\cite{terLS,
terPA} 
and correspond
to 
the orthogonal polynomials that make up the terminating
branch of the Askey scheme
\cite{
koekoek}. A tridiagonal pair is a mild generalization
of a Leonard pair
\cite[Definition~1.1]{TD00}.
Let $V$ denote a NonNil finite-dimensional irreducible
$\mathcal A_q$-module of type $(1,1)$. Then
the two $\mathcal A_q$-generators
act on $V$ as a tridiagonal pair
\cite[Example~1.7]{TD00}.
A tridiagonal pair
obtained in this way is said to be $q$-geometric
\cite[Section~1]{class}.
Now the bijection from the previous paragraph amounts to
the following. Let $V$ denote a finite-dimensional
irreducible $\boxtimes_q$-module of type 1. Then
for any pair of edges in the tetrahedron
that are opposite 
and singly oriented,  the two associated
edge-orientations act on $V$ as a tridiagonal pair
of $q$-geometric type. Moreover, every tridiagonal
pair of $q$-geometric type is
 obtained in this way.
For more detail on this correspondence
see \cite[Section~2]{nonnil} and
\cite[Section~10]{qtet}.
In 
\cite{darren}, Funk-Neubauer obtains a similar
correspondence between
finite-dimensional irreducible $\boxtimes_q$-modules
and certain tridiagonal pairs of $q$-Hahn type.

\medskip
\noindent 
We mention some other results concerning
$\boxtimes_q$.
In \cite{qinv}, Ito and Terwilliger characterize the
 finite-dimensional irreducible $\boxtimes_q$-modules
using $q$-inverting pairs of linear transformations.
In \cite{qtetDRG}, Ito and Terwilliger display
an action of $\boxtimes_q$ on the standard module of any
distance-regular graph 
that is self-dual and has classical parameters with base $q^2$.
In \cite{jykim}, Joohyung Kim provides more details
about this $\boxtimes_q$ action.
In \cite{qracDRG}, Ito and Terwilliger obtain 
a similar $\boxtimes_q$ action
for certain distance-regular graphs of $q$-Racah type.
\medskip

\noindent Turning to the present paper,
our topic is a family
of finite-dimensional irreducible $\boxtimes_q$-modules of type 1,
called evaluation modules.
These modules are important for the following reason.
In 
\cite{chari},
Chari and Pressley show that 
each finite-dimensional irreducible
$U_q(L(\mathfrak{sl}_2))$-module of type 1 is a tensor product of
evaluation modules. 
Earlier we described how each finite-dimensional irreducible
$\boxtimes_q$-module of type 1 corresponds to a 
 finite-dimensional irreducible
$U_q(L(\mathfrak{sl}_2))$-module of type 1.
The tensor product structure survives the correspondence
\cite[Theorem~8.1]{miki}
and consequently
each finite-dimensional irreducible
$\boxtimes_q$-module of type 1
 is a tensor product of
evaluation modules
\cite[Sections~8,~9]{miki}.
\medskip

\noindent 
This paper contains a comprehensive
description of the evaluation modules for $\boxtimes_q$.
Hoping to keep this description accessible,
we avoid Hopf algebra theory and use
 only linear algebra.
Our description is roughly analogous to the description given in
\cite{tetmodules} concerning
the evaluation modules for the tetrahedron algebra.
\medskip

\noindent In our description, each evaluation module
for $\boxtimes_q$ gets a notation of the form
$\mathbf V_d (t)$. Here
$d$ is a positive integer, and $t$ is a nonzero scalar in
the underlying field that is not among
$\lbrace q^{d-2n+1}\rbrace_{n=1}^d$.
The $\boxtimes_q$-module 
$\mathbf V_d (t)$
has dimension $d+1$. On $\mathbf V_d(t)$, each of the
eight $\boxtimes_q$-generators is diagonalizable with
eigenvalues $\lbrace q^{d-2n}\rbrace_{n=0}^d$. The $\boxtimes_q$-module
$\mathbf V_d(t)$ 
is determined up to isomorphism by $d$ and $t$.
 We obtain several polynomial identities
that hold on $\mathbf V_d(t)$;
these identities involve the eight $\boxtimes_q$-generators and also
 $t$.
\medskip

\noindent 
We display 24 bases for 
$\mathbf V_d (t)$ that we find attractive.
These bases are described as follows.
For each permutation $i,j,k,\ell$ of the vertices of
the tetrahedron, we define a basis 
for $\mathbf V_d (t)$
denoted $\lbrack i,j,k,\ell\rbrack$.
This basis diagonalizes each $\boxtimes_q$-generator
involving the vertices $k$ and $\ell$. Moreover, the sum of
the basis vectors is an eigenvector for
 each $\boxtimes_q$-generator involving the vertex $j$.
We display the matrices that represent the
 eight $\boxtimes_q$-generators with respect to 
$\lbrack i,j,k,\ell\rbrack$.
We also 
give the transition matrices 
from 
the basis $\lbrack i,j,k,\ell\rbrack$ 
to each of the bases
\begin{eqnarray*}
\lbrack j,i,k,\ell\rbrack,
\qquad \qquad
\lbrack i,k,j,\ell\rbrack,
\qquad \qquad
\lbrack i,j,\ell,k\rbrack.
\end{eqnarray*}
The first transition matrix is diagonal, the
second one is lower triangular, and the third one
is the identity matrix reflected about a vertical axis.
\medskip

\noindent  Recall that the group $\Z_4$
acts on $\boxtimes_q$ as a group of automorphisms.
We show that if $\mathbf V_d(t)$ is twisted via a generator for
$\Z_4$, 
then the resulting
$\boxtimes_q$-module is isomorphic to $\mathbf V_d(t^{-1})$.
Consider the element of $\Z_4$ that has order 2.
If $\mathbf V_d(t)$ is twisted via this element, 
then the resulting
$\boxtimes_q$-module is isomorphic to $\mathbf V_d(t)$.
A corresponding isomorphism of $\boxtimes_q$-modules
is called an exchanger.
We describe how these exchangers
act on the 24 bases for $\mathbf V_d(t)$.
We also characterize the exchangers in various ways.
\medskip

\noindent  Near the end of the paper we discuss how
Leonard pairs of $q$-Racah type are related to
evaluation modules for $\boxtimes_q$.
Given a Leonard pair of $q$-Racah type,
we consider a certain basis for the
underlying vector space, called the compact basis,
with respect to which the matrices
representing the pair are each tridiagonal with attractive
entries. 
Using the Leonard pair 
we turn the underlying vector space into an evaluation module for
$\boxtimes_q$.
On this $\boxtimes_q$-module, each element of the Leonard pair
acts as a linear combination of two $\boxtimes_q$-generators;
the associated edges in the
tetrahedron are adjacent and singly oriented.
We show that the compact basis
diagonalizes a pair of $\boxtimes_q$-generators that correspond
to a doubly oriented edge of the tetrahedron.

\medskip
\noindent The paper is organized as follows.
In Section 2 we review some preliminaries and fix our notation.
In Section 3 we recall some facts about
$U_q(\mathfrak{sl}_2)$ that will be used throughout the paper.
In Sections 4--6 we describe the finite-dimensional
irreducible $\boxtimes_q$-modules of type 1.
In Sections 7, 8 we further describe these $\boxtimes_q$-modules,
bringing in the dual space and $\Z_4$-twisting.
In Section 9 we define an evaluation module for
$\boxtimes_q$ called $\mathbf V_d(t)$, and we
obtain several polynomial identities that hold on
this module.
In Section 10 we describe 24 bases for $\mathbf V_d(t)$.
In Section 11 we describe how the eight $\boxtimes_q$-generators
act on the 24 bases.
In Sections 12, 13 we describe the transition matrices between
certain pairs of bases among the 24.
In Section 14 we obtain some identities that involve $\mathbf V_d(t)$
and its dual space.
Section 15 is about exchangers.
In Sections 16, 17 we describe how the evaluation modules for
$\boxtimes_q$ are related to Leonard pairs of $q$-Racah type.
Appendices 18, 19 contain some matrix definitions and related
identities.

\section{Preliminaries}
Our conventions are as follows.
Throughout the paper $\F$ denotes an algebraically closed field.
An algebra is meant to be associative and have a 1.
A subalgebra has the same 1 as the parent algebra.
Recall the natural numbers $\N=\lbrace 0,1,2,\ldots \rbrace$
and integers $\Z = \lbrace 0, \pm 1, \pm 2,\ldots \rbrace$.
For the duration of this paragraph fix $d \in \N$.
Let  $\lbrace u_n\rbrace_{n=0}^d$ denote
a sequence.
We call $u_n$ the {\it $n$th  component} of the
sequence. We call $d$ the {\it diameter} of the sequence.
By the {\it inversion} of the sequence
 $\lbrace u_n\rbrace_{n=0}^d$ we mean the sequence
 $\lbrace u_{d-n}\rbrace_{n=0}^d$.
Let $V$ denote a vector space over $\F$ with dimension $d+1$.
 Let
${\rm End}(V)$ denote the $\F$-algebra consisting of
the $\F$-linear maps from $V$ to $V$. An element $A \in
{\rm End}(V)$ is called {\it diagonalizable} whenever
$V$ is spanned by the eigenspaces of $A$.
The map $A$ is called {\it multiplicity-free} whenever
$A$ is diagonalizable and each eigenspace of $A$ has
dimension 1.
Let ${\rm Mat}_{d+1}(\F)$ denote the $\F$-algebra consisting of
the $d+1$ by $d+1$ matrices that have all entries in $\F$.
We index the rows and columns by $0,1,\ldots, d$.
Let $\lbrace v_n\rbrace_{n=0}^d$
denote a basis for $V$. For $A \in {\rm End}(V)$
and $M\in 
{\rm Mat}_{d+1}(\F)$, we say that {\it $M$ represents $A$
with respect to 
 $\lbrace v_n\rbrace_{n=0}^d$} whenever
$Av_n = \sum_{i=0}^d M_{in}v_i$ for $0 \leq n \leq d$.
For $M \in 
{\rm Mat}_{d+1}(\F)$,
$M$ is called
{\it upper bidiagonal} whenever 
 each nonzero entry lies on the diagonal or the superdiagonal.
The matrix $M$ is called 
{\it lower bidiagonal} whenever the transpose $M^t$ is upper
bidiagonal.
The matrix $M$ is called
{\it tridiagonal} whenever each nonzero entry
lies on the diagonal, the subdiagonal, or the superdiagonal.
Assume that $M$ is tridiagonal. Then $M$ is said to be
{\it irreducible} whenever each entry on the subdiagonal is
nonzero and each entry on the superdiagonal is nonzero.

\begin{definition}
\label{def:LPdef}
\rm 
 \cite[Definition~1.1]{terLS}.
Let $V$ denote a vector space over $\F$
with finite positive dimension. By a {\it Leonard pair} on $V$,
we mean an ordered pair $A,B$ of elements in 
${\rm End}(V)$ that satisfy the following conditions:
\begin{enumerate}
\item[\rm (i)] there exists a basis for $V$ with respect to which
the matrix representing $A$ is diagonal and the matrix representing
$B$ is irreducible tridiagonal;
\item[\rm (ii)] there exists a basis for $V$ with respect to which
the matrix representing $B$ is diagonal and the matrix representing
$A$ is irreducible tridiagonal.
\end{enumerate}
The above Leonard pair is said to be {\it over $\F$}. 
We call $V$ the {\it underlying vector space}.
\end{definition}

\begin{definition} 
\label{def:isoLP}
\rm
Let $A,B$ denote a Leonard pair over $\F$.
Let $A',B'$ denote a Leonard pair over $\F$.
By an {\it isomorphism of Leonard pairs} from
$A,B$ to $A',B'$ we mean an $\F$-linear bijection
$\mu$ from the vector space underlying $A,B$ to the
vector space underlying $A',B'$ such that
$\mu A=A' \mu$ and 
$\mu B=B' \mu $.
The Leonard pairs
$A,B$ 
and  $A',B'$ are called {\it isomorphic} whenever
there exists an isomorphism of Leonard pairs from
$A,B$ to $A',B'$.
\end{definition}

\begin{lemma}
\label{lem:LPQR}
{\rm \cite[Corollary~5.5]{ter:LPQR}}.
Let $A,B$ denote a Leonard pair
on $V$, as in Definition
\ref{def:LPdef}.
Then
the algebra ${\rm End}(V)$ is generated by $A,B$.
\end{lemma}

\noindent We refer the reader to
 \cite{terLS,
ter24,
ter:LPQR,
vidter} for background information on Leonard pairs.

\section{The equitable presentation for
$U_q(\mathfrak{sl}_2)$}

\noindent In this section we recall
the quantum  enveloping algebra
$U_q(\mathfrak{sl}_2)$. 
For background information on
$U_q(\mathfrak{sl}_2)$,
 we refer the
reader to the books by Jantzen
\cite{jantzen} and Kassel \cite{kassel}.
We will work with the equitable presentation of
$U_q(\mathfrak{sl}_2)$,
which was introduced in \cite{equit}.
\medskip

\noindent 
Throughout the paper,
fix a nonzero $q \in \F$ 
that is not a root of unity.
For $n \in \Z$ define
\begin{eqnarray*}
\lbrack n \rbrack_q = \frac{q^n-q^{-n}}{q-q^{-1}}
\end{eqnarray*}
and for $n\geq 0$  define
\begin{eqnarray*}
\lbrack n \rbrack^!_q 
=
\lbrack n \rbrack_q 
\lbrack n-1 \rbrack_q 
\cdots 
\lbrack 2 \rbrack_q 
\lbrack 1 \rbrack_q.
\end{eqnarray*}
We interpret 
$\lbrack 0 \rbrack^!_q=1$.

\begin{definition}
\label{def:uqdef}
\rm
\cite[Theorem 2.1]{equit}.
For the $\F$-algebra 
$U_q(\mathfrak{sl}_2)$
the equitable presentation
has generators $x, y^{\pm 1}, z$ and relations
$yy^{-1} = 1, y^{-1}y=1$,
\begin{eqnarray}
\label{eq:uqdef}
\frac{qxy-q^{-1}yx}{q-q^{-1}} = 1,\qquad 
\frac{qyz-q^{-1}zy}{q-q^{-1}} = 1,\qquad 
\frac{qzx-q^{-1}xz}{q-q^{-1}} = 1.
\end{eqnarray}
We call
$x, y^{\pm 1}, z$ the {\it equitable generators} for
$U_q(\mathfrak{sl}_2)$.
\end{definition}

\noindent In the next three lemmas we comment on the relations
(\ref{eq:uqdef}).
\begin{lemma}
\label{lem:uv}
Let $u,v$ denote elements in any $\F$-algebra, such that
\begin{eqnarray*}
\frac{quv-q^{-1}vu}{q-q^{-1}}=1.
\end{eqnarray*}
Then
\begin{eqnarray}
q(1-uv) = q^{-1}(1-vu) = \frac{\lbrack u,v\rbrack}{q-q^{-1}},
\label{lem:commute}
\end{eqnarray}
where $\lbrack u,v\rbrack$ means $uv-vu$.
\end{lemma}
\noindent {\it Proof:}  
Routine.
\hfill $\Box$ \\

\begin{lemma}
\label{lem:uvtwo}
Let $u$, $v$ denote elements in any $\F$-algebra,
such that
\begin{eqnarray*}
\frac{quv-q^{-1}vu}{q-q^{-1}}=1.
\end{eqnarray*}
\begin{enumerate}
\item[\rm (i)]
Assume $u^{-1}$ exists.  Then
\begin{eqnarray}
\bigl \lbrack u^{-1}, \lbrack u,v \rbrack \bigr \rbrack
= (q-q^{-1})^2 (u^{-1}-v).
\label{eq:uuv}
\end{eqnarray}
\item[\rm (ii)]
Assume $v^{-1}$ exists.  Then
\begin{eqnarray}
\bigl \lbrack \lbrack u,v \rbrack, v^{-1} \bigr \rbrack
= (q-q^{-1})^2 (v^{-1}-u).
\label{eq:uvv}
\end{eqnarray}
\end{enumerate}
\end{lemma}
\noindent {\it Proof:}  
(i) Using Lemma
\ref{lem:uv},
\begin{eqnarray*}
\frac{\bigl \lbrack u^{-1}, \lbrack u,v \rbrack \bigr \rbrack}{q-q^{-1}}
= qu^{-1}(1-uv)- q^{-1}(1-vu)u^{-1}
= (q-q^{-1}) (u^{-1}-v).
\end{eqnarray*}
\noindent (ii) Similar to the proof of (i) above.
\hfill $\Box$ \\

\begin{lemma}
\label{lem:uvw} 
Let $u$, $v$, $w$ denote elements in any $\F$-algebra, such that
both
\begin{eqnarray*}
\frac{quv-q^{-1}vu}{q-q^{-1}}=1,
\qquad
\qquad 
\frac{qvw-q^{-1}wv}{q-q^{-1}}=1.
\end{eqnarray*}
Then  both
\begin{eqnarray}
\lbrack v,uw\rbrack = q(q-q^{-1})(u-w),
\qquad \qquad
\lbrack v,wu\rbrack= q^{-1}(q-q^{-1})(u-w).
\label{eq:uvw}
\end{eqnarray}
Moreover
\begin{eqnarray}
\bigl \lbrack v, \lbrack u, w\rbrack \bigr \rbrack =(q-q^{-1})^2(u-w).
\label{eq:lieuvw}
\end{eqnarray}
\end{lemma}
\noindent {\it Proof:}  
To obtain
(\ref{eq:uvw}), observe
\begin{eqnarray*}
\frac{\lbrack v,uw\rbrack}{q-q^{-1}}
= qu\frac{qvw-q^{-1}wv}{q-q^{-1}} - q
\frac{quv-q^{-1}vu}{q-q^{-1}} w 
= q(u-w),
\end{eqnarray*}
and
\begin{eqnarray*}
\frac{\lbrack v,wu\rbrack}{q-q^{-1}}
= q^{-1}\frac{qvw-q^{-1}wv}{q-q^{-1}}u - 
q^{-1}w\frac{quv-q^{-1}vu}{q-q^{-1}}  
= q^{-1}(u-w).
\end{eqnarray*}
We have obtained
(\ref{eq:uvw}), and
(\ref{eq:lieuvw}) follows.
\hfill $\Box$ \\


\noindent In the literature on
$U_q(\mathfrak{sl}_2)$ there is a certain central
element called the Casimir element
\cite[Section~2.7]{jantzen},
\cite[p.~130]{kassel}.
We now recall how the Casimir element looks
from the equitable point of view.

\begin{definition}
\label{def:nce}
\rm
\cite[Lemma~2.15]{uawe}.
Let $\Lambda$ denote the following element
in
$U_q(\mathfrak{sl}_2)$:
\begin{eqnarray}
\Lambda = 
qx + q^{-1} y + q z - q xyz.
\label{eq:cas}
\end{eqnarray}
We call $\Lambda$  the {\it (normalized) Casimir element}.
\end{definition}

\begin{note} \rm The element $\Lambda
(q-q^{-1})^{-2}$ is equal to the Casimir
element of
$U_q(\mathfrak{sl}_2)$ discussed in
\cite[Section~2.7]{jantzen}.
\end{note}


\begin{lemma} {\rm \cite[Lemma~2.7, Proposition~2.18]{jantzen}. } 
The elements $\lbrace \Lambda^n \rbrace_{n \in \N}$
form a basis for the 
 center of 
$U_q(\mathfrak{sl}_2)$.
\end{lemma}

\begin{lemma} 
{\rm \cite[Lemma~2.15]{uawe}. } 
The element $\Lambda$ is equal to
each of the following:
\begin{eqnarray*}
&&
qx + q^{-1} y + q z - q xyz, \qquad \qquad
q^{-1}x + qy + q^{-1} z - q^{-1} zyx, 
\\
&&
q y + q^{-1} z + q x - q yzx, \qquad \qquad
q^{-1}y + qz + q^{-1} x - q^{-1} xzy, 
\\
&&
q z + q^{-1} x + q y - q zxy, \qquad \qquad
q^{-1}z + qx + q^{-1} y - q^{-1} yxz.
\end{eqnarray*}
\end{lemma}

\noindent The 
 finite-dimensional irreducible
$U_q(\mathfrak{sl}_2)$-modules
are described in 
\cite[Section~2]{jantzen}.
We now recall how these modules
look from the equitable point of view
\cite{equit, fduqe}.

\begin{lemma} 
\label{lem:xyz}
{\rm \cite[Lemma~4.2]{equit}},
{\rm \cite[Theorem~2.6]{jantzen}}.
There exists a family of finite-dimensional irreducible
$U_q(\mathfrak{sl}_2)$-modules 
\begin{eqnarray}
\label{eq:Vde}
{\mathbf V_{d,\varepsilon}}
\qquad \qquad \varepsilon \in \lbrace 1,-1\rbrace,
\qquad \qquad d \in \N
\end{eqnarray} with the following property:
$\mathbf V_{d,\varepsilon}$ has a basis with respect to which 
the
matrices representing $x,y, z$ are
    \begin{eqnarray*}
     &&x: \quad \varepsilon \left(
   \begin{array}{ccccc}
  q^{-d} & q^d-q^{-d}   & & & {\bf 0}  \\
  & q^{2-d} & q^d-q^{2-d}    & &   \\
 && q^{4-d} & \ddots &
 \\
  & & & \ddots  & q^d-q^{d-2} \\
 {\bf 0}  & &  & &q^{d} 
 \end{array}
 \right),
\\
 \\
 &&y: \quad
       \varepsilon \; {\rm diag} (q^d, q^{d-2}, q^{d-4}, \ldots, q^{-d}),
 \\
 \\
 &&
 z: \quad 
 \varepsilon \left(
 \begin{array}{ccccc}
 q^{-d} &  & & & {\bf 0}  \\
 q^{-d}-q^{2-d} &  q^{2-d}&  & &\\
 & q^{-d}-q^{4-d} &  q^{4-d} &   & \\
 & & \ddots &  \ddots &   \\
 {\bf 0} & & & q^{-d}-q^d   & q^{d} 
 \end{array}
 \right).
 \end{eqnarray*}
 Every finite-dimensional irreducible 
 $U_q(\mathfrak{sl}_2)$-module is isomorphic to exactly one
 of the modules 
 {\rm (\ref{eq:Vde})}.
\end{lemma}

\begin{note}\rm 
For ${\rm Char} (\F)=2$  we interpret 
$\lbrace 1,-1\rbrace$ to have a single element.
\end{note}

\begin{note}\rm The dimension of
$\mathbf V_{d,\varepsilon}$ is $d+1$.
\end{note}

\begin{definition}\rm 
For $\mathbf V_{d,\varepsilon}$
the parameter $d$ is called the 
{\it diameter}.
The parameter  
$\varepsilon$ is called the
{\it type}.
We sometimes abbreviate $\mathbf V_d = \mathbf V_{d,1}$.
\end{definition}

\begin{note}\rm
For each of $x,y^{\pm 1},z$ the action on
$\mathbf V_{d,\varepsilon}$ 
is multiplicity-free with 
 eigenvalues
$\lbrace \varepsilon q^{d-2n}|0 \leq n \leq d  \rbrace$.
\end{note}

\begin{note}\rm In Lemma
\ref{lem:xyz} the matrix representing $x$ (resp. $z$) 
has constant
row sum $ \varepsilon q^d$ (resp.
 $\varepsilon q^{-d}$). This reflects the fact that
 $xv=\varepsilon q^d v$
 (resp. $zv=\varepsilon q^{-d} v$), where
 $v$ denotes the sum of the basis vectors.
\end{note}

\begin{note}\rm
\label{note:names}
In Appendix I we define some
matrices in 
${\rm Mat}_{d+1}(\F)$ called $E_q$, $K_q$, $Z$.
The displayed matrices from
 Lemma
\ref{lem:xyz} 
that represent 
$x,y,z$ for $\mathbf V_d$ are
$E_q$,
$K_q$, 
$ZE_{q^{-1}}Z$ respectively.
\end{note}

\begin{lemma}
\label{lem:nce}
{\rm \cite[Section~2.7]{jantzen}.}
The Casimir element
$\Lambda $ acts on $\mathbf V_d $
as $ (q^{d+1} + q^{-d-1})I$.
\end{lemma}

\begin{lemma}
\label{lem:rowbasis}
{\rm \cite[Lemma~9.8]{fduqe}}.
Pick $\xi \in \lbrace x,y,z\rbrace$. 
There exists a basis
$\lbrace v_n\rbrace_{n=0}^d$ 
for $\mathbf V_d$  such that
\begin{enumerate}
\item[\rm (i)]
$\xi v_n = q^{d-2n}v_n$ for
$0 \leq n \leq d$;
\item[\rm (ii)] $\sum_{n=0}^d v_n$ is a common eigenvector
for the two elements among
$x,y,z$ other than $\xi$.
\end{enumerate}
\end{lemma}

%

\begin{definition}
\label{def:rowbasisinv}
\rm
{\cite[Definition~9.5]{fduqe}}.
Pick $\xi \in \lbrace x,y,z\rbrace$.
By a 
{\it $\lbrack \xi \rbrack_{row}$-basis} for
$\mathbf V_d$ we mean a basis for $\mathbf V_d$ 
from Lemma
\ref{lem:rowbasis}.
\end{definition}

\begin{note}\rm
The basis for
$\mathbf V_d$ in Lemma
\ref{lem:xyz} is a $\lbrack y \rbrack_{row}$-basis.
\end{note}

\noindent We comment on the uniqueness of the bases in
Definition \ref{def:rowbasisinv}.

\begin{lemma} 
{\rm 
\cite[Lemma~9.12]{fduqe}}.
Pick $\xi \in \lbrace x,y,z\rbrace$ and let
$\lbrace v_n\rbrace_{n=0}^d$ denote a $\lbrack \xi \rbrack_{row}$-basis
for  $\mathbf V_d$. Let 
$\lbrace v'_n\rbrace_{n=0}^d$ denote any vectors in $\mathbf V_d$.
Then the following are equivalent:
\begin{enumerate}
\item[\rm (i)] the sequence
$\lbrace v'_n\rbrace_{n=0}^d$ is a 
$\lbrack \xi \rbrack_{row}$-basis
for  $\mathbf V_d$;
\item[\rm (ii)] there exists
$0 \not= \alpha \in \F$ such that
$v'_n = \alpha v_n$ for $0 \leq n \leq d$.
\end{enumerate}
\end{lemma}


\begin{definition}\rm
\label{def:rowbasisinv2}
\rm
{\cite[Definition~9.5]{fduqe}}.
Pick $\xi \in \lbrace x,y,z\rbrace$.
By 
a {\it $\lbrack \xi \rbrack^{inv}_{row}$-basis} for $\mathbf V_d$
we mean the inversion of a 
 $\lbrack \xi \rbrack_{row}$-basis for $\mathbf V_d$.
\end{definition}

\noindent 
In Definition
\ref{def:rowbasisinv}  and
Definition
\ref{def:rowbasisinv2}
we gave the following six bases for $\mathbf V_d$:
\begin{eqnarray}
&&
\lbrack x \rbrack_{row}, \qquad 
\lbrack y\rbrack_{row}, \qquad 
\lbrack z \rbrack_{row},
\label{eq:urow}
\\
&&
\lbrack x \rbrack^{inv}_{row}, \qquad 
\lbrack y\rbrack^{inv}_{row}, \qquad 
\lbrack z \rbrack^{inv}_{row}.
\label{eq:urowinv}
\end{eqnarray}

\begin{lemma} 
\label{lem:sixb}
{\rm \cite[Theorem~10.12]{fduqe}}.
Consider the elements $x$, $y$, $z$ of 
$U_q(\mathfrak{sl}_2)$. In the table below we display
the matrices in 
${\rm Mat}_{d+1}(\F)$
that represent these elements with respect to
the six bases 
{\rm (\ref{eq:urow}), 
(\ref{eq:urowinv})}
for $\mathbf V_d$.

\medskip

\centerline{
\begin{tabular}[t]{c |ccc}
 {\rm basis} &
  {\rm $x$} &
  {\rm  $y$} &
    {\rm  $z$}
  \\ \hline  \hline
   $\lbrack x \rbrack_{row}$
   &
    $K_q$
    &
     $ZE_{q^{-1}}Z$
&
 $E_q$
  \\
   $\lbrack x \rbrack^{inv}_{row}$
    &
    $K_{q^{-1}}$
      &
      $E_{q^{-1}}$
 &
 $ZE_qZ$
   \\
   \hline
 $\lbrack y \rbrack_{row}$
  &
   $E_q$
    &
     $K_q$
      &
       $ZE_{q^{-1}}Z$
        \\
    $\lbrack y \rbrack^{inv}_{row}$
&
 $ZE_qZ$
  &
    $K_{q^{-1}}$
    &
    $E_{q^{-1}}$
      \\
 \hline
 $\lbrack z \rbrack_{row}$
  &
   $ZE_{q^{-1}}Z$
    &
     $E_q$
      &
       $K_q$
        \\
      $\lbrack z \rbrack^{inv}_{row}$
 &
 $E_{q^{-1}}$
  &
   $ZE_qZ$
    &
    $K_{q^{-1}}$
      \\
           \end{tabular}}
        \medskip

\end{lemma}
\noindent For more background information on the
bases
{\rm (\ref{eq:urow}), 
(\ref{eq:urowinv})}
we refer the reader to
\cite{fduqe}.



\section{The $q$-tetrahedron algebra $\boxtimes_q$}

\noindent In this section we recall the $q$-tetrahedron algebra
$\boxtimes_q$ and review some of its properties.

\medskip
\noindent 
Let $\Z_4=\Z /4\Z$ denote the cyclic group of order 4.
\begin{definition} 
\label{def:qtet}
\rm 
\cite[Definition~6.1]{qtet}.
Let $\boxtimes_q$ denote the $\F$-algebra defined by generators
\begin{eqnarray}
\lbrace x_{ij}\;|\;i,j \in 
\Z_4,
\;\;j-i=1 \;\mbox{\rm {or}}\;
j-i=2\rbrace
\label{eq:gen}
\end{eqnarray}
and the following relations:
\begin{enumerate}
\item[\rm (i)] For $i,j \in 
\Z_4$ such that $j-i=2$,
\begin{eqnarray}
 x_{ij}x_{ji} =1.
\label{eq:tet1}
\end{eqnarray}
\item[\rm (ii)] For $i,j,k \in 
\Z_4$ such that $(j-i,k-j)$ is one of $(1,1)$, $(1,2)$, 
$(2,1)$,
\begin{eqnarray}
\label{eq:tet2}
\frac{qx_{ij}x_{jk} - q^{-1} x_{jk}x_{ij}}{q-q^{-1}}=1.
\end{eqnarray}
\item[\rm (iii)] For $i,j,k,\ell \in 
\Z_4$ such that $j-i=k-j=\ell-k=1$,
\begin{eqnarray}
\label{eq:tet3}
x^3_{ij}x_{k\ell} 
- 
\lbrack 3 \rbrack_q
x^2_{ij}x_{k\ell} x_{ij}
+ 
\lbrack 3 \rbrack_q
x_{ij}x_{k\ell} x^2_{ij}
-
x_{k\ell} x^3_{ij} = 0.
\end{eqnarray}
\end{enumerate}
We call 
$\boxtimes_q$ the {\it $q$-tetrahedron algebra}.
The elements
(\ref{eq:gen})
 are called the {\it standard generators} for $\boxtimes_q$.
\end{definition}


\noindent We have some comments.

\begin{note}\rm We find it illuminating to view
$\boxtimes_q$ as follows.
Identify $\mathbb Z_4$ with the vertex set of a tetrahedron.
View each standard generator $x_{ij}$  as an orientation
$i\to j$ of the edge 
in the tetrahedron that involves vertices $i$ and $j$.
\end{note}

\begin{lemma} 
\label{lem:rho}
There exists an 
automorphism
$\rho $ of $\boxtimes_q$ that sends each
standard generator $x_{ij}$ to
$x_{i+1,j+1}$. Moreover
$\rho^4 = 1$.
\end{lemma}

\begin{lemma} There exists an 
automorphism
$\sigma$ of  $\boxtimes_q$ that sends each
standard generator $x_{ij}$ to
$-x_{ij}$. We have $\sigma^2=1$ if ${\rm Char}(\F) \not=2$
and $\sigma = 1$ if 
 ${\rm Char}(\F) =2$.
\end{lemma}

\begin{lemma}
\label{lem:miki}
{\rm \cite[Proposition~4.3]{miki}.}
For $i \in \Z_4$ there exists an
$\F$-algebra homomorphism 
$\kappa_i: U_q(\mathfrak{sl}_2)\to \boxtimes_q$ that sends
\begin{eqnarray*}
x \mapsto x_{i+2,i+3},
\qquad 
y \mapsto x_{i+3,i+1},
\qquad
y^{-1} \mapsto x_{i+1,i+3},
\qquad 
z \mapsto x_{i+1,i+2}.
\end{eqnarray*}
This homomorphism is injective.
\end{lemma}

\noindent Recall the Casimir element $\Lambda $
of $U_q(\mathfrak{sl}_2)$, from
Definition \ref{def:nce}.

\begin{definition} 
\label{def:casi}
\rm For $i \in \Z_4$ let $\Upsilon_i$ denote
the image of $\Lambda$ under the injection $\kappa_i$ 
from Lemma
\ref{lem:miki}.
\end{definition}

\noindent The elements $\Upsilon_i$ from Definition
\ref{def:casi} are not central in $\boxtimes_q$. However
we do have the following.

\begin{lemma} 
For $i \in \Z_4$ the element $\Upsilon_i$ commutes with
each of
\begin{eqnarray*}
x_{i+2,i+3},\qquad
x_{i+3,i+1},\qquad
x_{i+1,i+3},\qquad
x_{i+1,i+2}.
\end{eqnarray*}
\end{lemma}
\noindent {\it Proof:}  
By Lemma
\ref{lem:miki} and since 
$\Lambda$ is central in
 $U_q(\mathfrak{sl}_2)$.
\hfill $\Box$ \\

\section{Comparing $\boxtimes_q$ and 
$\boxtimes_{q^{-1}}$}

\noindent In this section we compare the 
algebras
$\boxtimes_q$ and 
$\boxtimes_{q^{-1}}$. For both algebras
we use same notation
$x_{ij}$ for the standard generators.

\begin{lemma} There exists an $\F$-algebra isomorphism
$\vartheta: \boxtimes_q \to \boxtimes_{q^{-1}}$ that sends
\begin{eqnarray*}
&&x_{01} \mapsto x_{01}, \qquad
x_{12} \mapsto x_{30}, \qquad
x_{23} \mapsto x_{23}, \qquad
x_{30} \mapsto x_{12},
\\
&&x_{02} \mapsto x_{31}, \qquad
x_{13} \mapsto x_{20}, \qquad
x_{20} \mapsto x_{13}, \qquad
x_{31} \mapsto x_{02}.
\end{eqnarray*}
\end{lemma}
\noindent {\it Proof:}  
Routine.
\hfill $\Box$ \\

\noindent We recall the notion of antiisomorphism.
Given $\F$-algebras $\mathcal A$, $\mathcal B$
a map $\gamma: \mathcal A \to \mathcal B$ is called
an {\it antiisomorphism of $\F$-algebras} whenever
$\gamma$ is an isomorphism of $\F$-vector spaces
and $(ab)^\gamma=b^\gamma a^\gamma$ for all
$a,b \in \mathcal A$. An antiisomorphism can be interpreted
as follows. The $\F$-vector space $\mathcal B$ has
an $\F$-algebra structure ${\mathcal B}^{opp}$ such that
for all $a,b \in \mathcal B$ the product
$ab$ (in 
${\mathcal B}^{opp}$) is equal to
$ba$ (in 
${\mathcal B}$).
A map $\gamma: \mathcal A \to \mathcal B$ is an antiisomorphism
of $\F$-algebras if and only if 
$\gamma: \mathcal A \to {\mathcal B}^{opp}$ is an isomorphism
of $\F$-algebras.

\begin{proposition} 
\label{prop:antiiso}
There exists an $\F$-algebra
antiisomorphism
$\tau : \boxtimes_q \to \boxtimes_{q^{-1}}$ that sends
each standard generator $x_{ij}$ to $x_{ij}$.
\end{proposition}
\noindent {\it Proof:}  
In Definition
\ref{def:qtet} we gave a presentation for
$\boxtimes_q$ by generators and relations.
We now modify this presentation by adjusting the
relations as follows.
For each relation in the presentation,
replace $q$ by $q^{-1}$ and invert the
order of multiplication. In each case, the adjusted
relation coincides with the original one.
Now on one hand, the modified presentation
 is a presentation for 
$(\boxtimes_{q^{-1}})^{opp}$ by generators and relations.
On the other hand, the modified   
 presentation coincides with the original one.
Therefore
there exists an $\F$-algebra  isomorphism 
$\tau: \boxtimes_q \to (\boxtimes_{q^{-1}})^{opp}$ that sends
each standard generator $x_{ij}$ to $x_{ij}$.
The result follows by the sentence prior to the theorem statement.
\hfill $\Box$ \\

\section{Finite-dimensional irreducible $\boxtimes_q$-modules}

\noindent In this section we review some basic facts 
and notation concerning finite-dimensional irreducible
$\boxtimes_q$-modules. This material is summarized from
\cite{qtet}.


\medskip
\noindent 
Let $V$ denote a vector space over $\F$ with finite positive
dimension. Let $\lbrace s_n\rbrace_{n=0}^d$ denote a sequence
of positive integers whose sum is the dimension of $V$.
By a {\it decomposition of $V$ of shape
 $\lbrace s_n\rbrace_{n=0}^d$} 
 we mean a sequence $\lbrace V_n \rbrace_{n=0}^d$
 of subspaces for $V$  such that
 $V_n$ has dimension $s_n$ for $0 \leq n \leq d$
 and $V=\sum_{n=0}^d V_n$ (direct sum). 
For notational convenience define $V_{-1}=0$ and $V_{d+1}=0$.

\medskip
\noindent Now let $V$ denote a finite-dimensional irreducible
$\boxtimes_q$-module. By \cite[Theorem~12.3]{qtet} each standard generator
$x_{ij}$ of $\boxtimes_q$ is diagonalizable on $V$.
Also by 
\cite[Theorem~12.3]{qtet}
there exists  $d \in \N$ and 
$\varepsilon \in \lbrace 1,-1\rbrace$ such that
for each $x_{ij}$ the set of distinct eigenvalues on 
$V$ is $\lbrace \varepsilon q^{d-2n} | 0 \leq n \leq d\rbrace$.
We call $d$ the {\it diameter} of $V$.
We call $\varepsilon$
the {\it type} of $V$.
Replacing each $x_{ij}$ by $\varepsilon x_{ij}$
the type becomes 1.
So without loss of generality we may
assume that $V$ has type 1.
From now on we adopt this assumption.
For distinct $i,j \in \Z_4$ 
we now define a decomposition of
$V$ called 
$\lbrack i,j\rbrack$. 
First assume
$j-i=1$ or $j-i=2$, so that $x_{ij}$ exists.
The decomposition
$\lbrack i,j\rbrack$
has
diameter $d$,
and for $0 \leq n \leq d$ the $n$th component of
$\lbrack i,j\rbrack$ is the eigenspace of $x_{ij}$ 
with eigenvalue $q^{d-2n}$. 
Next assume $j-i=3$. In this case the decomposition
$\lbrack i,j\rbrack$ is defined as
 the inversion of
$\lbrack j,i\rbrack$.
By the construction and Definition
\ref{def:qtet}(i), 
for
 distinct $i,j \in \Z_4$
the decomposition 
$\lbrack i,j\rbrack$ is the inversion of
$\lbrack j,i\rbrack$.
By \cite[Proposition~13.3]{qtet}, 
for 
 distinct $i,j \in \Z_4$
the shape of 
$\lbrack i,j\rbrack$ is independent of
the choice of $i,j$.  Denote this shape by
$\lbrace \rho_n \rbrace_{n=0}^d$ and
note that
$\rho_n = \rho_{d-n}$ for $0 \leq n \leq d$.
By the {\it shape of $V$} we mean the sequence
$\lbrace \rho_n \rbrace_{n=0}^d$
\cite[Definition~13.4]{qtet}. 

\medskip
\noindent
One feature of the shape
$\lbrace \rho_n \rbrace_{n=0}^d$ is that
$\rho_{n-1}\leq \rho_n$ for
$1 \leq n \leq d/2$. This feature is obtained as follows.  
 Pick
 $i \in \Z_4$ and consider the homomorphism
$\kappa_i :
U_q(\mathfrak{sl}_2) \to \boxtimes_q$
from Lemma
\ref{lem:miki}. Using
 $\kappa_i$ we
pull back the
$\boxtimes_q$-module structure on $V$ to
obtain a
$U_q(\mathfrak{sl}_2)$-module structure on $V$.
The 
$U_q(\mathfrak{sl}_2)$-module $V$ is completely 
reducible; this means that
$V$
 is 
a direct sum of irreducible 
$U_q(\mathfrak{sl}_2)$-submodules 
\cite[Theorem 2.9]{jantzen}.
For this sum consider the summands.
Each summand has type 1.
For each summand the diameter is at most $d$
and has the same parity as $d$.
For $0 \leq n \leq d/2$ the following coincide:
(i) 
the multiplicity with which the 
$U_q(\mathfrak{sl}_2)$-module
$\mathbf V_{d-2n}$ appears as a summand;
(ii) the integer $\rho_n - \rho_{n-1}$, where $\rho_{-1}=0$.
Therefore $\rho_{n-1}\leq \rho_n$ for
$1 \leq n \leq d/2$.
We just mentioned some multiplicities.
These multiplicities
are independent of the $i \in \Z_4$ that we initially picked.
Therefore, 
up to isomorphism of 
$U_q(\mathfrak{sl}_2)$-modules
the $U_q(\mathfrak{sl}_2)$-module $V$ is independent of
$i$.

\medskip
\noindent Returning to the $\boxtimes_q$-module $V$,
for distinct $i,j\in \Z_4$
and each standard generator $x_{rs}$  we now describe the action of
$x_{rs} $ on the decomposition 
$\lbrack i,j\rbrack$ of $V$.
Denote this decomposition by
$\lbrace V_n\rbrace_{n=0}^d$. 
For the case $s-r=1$,
by \cite[Theorem~14.1]{qtet}
the
action of $x_{rs}$ on $V_n$ is described in the table
below:
\bigskip

\centerline{
\begin{tabular}[t]{c|c}
    {\rm decomposition} & {\rm action of $x_{r,r+1}$ on $V_n$}
    \\ \hline  \hline
    $\lbrack r,r+1\rbrack$ & $(x_{r,r+1}-q^{d-2n}I)V_n=0$
     \\
    $\lbrack r+1,r\rbrack$ & $(x_{r,r+1}-q^{2n-d}I)V_n=0$
    \\
    \hline
     $\lbrack r+1,r+2\rbrack$ &
       $(x_{r,r+1}-q^{2n-d}I)V_n \subseteq V_{n-1}$   \\
     $\lbrack r+2,r+1\rbrack$ &
       $(x_{r,r+1}-q^{d-2n}I)V_n \subseteq V_{n+1}$   \\
    \hline 
     $\lbrack r+2,r+3\rbrack$ &
      $x_{r,r+1}V_n\subseteq V_{n-1}+V_n+V_{n+1}$  \\
     $\lbrack r+3,r+2\rbrack$ &
      $x_{r,r+1}V_n\subseteq V_{n-1}+V_n+V_{n+1}$  \\
    \hline 
     $\lbrack r+3,r\rbrack$ &
          $(x_{r,r+1}-q^{2n-d}I)V_n\subseteq V_{n+1}$
           \\
     $\lbrack r,r+3\rbrack$ &
          $(x_{r,r+1}-q^{d-2n}I)V_n\subseteq V_{n-1}$
           \\
    \hline 
    $\lbrack r,r+2\rbrack$ &
             $
  (x_{r,r+1}-q^{d-2n}I)V_n\subseteq V_{n-1}$
    \\
    $\lbrack r+2,r\rbrack$ &
             $
  (x_{r,r+1}-q^{2n-d}I)V_n\subseteq V_{n+1}$
    \\
   \hline
    $ \lbrack r+1,r+3\rbrack $ &
        $
    (x_{r,r+1}-q^{2n-d}I)V_n\subseteq V_{n-1}$
     \\
    $ \lbrack r+3,r+1\rbrack $ &
        $
    (x_{r,r+1}-q^{d-2n}I)V_n\subseteq V_{n+1}$
      \end{tabular}}
   \bigskip

\noindent For the case $s-r=2$,
by \cite[Theorem~14.2]{qtet}
the action of $x_{rs}$ on $V_n$ is
described in the table below:
\bigskip

\centerline{
\begin{tabular}[t]{c|c}
    {\rm decomposition} & {\rm action of $x_{r,r+2}$ on $V_n$}
   \\ \hline  \hline
   $\lbrack r,r+1\rbrack$ &
$(x_{r,r+2}-q^{d-2n}I)V_n \subseteq V_0+\cdots+V_{n-1}$
       \\
   $\lbrack r+1,r\rbrack$ &
$(x_{r,r+2}-q^{2n-d}I)V_n \subseteq V_{n+1}+\cdots+V_d$
       \\
\hline
    $\lbrack r+1,r+2\rbrack$ &
   $(x_{r,r+2}-q^{d-2n}I)V_n \subseteq V_{n+1}+\cdots+V_{d}$     \\
    $\lbrack r+2,r+1\rbrack$ &
   $(x_{r,r+2}-q^{2n-d}I)V_n \subseteq V_0+\cdots+V_{n-1}$     \\
   \hline 
    $\lbrack r+2,r+3\rbrack$ &
      $(x_{r,r+2}-q^{2n-d}I)V_n\subseteq V_{n-1}$
       \\
    $\lbrack r+3,r+2\rbrack$ &
      $(x_{r,r+2}-q^{d-2n}I)V_n\subseteq V_{n+1}$
       \\
 \hline 
  $\lbrack r+3,r\rbrack$ &
$(x_{r,r+2}-q^{2n-d}I)V_n\subseteq V_{n+1}$
    \\
  $\lbrack r,r+3\rbrack$ &
$(x_{r,r+2}-q^{d-2n}I)V_n\subseteq V_{n-1}$
    \\
 \hline 
  $\lbrack r,r+2\rbrack$ &
    $
   (x_{r,r+2}-q^{d-2n}I)V_n=0$
    \\
  $\lbrack r+2,r\rbrack$ &
    $
   (x_{r,r+2}-q^{2n-d}I)V_n=0$
    \\
\hline 
 $ \lbrack r+1,r+3\rbrack $ &
   $
  x_{r,r+2}V_n\subseteq V_{n-1}+\cdots+V_d$
  \\
 $ \lbrack r+3,r+1\rbrack $ &
   $
  x_{r,r+2}V_n\subseteq V_0+\cdots+V_{n+1}$
 \end{tabular}}

\bigskip

\noindent We recall the notion of a flag. For the moment let
$V$ denote a vector space over $\F$ with finite positive
dimension and let 
$\lbrace s_n\rbrace_{n=0}^d$ denote a sequence
of positive integers whose sum is the dimension of $V$.
By a {\it flag on $V$ of shape
$\lbrace s_n\rbrace_{n=0}^d$}
we mean a sequence $\lbrace U_n\rbrace_{n=0}^d$ of subspaces for
$V$ such that $U_{n-1} \subseteq U_n$ for $1 \leq n \leq d$
and $U_n$ has dimension $s_0+\cdots + s_n$ for
$0 \leq n \leq d$. Observe that $U_d=V$. The following construction
yields a flag on $V$. Let $\lbrace V_n\rbrace_{n=0}^d$ denote
a decomposition of $V$ of shape $\lbrace s_n\rbrace_{n=0}^d$.
Define
$U_n= V_0+\cdots + V_n$ for $0 \leq n \leq d$. Then the sequence
$\lbrace U_n\rbrace_{n=0}^d$ is a flag on $V$ of shape
$\lbrace s_n\rbrace_{n=0}^d$. This flag is said to be
{\it induced} by the decomposition
 $\lbrace V_n\rbrace_{n=0}^d$.
We now recall what it means for two flags to be opposite.
Suppose we are given two flags on $V$ with the same diameter:
$\lbrace U_n\rbrace_{n=0}^d$ and
$\lbrace U'_n\rbrace_{n=0}^d$.
These flags are called {\it opposite} whenever 
there exists a decomposition $\lbrace V_n\rbrace_{n=0}^d$
of $V$ that induces 
$\lbrace U_n\rbrace_{n=0}^d$ and
whose inversion induces 
$\lbrace U'_n\rbrace_{n=0}^d$.
In this case $V_n = U_n \cap U'_{d-n}$ for $0 \leq n \leq d$, 
and also
$U_r \cap U'_s = 0$ if $r+s < d$ $(0 \leq r,s\leq d)$
\cite[p.~846]{ter24}.

\medskip
\noindent We now return our attention to $\boxtimes_q$-modules.
Let $V$ denote a finite-dimensional irreducible $\boxtimes_q$-module
of type 1 and diameter $d$. By 
\cite[Theorem 16.1]{qtet} there
exists a collection of 
flags on $V$, denoted
$\lbrack i \rbrack$,
$i \in \Z_4$, such that for distinct
$i,j \in \Z_4$ the decomposition 
$\lbrack i,j\rbrack$  induces
$\lbrack i \rbrack$.
By construction, the shape of the flag 
$\lbrack i \rbrack$ coincides with the shape of $V$.
By construction, the flags 
$\lbrack i \rbrack$,
$i \in \Z_4$ are mutually opposite.
Also by construction, for distinct 
$i,j \in \Z_4$ and $0 \leq n \leq d$
the $n$th component of 
$\lbrack i,j\rbrack$
is the intersection of the following two sets:
\begin{enumerate}
\item[\rm (i)]
component $n$ of flag $\lbrack i \rbrack$;
\item[\rm (ii)]
component $d-n$ of flag $\lbrack j \rbrack$.
\end{enumerate}

\noindent We comment on a very special case.
Up to isomorphism there exists a unique $\boxtimes_q$-module
of dimension one, on which each standard generator acts as the
identity. This $\boxtimes_q$-module is said to be {\it trivial}.
For a $\boxtimes_q$-module $V$ the following are equivalent:
(i) $V$ is trivial;
(ii) 
$V$ is finite-dimensional and irreducible,
with type 1 and diameter $0$.

\section{ The dual space}
Throughout this section  $V$ denotes a finite-dimensional irreducible
$\boxtimes_q$-module of type 1 and diameter $d$.
By definition, the dual space $V^*$ is the vector space over $\F$
consisting of the $\F$-linear maps from $V$ to $\F$. The vector spaces
$V$ and $V^*$ have the same dimension. 
In this section we will 
turn $V^*$
into 
a $\boxtimes_{q^{-1}}$-module, and discuss
how this 
module 
is related to the original $\boxtimes_q$-module $V$.

\begin{definition}
\label{def:bil}
\rm
Define a bilinear form $(\,,\,): V \times V^*\to \F$
such that $(u,f)=f(u)$ for all $u \in V$ and $f \in V^*$.
The form 
$(\,,\,)$ is nondegenerate.
\end{definition}
\noindent
Vectors $u \in V$ and $v \in V^*$ are called {\it orthogonal}
whenever $(u,v)=0$.

\medskip
\noindent We recall the adjoint map 
\cite[p.~227]{roman}. Let
$A \in
{\rm End}(V)$.
The {\it adjoint} of $A$, denoted $A^{adj}$,
is the unique element of 
${\rm End}(V^*)$ such that $(Au,v)=(u,A^{adj}v)$ for all 
$u \in V$ and $v \in V^*$. The adjoint map
${\rm End}(V)\to 
{\rm End}(V^*)$,
$A \mapsto A^{adj}$
is an antiisomorphism of $\F$-algebras.

\medskip
\noindent Recall the antiisomorphism $\tau:
 \boxtimes_q\to
 \boxtimes_{q^{-1}}$
from Proposition
\ref{prop:antiiso}.

\begin{proposition} 
\label{prop:adjointaction}
There exists a unique 
$\boxtimes_{q^{-1}}$-module structure on $V^*$ such that
\begin{eqnarray}
(\zeta u,v) = (u,\zeta^\tau v) 
\qquad u \in V, \quad v \in V^*, \quad \zeta \in \boxtimes_q.
\label{eq:requirement}
\end{eqnarray}
\end{proposition} 
\noindent {\it Proof:} 
The action of $\boxtimes_q$ on $V$ induces
 an $\F$-algebra homomorphism
 $\boxtimes_q\to 
{\rm End}(V)$. Call this homomorphism $\psi$.
The composition
\[
\begin{CD}
 \boxtimes_{q^{-1}}  @>>\tau^{-1} >  
 \boxtimes_q  @>>\psi> {\rm End}(V) 
                 @>>adj > {\rm End}(V^*)
                  \end{CD} 
              \]
is an $\F$-algebra homomorphism.
This homomorphism gives 
   $V^*$ a
$\boxtimes_{q^{-1}}$-module structure.
By construction 
the $\boxtimes_{q^{-1}}$-module $V^*$
satisfies  (\ref{eq:requirement}).
We have shown that the desired
$\boxtimes_{q^{-1}}$-module structure exists.
One routinely checks that this structure is unique.
\hfill $\Box$ \\


\begin{proposition} 
\label{prop:tetqtetqi}
For all $\zeta \in 
 \boxtimes_{q}$,
$\zeta^\tau$ acts on $V^*$
as the adjoint of the action of $\zeta $ on
$V$.
\end{proposition}
\noindent {\it Proof:}
By 
(\ref{eq:requirement}) and the definition 
of adjoint from above Proposition
\ref{prop:adjointaction}.
\hfill $\Box$ \\

\begin{proposition} 
\label{prop:howuquqirel}
For each standard generator $x_{ij}$,
\begin{eqnarray*}
&&(x_{ij}u,v) = (u,x_{ij} v)
\qquad \qquad
u \in V, \qquad v \in V^*.
\end{eqnarray*}
\end{proposition}
\noindent {\it Proof:}
Evaluate (\ref{eq:requirement}) 
using
Proposition
\ref{prop:antiiso}.
\hfill $\Box$ \\

\noindent
Given a subspace $W$ of $V$ (resp. $V^*$)
let $W^\perp$ denote the set of vectors in $V^*$ (resp. $V$)
that are orthogonal to everything
in $W$. The space $W^\perp$ is called the {\it orthogonal complement} of
$W$. 
For $W, W^\perp$ the sum of the dimensions
is equal to the common dimension of $V, V^*$.
Note that $(W^\perp)^\perp = W$.

\begin{lemma} 
\label{lem:wperp}
For a subspace $W \subseteq V$ and
 $\zeta \in 
 \boxtimes_q$,
$W$ is $\zeta$-invariant if and only if
 $W^\perp$ is $\zeta^\tau$-invariant.
\end{lemma}
\noindent {\it Proof:}
Use 
(\ref{eq:requirement}).
\hfill $\Box$ \\

\begin{lemma}
\label{lem:dualirred}
The 
 $\boxtimes_{q^{-1}}$-module $V^*$
is irreducible.
\end{lemma}
\noindent {\it Proof:}
Let $W$ denote a 
 $\boxtimes_{q^{-1}}$-submodule of $V^*$.
We show that $W=0$ or $W=V^*$. 
 Consider the orthogonal complement $W^\perp\subseteq V$.
By Lemma
\ref{lem:wperp} 
$W^\perp$
is a 
 $\boxtimes_q$-submodule of $V$.
 The 
 $\boxtimes_q$-module $V$ is irreducible
 so
$W^\perp=V$ or  
$W^\perp=0$.
It follows that
$W=0$ or $W=V^*$. 
\hfill $\Box$ \\

\begin{lemma} 
\label{lem:min}
For $\zeta \in 
 \boxtimes_q$ the following coincide:
\begin{enumerate}
\item[\rm (i)] the minimal polynomial for the action of $\zeta$ on  $V$;
\item[\rm (ii)] the minimal polynomial for the action of $\zeta^\tau$
on  $V^*$.
\end{enumerate}
\end{lemma}
\noindent {\it Proof:}
Use 
(\ref{eq:requirement}).
\hfill $\Box$ \\

\begin{proposition} 
\label{lem:vvs}
The $\boxtimes_{q^{-1}}$-module $V^*$ has type 1 and diameter $d$.
\end{proposition} 
\noindent {\it Proof:} 
By assumption the $\boxtimes_q$-module $V$
has type 1 and diameter $d$.
Therefore each $x_{ij}$ 
is diagonalizable on $V$ with eigenvalues
$\lbrace q^{d-2n}| 0 \leq n \leq d\rbrace $.
Now by
Lemma
\ref{lem:min} and since $x^\tau_{ij} = x_{ij}$,
the element $x_{ij}$ is diagonalizable on
$V^*$ with eigenvalues
$\lbrace q^{d-2n}| 0 \leq n \leq d\rbrace$.
The result follows.
\hfill $\Box$ \\

\medskip
\noindent 
Recall that $V$ is an irreducible $\boxtimes_q$-module
of type 1 and diameter $d$.
We described $V$ in Section 6.
In view of Lemma
\ref{lem:dualirred} and
Proposition
\ref{lem:vvs}, this description  (with $q$ replaced by $q^{-1}$)
applies to  $V^*$.
In particular, for distinct
$i,j\in \Z_4$
we may speak of  the
decomposition $\lbrack i,j\rbrack$ of $V^*$,
and for $i \in \Z_4$ we may speak of the
flag $\lbrack i\rbrack $ on $V^*$.

\medskip
\noindent 
Suppose we are given a decomposition of $V$
and a decomposition of $V^*$. These decompositions
are said to be {\it dual} whenever (i) they have the
same diameter $\delta$; (ii) for distinct $0 \leq i,j\leq \delta$
the $i$th component of the one is orthogonal to the 
$j$th component of the other.
Each decomposition of $V$ (resp. $V^*$) is dual to
a unique decomposition of $V^*$ (resp. $V$).
Dual decompositions have the same shape.

\begin{proposition}
\label{prop:dualdec}
For distinct $i,j\in \Z_4$ the decomposition
$\lbrack i,j\rbrack$ of $V$ is dual to the decomposition
$\lbrack j,i\rbrack$ of $V^*$.
\end{proposition}
\noindent {\it Proof:} 
First assume
$j-i=1$ or $j-i=2$, so that $x_{ij}$ exists.
Pick distinct integers $r,s$ $(0 \leq r,s\leq d)$.
Let $u$ denote a vector in component $r$
of the decomposition 
$\lbrack i,j\rbrack$ of $V$, and let
 $v$ denote a vector in component $s$
of the decomposition 
$\lbrack j,i \rbrack$ of $V^*$.
We show that $(u,v)=0$.
By Proposition
\ref{prop:howuquqirel}
$(x_{ij}u,v)=(u,x_{ij}v)$. By
construction
$x_{ij}u=q^{d-2r}u$ and $x_{ij}v=q^{d-2s}v$.
Observe that $q^{d-2r}\not=q^{d-2s}$  since $q$ is 
not a root of unity. By these comments
$(u,v)=0$.
Next assume
that $j-i=3$.
By the first part of this proof,
the decomposition $\lbrack j,i\rbrack$ of $V$ is dual to the decomposition
$\lbrack i,j\rbrack$ of $V^*$. Invert both decompositions to find that
the decomposition $\lbrack i,j\rbrack$ of $V$ is dual to the decomposition
$\lbrack j,i\rbrack$ of $V^*$.
\hfill $\Box$ \\

\begin{proposition} 
\label{prop:sameshape}
The 
$\boxtimes_q$-module $V$ and the
$\boxtimes_{q^{-1}}$-module $V^*$ have the same shape.
\end{proposition}
\noindent {\it Proof:} 
Pick distinct $i,j \in \Z_4$.
By definition, the shape of the $\boxtimes_q$-module
$V$ is equal to the shape of
the decomposition $\lbrack i,j\rbrack$ of $V$.
Similarly the shape of the $\boxtimes_{q^{-1}}$-module
$V^*$ is equal to the shape of
the decomposition $\lbrack j,i\rbrack$ of $V^*$.
All these  shapes are the same, by
 Proposition
\ref{prop:dualdec} and the sentence prior to
it. 
\hfill $\Box$ \\

\begin{proposition}
\label{prop:flagdual}
For $i \in \Z_4$ and $0 \leq n \leq d-1$
the following are orthogonal complements
with respect to the bilinear form $(\,,\,)$:
\begin{enumerate}
\item[\rm (i)] component $n$ of the flag $\lbrack i \rbrack$ on $V$;
\item[\rm (ii)] component $d-n-1$ of the flag
$\lbrack i \rbrack$ on $V^*$.
\end{enumerate}
\end{proposition}
\noindent {\it Proof:} 
Pick $j \in \Z_4$ with $j\not=i$.
By construction, the decomposition $\lbrack i,j\rbrack $ of $V$ 
(resp. $V^*$) induces the
flag $\lbrack i\rbrack $ of
$V$ 
(resp. $V^*$).
The decomposition $\lbrack i,j\rbrack $ of $V^*$ 
is the inversion of 
the decomposition $\lbrack j,i\rbrack $ of $V^*$.
The result follows from these comments and
Proposition
\ref{prop:dualdec}.
\hfill $\Box$ \\

\section{Twisting}

\noindent Recall the automorphism $\rho$ of $\boxtimes_q$
from Lemma
\ref{lem:rho}.
In this section we investigate what happens
when we twist a $\boxtimes_q$-module via
$\rho$.

\begin{definition}
\label{def:twist}
\rm
Let $V$ denote a $\boxtimes_q$-module and let $\pi$ denote
an  automorphism of $\boxtimes_q$. Then
there exists a $\boxtimes_q$-module structure on 
$V$, called {\it $V$ twisted via $\pi$}, that behaves as
follows: for all $\zeta \in \boxtimes_q$ and
$v \in V$, the vector $\zeta v$ computed in
$V$ twisted via $\pi$ coincides with the vector
$\pi^{-1}(\zeta) v$ computed in the original $\boxtimes_q$-module
$V$. Sometimes we abbreviate 
${}^\pi V$ for $V$ twisted via $\pi$.
\end{definition}

\begin{lemma} 
\label{lem:rhot}
Referring to Definition
\ref{def:twist}, the $\boxtimes_q$-module $V$ is irreducible
if and only if the 
 $\boxtimes_q$-module 
${}^\pi V$ is irreducible.
\end{lemma}
\noindent {\it Proof:} Immediate from Definition
\ref{def:twist}.
\hfill $\Box$ \\

\noindent Referring to Definition 
\ref{def:twist}, we now consider the case $\pi =\rho$.

\begin{lemma} 
\label{lem:rhoxijtwist}
Let $V$ denote a $\boxtimes_q$-module.
For each standard generator $x_{ij}$ of $\boxtimes_q$,
the following are the same:
\begin{enumerate}
\item[\rm (i)]
the action of
$x_{ij}$ on the $\boxtimes_q$-module $V$;
\item[\rm (ii)]
the action of
$x_{i+1,j+1}$ on the $\boxtimes_q$-module 
${}^\rho V$.
\end{enumerate}
\end{lemma}
\noindent {\it Proof:}
By
Lemma \ref{lem:rho} and
Definition
\ref{def:twist}.
\hfill $\Box$ \\

\begin{lemma}
\label{lem:typediam}
Let $V$ denote a finite-dimensional
irreducible $\boxtimes_q$-module
of type 1 and diameter $d$. Then the $\boxtimes_q$-module
${}^\rho V$ is irreducible, with type 1 and diameter $d$.
\end{lemma}
\noindent {\it Proof:}
By Lemma
\ref{lem:rhot} and 
Lemma
\ref{lem:rhoxijtwist},
along with the meaning of
type and diameter.
\hfill $\Box$ \\

\noindent For the following three lemmas the proof is
routine and left to the reader.

\begin{lemma} 
\label{lem:twistdec}
Let $V$ denote a finite-dimensional
irreducible $\boxtimes_q$-module of type 1. Then
for distinct $i,j \in \Z_4$ the following coincide:
\begin{enumerate}
\item[\rm (i)] the decomposition $\lbrack i,j\rbrack$ of $V$;
\item[\rm (ii)] the decomposition $\lbrack i+1,j+1\rbrack$ of
${}^\rho V$.
\end{enumerate}
\end{lemma}

\begin{lemma}
\label{lem:twistflag}
Let $V$ denote a finite-dimensional
irreducible $\boxtimes_q$-module of type 1. Then
for $i\in \Z_4$ the following coincide:
\begin{enumerate}
\item[\rm (i)] the flag $\lbrack i\rbrack$ on $V$;
\item[\rm (ii)] the flag $\lbrack i+1\rbrack$ on
${}^\rho V$.
\end{enumerate}
\end{lemma}

\begin{lemma}
\label{lem:twistshape}
Let $V$ denote a finite-dimensional
irreducible $\boxtimes_q$-module
of type 1. 
Then the following coincide:
\begin{enumerate}
\item[\rm (i)] the shape of $V$;
\item[\rm (ii)] the shape of
${}^\rho V$.
\end{enumerate}
\end{lemma}

\section{Evaluation modules}

\noindent We have been discussing the finite-dimensional
irreducible $\boxtimes_q$-modules. We now restrict our
attention to a special case, called an evaluation module.

\begin{definition}\rm By an {\it evaluation module}
for $\boxtimes_q$ we mean a nontrivial, finite-dimensional,
irreducible $\boxtimes_q$-module, of type 1 
and shape $(1,1,\ldots, 1)$.
\end{definition}

\noindent Let $V$ denote an evaluation module for
$\boxtimes_q$. Since $V$ is nontrivial
the diameter $d\geq 1$. Each standard generator 
$x_{ij}$ is multiplicity-free on $V$, with
eigenvalues $\lbrace q^{d-2n}| 0 \leq n \leq d\rbrace$.
In Section 6, for each $i \in \Z_4$ we
used the homomorphism
$\kappa_i : 
U_q(\mathfrak{sl}_2) \to \boxtimes_q$
to turn
$V$ into a
$U_q(\mathfrak{sl}_2)$-module.
Each of these
$U_q(\mathfrak{sl}_2)$-modules is isomorphic to
 $\mathbf V_d$. In order to
 recover the $\boxtimes_q$-module
 $V$ from $\mathbf V_d$, we add extra structure involving
 a parameter $t$. This is done as follows.

\begin{proposition} 
\label{prop:class1}
Pick an integer $d\geq 1$ and a nonzero
$t \in \F$ that is not among
$\lbrace q^{d-2n+1}\rbrace_{n=1}^d$.
Then there exists an evaluation module $\mathbf V_d(t)$ for $\boxtimes_q$
with the following property: $\mathbf V_d(t)$ has a basis with respect to
which the standard generators $x_{ij}$ are represented by the
following matrices in
${\rm Mat}_{d+1}(\F)$:

\medskip
\centerline{
\begin{tabular}[t]{c|cccc}
 {\rm generator} & 
 $x_{01}$ 
  &
  $x_{12}$ 
  &
  $x_{23}$ 
  &
  $x_{30}$ 
\\
\hline
{\rm representing matrix} 
&
$ZS_{q^{-1}}(t^{-1})Z$
&
$E_q$
&
$K_q$
&
$Z G_{q^{-1}}(t)Z$
     \end{tabular}}

\bigskip

\centerline{
\begin{tabular}[t]{c|cccc}
{\rm generator} & 
  $x_{02}$ 
  &
  $x_{13}$ 
  &
  $x_{20}$ 
  &
  $x_{31}$ 
\\
\hline
{\rm representing matrix}
&
$L_q(t)$ 
&
$(ZE_{q^{-1}}Z)^{-1}$
&
$(L_q(t))^{-1}$
&
$ZE_{q^{-1}}Z$
     \end{tabular}}
\medskip

\noindent
(The above matrices are defined in Appendix I).
\end{proposition}
\noindent {\it Proof:} It is routine (but tedious) to
check that the above matrices satisfy the defining
relations for $\boxtimes_q$. This makes $\mathbf V_d(t)$
a $\boxtimes_q$-module. The
$\boxtimes_q$-module $\mathbf V_d(t)$ is irreducible
by Note
\ref{note:names}  and since
the matrices $E_q$, $K_q$, $ZE_{q^{-1}}Z$ are 
included in the above tables.
By construction the $\boxtimes_q$-module 
$\mathbf V_d(t)$ is nontrivial, with type $1$ and shape
$(1,1,\ldots, 1)$.
\hfill $\Box$ \\

\noindent We have a comment.
\begin{lemma} 
\label{def:evalcom}
On the evaluation module
$\mathbf V_d(t)$ from 
Proposition
\ref{prop:class1}, 
the actions of the
 standard generators
{\rm (\ref{eq:gen})}
are mutually distinct.
\end{lemma}
\noindent {\it Proof:} The eight
representing matrices
in
the tables of Proposition
\ref{prop:class1} are mutually distinct.
This is checked using the definitions in
Appendix I.
\hfill $\Box$ \\

\noindent 
Consider the $\boxtimes_q$-module $\mathbf V_d(t)$ 
from Proposition
\ref{prop:class1}. This module
has dimension $d+1$.
Now let $V$ denote any evaluation module for $\boxtimes_q$
that has dimension $d+1$. 
Shortly we will show that the $\boxtimes_q$-module $V$
is isomorphic to 
$\mathbf V_d(t)$ for a unique $t$.

\begin{lemma}
\label{lem:eval}
On the  $\boxtimes_q$-module $\mathbf V_d(t)$,
\begin{eqnarray}
t(x_{01}-x_{23}) =
\frac{\lbrack x_{30}, x_{12}\rbrack}{q-q^{-1}},
\qquad \qquad
t^{-1}(x_{12}-x_{30}) =
\frac{\lbrack x_{01}, x_{23}\rbrack}{q-q^{-1}}.
\label{eq:0123}
\end{eqnarray}
\end{lemma}
\noindent {\it Proof:} 
To get the equation on the right,
represent 
the standard generators by matrices as in
Proposition
\ref{prop:class1}.
To get the equation on the left,
start with the equation on the right.
In this equation
take the commutator of $x_{12}$ with each side,
and evaluate the result using
Lemma \ref{lem:uvw} (with $u=x_{01}$, $v=x_{12}$, $w=x_{23}$).
\hfill $\Box$ \\


\begin{lemma} 
\label{lem:teq2}
\label{lem:teq3}
On the $\boxtimes_q$-module  $\mathbf V_d(t)$,
\begin{eqnarray}
&&
t (x_{01}-x_{02}) = \frac{\lbrack x_{30}, x_{02}\rbrack}{q-q^{-1}},
\qquad \qquad 
t^{-1} (x_{12}-x_{13}) = 
\frac{\lbrack x_{01}, x_{13}\rbrack}{q-q^{-1}},
\label{eq:four1}
\\
&&
t (x_{23}-x_{20}) = \frac{\lbrack x_{12}, x_{20}\rbrack}{q-q^{-1}},
\qquad \qquad 
t^{-1} (x_{30}-x_{31}) = \frac{\lbrack x_{23}, x_{31}\rbrack}{q-q^{-1}},
\label{eq:four2}
\end{eqnarray}
and also
\begin{eqnarray}
&&
t^{-1}(x_{30}-x_{20}) = \frac{\lbrack x_{20},x_{01}\rbrack}{q-q^{-1}},
\qquad \qquad 
t(x_{01}-x_{31}) = \frac{\lbrack x_{31},x_{12}\rbrack}{q-q^{-1}},
\label{eq:four3}
\\
&&
t^{-1}(x_{12}-x_{02}) = \frac{\lbrack x_{02},x_{23}\rbrack}{q-q^{-1}},
\qquad \qquad 
t(x_{23}-x_{13}) = \frac{\lbrack x_{13},x_{30}\rbrack}{q-q^{-1}}.
\label{eq:four4}
\end{eqnarray}
\end{lemma}
\noindent {\it Proof:}
We first verify the equations on the left 
in (\ref{eq:four1})--(\ref{eq:four4}).
Call these equations 
(\ref{eq:four1}L)--(\ref{eq:four4}L).
To obtain 
(\ref{eq:four4}L),
represent 
the standard generators by matrices as in
Proposition
\ref{prop:class1}.
To obtain
(\ref{eq:four1}L),
in
(\ref{eq:four4}L)
take the commutator of $x_{30}$ with each side,
 and evaluate the result using
Lemma
\ref{lem:uvw} (with $u=x_{23}$,
$v= x_{30}$, $w=x_{02}$)
along with the equation on the left in
(\ref{eq:0123}).
To obtain
(\ref{eq:four3}L),
in 
(\ref{eq:four1}L)
take the commutator of 
each side with
 $x_{20}$,
and evaluate the result using
Lemma
\ref{lem:uvtwo}(ii).
To  obtain
(\ref{eq:four2}L),
 in
(\ref{eq:four3}L)
take the commutator of
$x_{12}$ with each side, and evaluate
the result using Lemma
\ref{lem:uvw} (with $u=x_{01}$, $v=x_{12}$, $w=x_{20}$)
along with the equation on the left in
(\ref{eq:0123}).
We have verified 
the equations 
on the left in
(\ref{eq:four1})--(\ref{eq:four4}).
\medskip

\noindent 
We now verify the
 equations on the right in
(\ref{eq:four1})--(\ref{eq:four4}).
Call these equations
(\ref{eq:four1}R)--(\ref{eq:four4}R).
To obtain
(\ref{eq:four2}R),
represent 
the standard generators by matrices as in
Proposition
\ref{prop:class1}.
To obtain
(\ref{eq:four4}R),
in
(\ref{eq:four2}R)
take the commutator
of each side with $x_{13}$, and
evaluate the result using
Lemma
\ref{lem:uvtwo}(ii).
To obtain
(\ref{eq:four1}R),
in (\ref{eq:four4}R)
take the commutator of $x_{01}$ with each side,
 and evaluate the result using
Lemma
\ref{lem:uvw} 
(with
$u=x_{30}$, $v=x_{01}$, $w=x_{13}$)
along with the equation on the right in
(\ref{eq:0123}).
To obtain
 (\ref{eq:four3}R),
in (\ref{eq:four1}R)
take the commutator
of each side with $x_{31}$, and
evaluate the result using
Lemma
\ref{lem:uvtwo}(ii).
We have verified the equations on the right
in
(\ref{eq:four1})--(\ref{eq:four4}).
The result follows.
\hfill $\Box$ \\

\begin{proposition}  
\label{prop:class2}
Let $V$ denote an evaluation module for $\boxtimes_q$ that
has 
diameter $d$. Then there exists a unique $t \in \F$ 
such that:
\begin{enumerate}
\item[\rm (i)]
$t$ is nonzero and not among
$\lbrace q^{d-2n+1}\rbrace_{n=1}^d$;
\item[\rm (ii)]
the $\boxtimes_q$-module $V$ is isomorphic to $\mathbf V_d(t)$.
\end{enumerate}
\end{proposition}
\noindent {\it Proof:} 
We first show that $t$ exists.
In Section 6, above the first table,
for each $i \in \Z_4$
we used the homomorphism
$\kappa_i: U_q(\mathfrak{sl}_2) 
\to 
\boxtimes_q$ to
turn
$V$ into  a
$U_q(\mathfrak{sl}_2)$-module
 isomorphic to $\mathbf V_d$. Let us take $i=0$.
For this 
$U_q(\mathfrak{sl}_2)$-module 
$V$ let 
 $\lbrace v_n\rbrace_{n=0}^d $ denote an
 $\lbrack x \rbrack_{row}$-basis from
Lemma
\ref{lem:rowbasis}
and Definition
\ref{def:rowbasisinv}.
By Lemma
\ref{lem:sixb} and the construction,
with respect to
 $\lbrace v_n\rbrace_{n=0}^d $ the matrices in
${\rm Mat}_{d+1}(\F)$ 
 that represent
 $x_{12}$, $x_{23}$, $x_{31}$ are
$E_q$,
$K_q$,
$ZE_{q^{-1}}Z$
respectively.
Let $G$ (resp. $S$) (resp. $L$) 
denote the matrix in
${\rm Mat}_{d+1}(\F)$  that represents 
$x_{30}$ (resp. $x_{01}$) (resp. $x_{02}$)
with respect to the basis
 $\lbrace v_n\rbrace_{n=0}^d$.
Using the tables in Section 6
 we find that
(i) $G$ is lower bidiagonal with $(n,n)$-entry
$q^{2n-d}$ for $0 \leq n \leq d$;
(ii) $S$ is tridiagonal;
(iii) $L$ is upper bidiagonal with $(n,n)$-entry
$q^{2n-d}$ for $0 \leq n \leq d$.
On the above matrices we impose
 the defining relations for
$\boxtimes_q$, and solve for $G$, $S$, $L$.
After a brief calculation we find that there exists 
a nonzero $t \in \F$ not among
$\lbrace q^{d-2n+1}\rbrace_{n=1}^d $ such that
\begin{eqnarray*}
G=ZG_{q^{-1}}(t)Z,\qquad 
S=ZS_{q^{-1}}(t^{-1})Z,
\qquad 
L=L_q(t).
\end{eqnarray*}
\noindent Therefore the $\boxtimes_q$-module $V$ is isomorphic
to $\mathbf V_d(t)$.
We have shown that $t$ exists. The scalar $t$ is unique
by
Lemma
\ref{def:evalcom}
and either equation in 
(\ref{eq:0123}).
\hfill $\Box$ \\

\begin{definition}
\label{def:evalp}
\rm 
Let $V$ denote an evaluation module
for $\boxtimes_q$. By the {\it evaluation parameter} of
$V$ we mean the scalar $t$ in Proposition
\ref{prop:class2}.
\end{definition}

\begin{lemma} 
\label{lem:match}
Two evaluation modules for
$\boxtimes_q$ are isomorphic if and only if 
they have the same diameter and same evaluation parameter.
\end{lemma}
\noindent {\it Proof:} 
Consider two isomorphic evaluation modules for $\boxtimes_q$.
They 
have the same diameter, since they have the same 
dimension and the dimension is one more than the
diameter. They have the same evaluation 
parameter by Proposition
\ref{prop:class2}  and Definition
\ref{def:evalp}.
\hfill $\Box$ \\


\medskip
\noindent
Recall the elements $\Upsilon_i$ of $\boxtimes_q$, from
Definition
\ref{def:casi}.

\begin{lemma}
\label{lem:upsi}
Let $V$ denote an evaluation module
for $\boxtimes_q$ that has diameter $d$.
Then
for $i \in \Z_4$ the element
                 $\Upsilon_i$ acts on
$V$ as 
                 $(q^{d+1}+q^{-d-1})I$.
\end{lemma}
\noindent {\it Proof:}
Using the homomorphism $\kappa_i : 
U_q(\mathfrak{sl}_2) \to \boxtimes_q$
we give $V$ 
a $U_q(\mathfrak{sl}_2)$-module structure as in Section 6.
By 
Definition
\ref{def:casi} and the construction,
on $V$ the element
$\Upsilon_i$ agrees with the  
Casimir element of 
$U_q(\mathfrak{sl}_2)$.
The result now follows via Lemma
\ref{lem:nce} and since the
$U_q(\mathfrak{sl}_2)$-module $V$ is isomorphic to 
$\mathbf V_d$.
\hfill $\Box$ \\

\begin{definition}
\label{def:upsilon}
\rm
Let $V$ denote an evaluation module for $\boxtimes_q$ that has
diameter $d$.
By Lemma 
\ref{lem:upsi},
for $i \in \Z_4$ 
the action of 
$\Upsilon_i$ on $V$ is independent
of $i$. Denote this common action by $\Upsilon$. 
Thus on $V$,
\begin{eqnarray}
\label{eq:upsmeaning}
\Upsilon = (q^{d+1}+q^{-d-1})I.
\end{eqnarray}
\end{definition}

\noindent 
We now give 
some identities that involve 
$\Upsilon $.


\begin{lemma} 
\label{lem:upsfour}
Let $V$ denote an evaluation module for $\boxtimes_q$,
with evaluation parameter $t$.
Then on $V$,
\begin{eqnarray*}
&&
\Upsilon = t(x_{01}x_{23}-1)+qx_{30}+q^{-1}x_{12},
 \qquad 
\Upsilon = t^{-1}(x_{12}x_{30}-1)+qx_{01}+q^{-1}x_{23},
\\
&&
\Upsilon = t(x_{23}x_{01}-1)+qx_{12}+q^{-1}x_{30},
 \qquad 
\Upsilon = t^{-1}(x_{30}x_{12}-1)+qx_{23}+q^{-1}x_{01}.
\end{eqnarray*}
\end{lemma}
\noindent {\it Proof:}
Pick $i \in \Z_4$.
In (\ref{eq:cas}) we defined
the Casimir element $\Lambda$ for
$U_q(\mathfrak{sl}_2)$.
In Lemma
\ref{lem:miki}
we described the map
$\kappa_i: 
U_q(\mathfrak{sl}_2) \to \boxtimes_q$.
Apply $\kappa_i$ to
each side of
(\ref{eq:cas}), and evaluate the result
using Definition
\ref{def:casi}
and Definition
\ref{def:upsilon}. 
This shows that on $V$,
\begin{eqnarray}
\Upsilon &=& qx_{i+2,i+3}+q^{-1}x_{i+3,i+1}+qx_{i+1,i+2} 
- qx_{i+2,i+3}x_{i+3,i+1}x_{i+1,i+2}
\nonumber
\\
&=& qx_{i+2,i+3}+q^{-1}x_{i+3,i+1}
+ q(1-x_{i+2,i+3}x_{i+3,i+1})x_{i+1,i+2}.
\label{eq:upswant}
\end{eqnarray}
By 
(\ref{lem:commute}) and 
(\ref{eq:four1}), 
(\ref{eq:four2}) 
the following hold on $V$:
\begin{eqnarray*}
 q(1-x_{i+2,i+3}x_{i+3,i+1})
&=& \frac{\lbrack x_{i+2,i+3}, x_{i+3,i+1} \rbrack}{q-q^{-1}}
\\
&=& t^s (x_{i+3,i}-x_{i+3,i+1}),
\end{eqnarray*}
where $s=(-1)^{i+1}$.
Evaluating 
(\ref{eq:upswant}) using
these comments we find that on $V$,
\begin{eqnarray}
\Upsilon &=& 
qx_{i+2,i+3}+q^{-1}x_{i+3,i+1}
+
 t^s (x_{i+3,i}-x_{i+3,i+1})x_{i+1,i+2}
\nonumber
\\
&=& qx_{i+2,i+3}+q^{-1}x_{i+3,i+1}
+
 t^s (x_{i+3,i}x_{i+1,i+2}-1)
+  t^s(1-x_{i+3,i+1}x_{i+1,i+2}).
\label{eq:upswant2}
\end{eqnarray}
By
(\ref{lem:commute}) and
(\ref{eq:four3}), 
(\ref{eq:four4}) 
the following hold on $V$:
\begin{eqnarray*}
 q(1-x_{i+3,i+1}x_{i+1,i+2})
&=& \frac{\lbrack x_{i+3,i+1}, x_{i+1,i+2} \rbrack}{q-q^{-1}}
\\
&=& t^{-s}(x_{i,i+1}-x_{i+3,i+1}).
\end{eqnarray*}
Evaluating 
(\ref{eq:upswant2}) using these comments
we find that on $V$,
\begin{eqnarray*}
\Upsilon &=& qx_{i+2,i+3}+q^{-1}x_{i+3,i+1}
+
 t^s (x_{i+3,i}x_{i+1,i+2}-1)
+  
q^{-1}(x_{i,i+1}-x_{i+3,i+1})
\\
&=&
 t^s (x_{i+3,i}x_{i+1,i+2}-1)
+
qx_{i+2,i+3}
+  
q^{-1}x_{i,i+1}.
\end{eqnarray*}
The result follows.
\hfill $\Box$ \\

\begin{lemma}
\label{lem:twopair}
Let $V$ denote an evaluation module for $\boxtimes_q$,
with evaluation parameter $t$.
Then on $V$, 
\begin{eqnarray}
&&
\Upsilon = (q+q^{-1})x_{30} + t
\biggl(\frac{q x_{01} x_{23} -q^{-1} x_{23} x_{01}}{q-q^{-1}} -1 \biggr),
\label{eq:firstone}
\\
&&
\Upsilon = (q+q^{-1})x_{01} + t^{-1}
\biggl(\frac{q x_{12} x_{30} -q^{-1} x_{30} x_{12}}{q-q^{-1}} -1\biggr),
\nonumber
\\
&&
\Upsilon
= (q+q^{-1})x_{12} + t
\biggl(\frac{q x_{23} x_{01} -q^{-1} x_{01} x_{23}}{q-q^{-1}} -1 \biggr),
\nonumber
\\
&&\Upsilon= (q+q^{-1})x_{23} + t^{-1}
\biggl(\frac{q x_{30} x_{12} -q^{-1} x_{12} x_{30}}{q-q^{-1}} -1 \biggr).
\nonumber
\end{eqnarray}
\end{lemma}
\noindent {\it Proof:} 
For each displayed equation,
 evaluate the parenthetical expression
using Lemma
\ref{lem:upsfour}
and simplify the result.
\hfill $\Box$ \\

\begin{note}
\label{note:lp} 
\rm 
Let $V$ denote an evaluation module for $\boxtimes_q$.
Among the 
standard generators for $\boxtimes_q$ consider the
following two pairs:
(i) $x_{01}$ and $x_{23}$;
(ii)
$x_{12}$ and $x_{30}$.
Each pair acts on $V$ as a Leonard pair; see
\cite[Example~1.7]{TD00} and
\cite[Theorem~10.3]{qtet}.
Lemma
\ref{lem:twopair} shows how
each Leonard pair determines the other.
\end{note}

\begin{lemma} 
\label{lem:aw}
Let $V$ denote an evaluation module for $\boxtimes_q$,
with evaluation parameter $t$.
Then $x_{01}, x_{23}$
satisfy the following on $V$:
\begin{eqnarray*}
&&
x_{01}^2 x_{23} 
-
(q^2+q^{-2})x_{01} x_{23} x_{01}
+
x_{23} x_{01}^2 \\
&& \qquad \qquad 
=\; -(q-q^{-1})^2 (1+ t^{-1} \Upsilon)x_{01}
+
(q-q^{-1})(q^2-q^{-2})t^{-1},
\\
&&
x_{23}^2 x_{01} 
-
(q^2+q^{-2})x_{23} x_{01} x_{23} 
+
x_{01} x_{23}^2 \\
&& \qquad \qquad 
=\; -(q-q^{-1})^2 (1+ t^{-1} \Upsilon)x_{23}
+
(q-q^{-1})(q^2-q^{-2})t^{-1}.
\end{eqnarray*}
Moreover $x_{12}, x_{30}$
satisfy the following on $V$:
\begin{eqnarray*}
&&
x_{12}^2 x_{30} 
-
(q^2+q^{-2})x_{12} x_{30} x_{12}
+
x_{30} x_{12}^2 \\
&& \qquad \qquad 
=\; -(q-q^{-1})^2 (1+ t \Upsilon)x_{12}
+
(q-q^{-1})(q^2-q^{-2})t,
\\
&&
x_{30}^2 x_{12} 
-
(q^2+q^{-2})x_{30} x_{12} x_{30} 
+
x_{12} x_{30}^2 \\
&& \qquad \qquad 
=\; -(q-q^{-1})^2 (1+ t \Upsilon)x_{30}
+
(q-q^{-1})(q^2-q^{-2})t.
\end{eqnarray*}
\end{lemma}
\noindent {\it Proof:} 
To obtain the first equation,
compute 
(\ref{eq:firstone}) times $qx_{01}$
minus $q^{-1} x_{01} $ times
(\ref{eq:firstone}), and simplify
the result using
\begin{eqnarray*}
\frac{q x_{30}x_{01}-q^{-1}x_{01} x_{30}}{q-q^{-1}}=1.
\end{eqnarray*}
The remaining equations are similarly obtained.
\hfill $\Box$ \\

\begin{note} \rm 
The equations in Lemma
\ref{lem:aw} are the Askey-Wilson
relations \cite[Theorem~1.5]{vidter} for the Leonard pairs
in Note
\ref{note:lp}.
\end{note}

\noindent We end this section with some comments
about the evaluation parameter.

\begin{lemma}
\label{lem:dualt}
Let $V$ denote an evaluation module for $\boxtimes_q$,
with evaluation parameter $t$. Then the $\boxtimes_{q^{-1}}$-module
$V^*$ is an evaluation module with evaluation parameter $t$.
\end{lemma}
\noindent {\it Proof:} 
The 
$\boxtimes_{q^{-1}}$-module $V^*$ is an 
evaluation module 
by Proposition
\ref{lem:vvs}
and
 Proposition 
\ref{prop:sameshape}. Let $t'$ denote the
corresponding evaluation parameter. We show that $t'=t$.
Applying  Lemma
\ref{lem:eval} to $V$ we find that on $V$,
\begin{eqnarray}
\label{eq:tstart}
t(x_{01}-x_{23}) = \frac{\lbrack x_{30}, x_{12} \rbrack}{q-q^{-1}}.
\end{eqnarray}
Applying Lemma
\ref{lem:eval} to the
$\boxtimes_{q^{-1}}$-module
$V^*$ we find that on $V^*$,
\begin{eqnarray}
\label{eq:tprime}
t'(x_{01}-x_{23}) = \frac{\lbrack x_{30}, x_{12} \rbrack}{q^{-1}-q}.
\end{eqnarray}
Let $\Delta $ denote the left-hand side of 
(\ref{eq:tstart}) minus the 
right-hand side of 
(\ref{eq:tstart}).
Then $\Delta V=0$, so $\Delta^\tau V^*=0$ in view of
(\ref{eq:requirement}).
 By this and Proposition
\ref{prop:antiiso} we find that on $V^*$,
\begin{eqnarray}
\label{eq:tend}
t(x_{01}-x_{23}) = \frac{\lbrack x_{12}, x_{30} \rbrack}{q-q^{-1}}.
\end{eqnarray}
We now compare 
(\ref{eq:tprime}) and
(\ref{eq:tend}).
For these equations 
the right-hand sides are the same,
so the left-hand sides agree on $V^*$.
In other words $(t-t')(x_{01}-x_{23})V^*=0$.
But $(x_{01}-x_{23})V^*\not=0$
by Lemma
\ref{def:evalcom} and since $V^*$ is an evaluation module.
Therefore $t=t'$.
\hfill $\Box$ \\


\begin{lemma}
\label{lem:rhoeffect2}
Let $V$ denote an evaluation module for $\boxtimes_q$,
with evaluation parameter $t$.
Then the 
$\boxtimes_q$-module 
${}^\rho V$ is an evaluation module
 with evaluation parameter $t^{-1}$.
\end{lemma}
\noindent {\it Proof:}
The $\boxtimes_q$-module
${}^\rho V$ is an evaluation module
by Lemma
\ref{lem:typediam} and
Lemma
\ref{lem:twistshape}.
To see that 
${}^\rho V$
has evaluation parameter
 $t^{-1}$,
compare the two equations in
Lemma
\ref{lem:eval}, and use
Lemma
\ref{def:evalcom}.
\hfill $\Box$ \\

\begin{corollary}
\label{cor:ttinv}
Let
 $V$ denote an evaluation module for $\boxtimes_q$.
 Then the following
 $\boxtimes_q$-modules are isomorphic: {\rm (i)} $V$;
{\rm (ii)} $V$ twisted via $\rho^2$.
\end{corollary}
\noindent {\it Proof:} Apply 
Lemma \ref{lem:rhoeffect2} 
twice to $V$, and use Lemma
\ref{lem:match}.
\hfill $\Box$ \\

\section{24 bases for an evaluation module}

\noindent 
Let $V$ denote an evaluation module
for $\boxtimes_q$ that has diameter $d$.
In this section we display 24 bases for $V$ that we find attractive.
In Section 11 we show how the standard generators for $\boxtimes_q$
act on these bases. 

\medskip
\noindent For $i \in \Z_4$ consider the flag $\lbrack i \rbrack$
on $V$. We will be discussing component $0$ of this flag.
This component has dimension one. It is an eigenspace 
for each standard generator listed in the table below.
In each case the corresponding eigenvalue is given.

\medskip
\centerline{
\begin{tabular}[t]{c | cccc}
{\rm generator} &  
$x_{i,i+1}$ & $x_{i,i+2}$
&
$x_{i-1,i}$ &
$x_{i+2,i}$
\\
\hline
   {\rm eigenvalue} 
&
$q^{d}$ & $q^d$ 
& $q^{-d}$
& $q^{-d}$
     \end{tabular}}
     \medskip

\begin{definition} 
\label{def:basis}
\rm 
Let $V$ denote an evaluation module for $\boxtimes_q$ that has
diameter $d$. Pick
mutually distinct $i,j,k,\ell$ in  $\Z_4$. 
A basis $\lbrace v_n\rbrace_{n=0}^d$
for $V$ is called an
{\it $\lbrack i,j,k,\ell\rbrack$-basis} whenever:
\begin{enumerate}
\item[\rm (i)] for $0 \leq n \leq d$ the
vector $v_n$ is contained in 
component $n$ of the decomposition
$\lbrack k,\ell\rbrack$ of $V$;
\item[\rm (ii)] $\sum_{n=0}^d v_n$ is contained in
component $0$ of the flag $\lbrack j \rbrack$ on $V$.
\end{enumerate}
\end{definition}

\noindent We will discuss the existence and uniqueness of
the bases in Definition
\ref{def:basis}. We start with uniqueness.

%

\begin{lemma}
\label{lem:basisprime}
Let $V$ denote an evaluation module for $\boxtimes_q$
that has diameter $d$.
Pick
mutually distinct $i,j,k,\ell$ in  $\Z_4$
and let
$\lbrace v_n \rbrace_{n=0}^d$ denote an 
$\lbrack i,j,k,\ell\rbrack$-basis 
 for $V$.
Let 
$\lbrace v'_n \rbrace_{n=0}^d$ denote any vectors in $V$.
Then the following are equivalent:
\begin{enumerate}
\item[\rm (i)] the sequence 
$\lbrace v'_n \rbrace_{n=0}^d$ is an
$\lbrack i,j,k,\ell\rbrack$-basis 
for $V$;
\item[\rm (ii)] there exists $0 \not=\alpha \in \F$
such that $v'_n = \alpha v_n $ for
$0 \leq n \leq d$.
\end{enumerate}
\end{lemma}
\noindent {\it Proof:}
${\rm (i)}\Rightarrow {\rm (ii)}$ 
By assumption $V$ has shape $(1,1,\ldots, 1)$.
\\
\noindent ${\rm (ii)}\Rightarrow {\rm (i)}$ 
Immediate from
Definition
\ref{def:basis}.
\hfill $\Box$ \\

\noindent Let $V$ denote an evaluation module for $\boxtimes_q$
that has
diameter $d$.
In Section 6, for
$i \in \Z_4$ 
we used the homomorphism $\kappa_i : 
 U_q(\mathfrak{sl}_2) \to
\boxtimes_q$
to turn $V$
 into a $U_q(\mathfrak{sl}_2)$-module 
 isomorphic to $\mathbf V_d$. 
 Six bases for  
 this
 $U_q(\mathfrak{sl}_2)$-module
 were displayed in
(\ref{eq:urow}), 
       (\ref{eq:urowinv}).

\begin{lemma} 
\label{lem:comparebasis}
Let $V$ denote an evaluation module for $\boxtimes_q$.
Pick $i \in \Z_4$ and
consider the corresponding
 $U_q(\mathfrak{sl}_2)$-module $V$ as above.
In each row of the table below
we display a basis
for the 
 $U_q(\mathfrak{sl}_2)$-module $V$ from 
{\rm (\ref{eq:urow}), 
       (\ref{eq:urowinv})}
 and a basis for the
$\boxtimes_q$-module $V$ from Definition
\ref{def:basis}. These two bases are the same.

\medskip

\centerline{
\begin{tabular}[t]{c|c}
    {\rm basis for $V$ from
(\ref{eq:urow}), 
       (\ref{eq:urowinv})
    } & {\rm basis for $V$ from Def. \ref{def:basis}}
    \\ \hline  \hline
$\lbrack x\rbrack_{row}$ & $\lbrack i,i+1,i+2,i+3 \rbrack$
\\
$\lbrack x\rbrack^{inv}_{row}$ & $\lbrack i,i+1,i+3,i+2 \rbrack$
\\
\hline
$\lbrack y\rbrack_{row}$ & $\lbrack i,i+2,i+3,i+1 \rbrack$
\\
$\lbrack y\rbrack^{inv}_{row}$ & $\lbrack i,i+2,i+1,i+3 \rbrack$
\\
\hline
$\lbrack z\rbrack_{row}$ & $\lbrack i,i+3,i+1,i+2 \rbrack$
\\
$\lbrack z\rbrack^{inv}_{row}$ & $\lbrack i,i+3,i+2,i+1 \rbrack$
      \end{tabular}}
   \medskip

\end{lemma}
\noindent {\it Proof:}
By Lemma
\ref{lem:miki} and the construction.
\hfill $\Box$ \\

\begin{lemma} 
\label{lem:24exist}
Let $V$ denote an evaluation module for $\boxtimes_q$, and
pick mutually distinct $i,j,k,\ell$ in $\Z_4$. Then
there exists an
$\lbrack i,j,k,\ell\rbrack$-basis
for $V$.
\end{lemma} 
\noindent {\it Proof:}
Immediate from Lemma
\ref{lem:comparebasis}.
\hfill $\Box$ \\

\begin{note}
\label{note:10pt5}
\rm The basis for
$\mathbf V_d(t)$ given in Proposition
\ref{prop:class1} is a $\lbrack 0,1,2,3\rbrack$-basis.
\end{note}

\noindent 
Let $V$ denote an evaluation module for $\boxtimes_q$.
In Definition
\ref{def:basis}
we gave 24 bases for $V$.
In Section 11 we will compute the
matrices that represent the standard generators
with respect to these bases.
We now mention some results that
will facilitate this computation.

\begin{lemma}
\label{lem:24Z}
Let $V$ denote an evaluation module for $\boxtimes_q$, and
pick mutually distinct $i,j,k,\ell$ in $\Z_4$. 
Then each $\lbrack i,j,k,\ell\rbrack$-basis
for $V$ is the inversion of an 
$\lbrack i,j,\ell,k\rbrack$-basis for 
$V$.
\end{lemma} 
\noindent {\it Proof:} 
By Definition
\ref{def:basis} and the meaning of inversion.
\hfill $\Box$ \\

\begin{corollary}
\label{cor:zcz}
Let $V$ denote an evaluation module for $\boxtimes_q$, and
pick mutually distinct $i,j,k,\ell$ in $\Z_4$. 
For each standard generator $x_{rs}$ consider the
following two matrices:
\begin{enumerate}
\item[\rm (i)] the matrix
 that represents $x_{rs}$
with respect to an
$\lbrack i,j,k,\ell \rbrack$-basis for 
$V$;
\item[\rm (ii)] the matrix 
that represents $x_{rs}$
with respect to an
$\lbrack i,j,\ell,k\rbrack$-basis for 
$V$.
\end{enumerate}
Each of these matrices is obtained from the other via 
 conjugation by $Z$.
\end{corollary}
\noindent {\it Proof:}
By Lemma
\ref{lem:24Z}
and linear algebra.
\hfill $\Box$ \\

\begin{lemma}
\label{lem:rhoeffect}
Let $V$ denote an evaluation module for $\boxtimes_q$, and
pick mutually distinct $i,j,k,\ell$ in $\Z_4$.
 Then the following are the same:
\begin{enumerate}
\item[\rm (i)] an 
$\lbrack i,j,k,\ell\rbrack$-basis
for $V$;
\item[\rm (ii)] an 
$\lbrack i+1,j+1,k+1,\ell+1\rbrack$-basis
for ${}^{\rho}V$.
\end{enumerate}
\end{lemma}
\noindent {\it Proof:} Use
Lemma
\ref{lem:twistdec} and
Lemma
\ref{lem:twistflag}.
\hfill $\Box$ \\

\begin{lemma}
\label{lem:tshift}
Consider the $\boxtimes_q$-module $\mathbf V_d(t)$.
Pick mutually distinct $i,j,k,\ell$ in  $\Z_4$. 
Then for each standard generator $x_{rs}$ the following are the same:
\begin{enumerate}
\item[\rm (i)] the matrix that represents $x_{rs}$ with
respect to an
$\lbrack i,j,k,\ell\rbrack$-basis
for $\mathbf V_d(t)$;
\item[\rm (ii)] the matrix that represents $x_{r+1,s+1}$ with
respect to an
$\lbrack i+1,j+1,k+1,\ell+1\rbrack$-basis
for
$\mathbf V_d(t^{-1})$.
\end{enumerate}
\end{lemma}
\noindent {\it Proof:} 
Use
Lemma
\ref{lem:rhoxijtwist},
Lemma
\ref{lem:rhoeffect2},
and Lemma
\ref{lem:rhoeffect}.
\hfill $\Box$ \\

\begin{corollary}
\label{cor:tshift2}
Let $V$ denote an evaluation module for
$\boxtimes_q$.
Pick mutually distinct $i,j,k,\ell$ in  $\Z_4$. 
Then for each standard generator $x_{rs}$ the following are the same:
\begin{enumerate}
\item[\rm (i)] the matrix that represents $x_{rs}$ with
respect to an
$\lbrack i,j,k,\ell\rbrack$-basis
for $V$;
\item[\rm (ii)] the matrix that represents $x_{r+2,s+2}$ with
respect to an
$\lbrack i+2,j+2,k+2,\ell+2\rbrack$-basis
for
$V$.
\end{enumerate}
\end{corollary}
\noindent {\it Proof:}  
Apply
Lemma \ref{lem:tshift} twice, or use
Corollary
\ref{cor:ttinv} along with
Lemma
\ref{lem:rhoeffect}.
\hfill $\Box$ \\

\section{The action of the standard 
generators on the
24 bases}

\noindent Let $V$ denote an evaluation module
for $\boxtimes_q$.
In
Definition
\ref{def:basis}
we defined 24 bases for $V$. In this section we
display the matrices that represent the standard generators 
$x_{ij}$ with respect to these bases. 

\medskip
\noindent The matrices displayed in the following 
theorem are defined in Appendix I.

\begin{theorem}
\label{thm:xijaction}
Let $V$ denote an evaluation module for $\boxtimes_q$,
with diameter $d$ and evaluation parameter $t$.
In the tables below,
we display the matrices in 
${\rm Mat}_{d+1}(\F)$ 
that represent the standard generators $x_{ij}$ 
with respect to the 24 bases for $V$ from Definition
\ref{def:basis}.
Pick $r \in \Z_4$, and first assume that $r$ is even. Then

\medskip
\centerline{
\begin{tabular}[t]{c|cccc}
{\rm basis} & 
  $x_{r,r+1}$ 
  &
  $x_{r+1,r+2}$ 
  &
  $x_{r+2,r+3}$ 
  &
  $x_{r+3,r}$ 
 %
\\
\hline
\hline
\\
$\lbrack r,r+1,r+2,r+3\rbrack$
&
$ZS_{q^{-1}}(t^{-1})Z$
&
$E_q$
&
$K_q$
&
$Z G_{q^{-1}}(t)Z$
\\
\\
$\lbrack r+1,r,r+2,r+3\rbrack$
&
$S_q(t^{-1})$
&
$G_q(t)$
&
$K_q$
&
$ZE_{q^{-1}}Z$
\\
\\
$\lbrack r,r+1,r+3,r+2\rbrack$
&
$S_{q^{-1}}(t^{-1})$
&
$ZE_qZ$
&
$K_{q^{-1}}$
&
$ G_{q^{-1}}(t)$
\\
\\
$\lbrack r+1,r,r+3,r+2\rbrack$
&
$ZS_q(t^{-1})Z$
&
$ZG_q(t)Z$
&
$K_{q^{-1}}$
&
$E_{q^{-1}}$
\\
\\
\hline
\\
$\lbrack r,r+2,r+1,r+3\rbrack$
&
$ F_q(t^{-1}) $
&
$ E_{q^{-1}} $
&
$Z E_q Z$
&
$Z F_{q^{-1}}(t) Z$
\\
\\
$\lbrack r+2,r,r+1,r+3\rbrack$
&
$E_q$  & $F_{q^{-1}}(t)$ & $ZF_q(t^{-1})Z$ & $ZE_{q^{-1}}Z$
     \end{tabular}}

\bigskip

\bigskip

\centerline{
\begin{tabular}[t]{c|cccc}
{\rm basis} & 
  $x_{r,r+2}$ 
  &
  $x_{r+1,r+3}$ 
  &
  $x_{r+2,r}$ 
  &
  $x_{r+3,r+1}$ 
\\
\hline
\hline
\\
$\lbrack r,r+1,r+2,r+3\rbrack$
&
$L_q(t)$ 
&
$(ZE_{q^{-1}}Z)^{-1}$
&
$(L_q(t))^{-1}$
&
$ZE_{q^{-1}}Z$
\\
\\
$\lbrack r+1,r,r+2,r+3\rbrack$
&
$E_q$
&
$(Z L_{q^{-1}}(t)Z)^{-1}$
&
$(E_q)^{-1}$
&
$Z L_{q^{-1}}(t)Z$
\\
\\
$\lbrack r,r+1,r+3,r+2\rbrack$
&
$ZL_q(t)Z$ 
&
$(E_{q^{-1}})^{-1}$
&
$(ZL_q(t)Z)^{-1}$
&
$E_{q^{-1}}$
\\
\\
$\lbrack r+1,r,r+3,r+2\rbrack$
&
$ZE_qZ$
&
$(L_{q^{-1}}(t))^{-1}$
&
$(ZE_qZ)^{-1}$
&
$ L_{q^{-1}}(t)$
\\
\\
\hline
\\
$\lbrack r,r+2,r+1,r+3\rbrack$
&
$ Z M_{q^{-1}}(t^{-1}) Z $
&
$ K_q $
&
$ M_q(t) $
&
$ K_{q^{-1}} $
\\
\\
$\lbrack r+2,r,r+1,r+3\rbrack$
&
$ZM_q(t)Z$ & $K_q$ & $M_{q^{-1}}(t^{-1})$ & $K_{q^{-1}}$
     \end{tabular}}
     \medskip

\noindent
Next assume that $r$ is odd. Then in the above tables
replace $t$ by $t^{-1}$.

\end{theorem}
\noindent {\it Proof:}  
By Lemma
\ref{lem:tshift} we may assume that $r=0$. We now consider
six cases. 
In the following discussion all bases mentioned are for $V$.
\\
\noindent Basis $\lbrack 0,1,2,3\rbrack$. By Proposition
\ref{prop:class1}
and
Note
\ref{note:10pt5}. 
\\
\noindent Basis $\lbrack 1,0,2,3\rbrack$.
We first verify the data for $x_{23}$, $x_{30}$, $x_{02}$.
Using the map $\kappa_1:
U_q(\mathfrak{sl}_2) \to \boxtimes_q$
we turn 
$V$ into a
$U_q(\mathfrak{sl}_2)$-module isomorphic to 
$\mathbf V_d$.
By Lemma
\ref{lem:comparebasis} the bases
$\lbrack 1,0,2,3\rbrack $
and
$\lbrack z \rbrack_{row}$ are the same.
By Lemma
\ref{lem:miki}
the equations
$x_{23}=z$,
$x_{30}=x$,
$x_{02}=y$ hold on $V$.
By Lemma
\ref{lem:sixb}
the matrices representing $z$, $x$, $y$ with respect to
$\lbrack z \rbrack_{row}$ are 
$K_q$, $ZE_{q^{-1}}Z$, $E_q$ respectively.
By these comments the matrices representing
$x_{23}$,
$x_{30}$,
$x_{02}$
with respect to
$\lbrack 1,0,2,3\rbrack $
are 
$K_q$, $ZE_{q^{-1}}Z$, $E_q$ respectively.
By Definition
\ref{def:qtet}(i) the generator $x_{20}$ is the inverse of
$x_{02}$. Therefore the matrix representing $x_{20}$ with
respect to $\lbrack 1,0,2,3\rbrack$ is $(E_q)^{-1}$. 
To get the
matrix representing $x_{01}$ with respect to 
$\lbrack 1,0,2,3\rbrack$ 
use the equation on the left in
(\ref{eq:four1}).
To get the
matrix representing $x_{12}$ with respect to 
$\lbrack 1,0,2,3\rbrack$ 
use the equation on the left in
(\ref{eq:four4}).
We now verify the data for $x_{31}$ and $x_{13}$. Let $L$ denote
the matrix representing $x_{31}$ with respect to
$\lbrack 1,0,2,3\rbrack$. By the second table of Section
6, $L$ is lower bidiagonal with $(n,n)$-entry $q^{2n-d}$ for
$0 \leq n \leq d$. On this matrix we impose the defining
relations for $\boxtimes_q$. After a brief calculation
we obtain $L=ZL_{q^{-1}}(t)Z$. By Definition 
\ref{def:qtet}(i)  $x_{13}$ is the inverse of $x_{31}$.
Therefore the matrix representing $x_{13}$ with respect to
$\lbrack 1,0,2,3\rbrack$ is
$(ZL_{q^{-1}}(t)Z)^{-1}$.
\\
\noindent Basis $\lbrack 0,1,3,2\rbrack$.
Use Corollary
\ref{cor:zcz}.
\\
\noindent Basis $\lbrack 1,0,3,2\rbrack$. 
Use Corollary
\ref{cor:zcz}.
\\
\noindent Basis $\lbrack 0,2,1,3\rbrack$.
We first verify the data for $x_{12}$, $x_{23}$, $x_{31}$.
Using the map $\kappa_0:
U_q(\mathfrak{sl}_2) \to \boxtimes_q$
we turn 
$V$ into a
$U_q(\mathfrak{sl}_2)$-module isomorphic to 
$\mathbf V_d$.
By Lemma
\ref{lem:comparebasis} the bases
$\lbrack 0,2,1,3\rbrack $
and
$\lbrack y \rbrack^{inv}_{row}$ are the same.
By Lemma
\ref{lem:miki}
the equations
$x_{12}=z$,
$x_{23}=x$,
$x_{31}=y$ hold on $V$.
By Lemma
\ref{lem:sixb}
the matrices representing $z$, $x$, $y$ with respect to
$\lbrack y \rbrack^{inv}_{row}$ are 
$E_{q^{-1}}$,
$ZE_qZ$, 
$K_{q^{-1}}$ respectively.
By these comments the matrices representing
$x_{12}$,
$x_{23}$,
$x_{31}$
with respect to
$\lbrack 0,2,1,3\rbrack $
are 
$E_{q^{-1}}$,
$ZE_qZ$, 
$K_{q^{-1}}$ respectively.
By Definition
\ref{def:qtet}(i) the generator $x_{13}$ is the inverse of
$x_{31}$. Therefore the matrix representing $x_{13}$ with
respect to $\lbrack 0,2,1,3\rbrack$ is $(K_{q^{-1}})^{-1}$, which
is equal to
 $K_q$.
To get the
matrix representing $x_{01}$ with respect to 
$\lbrack 0,2,1,3\rbrack$ 
use the equation on the right
in (\ref{eq:four3}).
To get the
matrix representing $x_{30}$ with respect to 
$\lbrack 0,2,1,3\rbrack$ 
use the equation on the right
in (\ref{eq:four2}).
We now verify the data for $x_{20}$ and $x_{02}$. Let $M$ 
(resp. $M'$) denote
the matrix representing $x_{20}$ (resp. $x_{02}$) with respect to
$\lbrack 0,2,1,3\rbrack$. 
On this matrix we impose the defining
relations for $\boxtimes_q$. After a brief calculation
we obtain 
$M=M_q(t)$ and
$M'=ZM_{q^{-1}}(t^{-1})Z$.
\\
\noindent Basis $\lbrack 2,0,1,3\rbrack$.
Use Corollary
\ref{cor:zcz} and
Corollary \ref{cor:tshift2}.
\hfill $\Box$ \\

\section{A normalization of the 24 bases}

\noindent Let $V$ denote an evaluation module for 
$\boxtimes_q$. 
In 
Definition
\ref{def:basis} we gave 
24 bases for $V$.
In Section 13 we will give the transition matrices
between certain pairs of bases among the 24. Before
doing this, we first normalize our bases.

\begin{definition}
\label{def:cornervector}
\rm
Let $V$ denote an evaluation module for $\boxtimes_q$.
For $i \in \Z_4$ let $\eta_i$ (resp. $\eta^*_i$)
denote a nonzero vector in component 0 of the flag
$\lbrack i \rbrack$ on $V$ (resp. $V^*$).
\end{definition}

\begin{lemma} 
Let $V$ denote an evaluation module for $\boxtimes_q$ that has
diameter $d$.
Then the following {\rm (i), (ii)} 
hold for $i \in \Z_4$.
\begin{enumerate}
\item[\rm (i)]
The vector $\eta_i$ is an eigenvector for
each standard generator listed in the table below. 
In each case the corresponding eigenvalue is given.

\medskip
\centerline{
\begin{tabular}[t]{c | cccc}
{\rm generator} &  
$x_{i,i+1}$ & $x_{i,i+2}$
&
$x_{i-1,i}$ &
$x_{i+2,i}$
\\
\hline
   {\rm eigenvalue} 
&
$q^{d}$ & $q^d$ 
& $q^{-d}$
& $q^{-d}$
     \end{tabular}}
     \medskip

\item[\rm (ii)]
The vector $\eta^*_i$ is an eigenvector for
each standard generator listed in the table below. 
In each case the corresponding eigenvalue is given.

\medskip
\centerline{
\begin{tabular}[t]{c | cccc}
{\rm generator} &  
$x_{i,i+1}$ & $x_{i,i+2}$
&
$x_{i-1,i}$ &
$x_{i+2,i}$
\\
\hline
   {\rm eigenvalue} 
&
$q^{-d}$ & $q^{-d}$ 
& $q^{d}$
& $q^{d}$
     \end{tabular}}
     \medskip
\end{enumerate}
\end{lemma}
\noindent {\it Proof:} (i) By Definition
\ref{def:cornervector} and the paragraph above 
Definition
\ref{def:basis}.
\\
\noindent (ii) Apply (i) above to the
$\boxtimes_{q^{-1}}$-module $V^*$.
\hfill $\Box$ \\

\begin{lemma} 
Referring to Definition
\ref{def:cornervector}
the following {\rm (i), (ii)} hold.
\begin{enumerate}
\item[\rm (i)] For distinct $i,j \in \Z_4$ we have
$(\eta_i,\eta^*_j)\not=0$.
\item[\rm (ii)] For $i \in \Z_4$
we have $(\eta_i, \eta^*_i ) = 0$.
\end{enumerate}
\end{lemma}
\noindent {\it Proof:}  
(i) The vector $\eta_i$ is a basis for component $0$ of the 
decomposition $\lbrack i,j\rbrack$ for $V$.
The vector $\eta^*_j$ is a basis for component $0$ of the 
decomposition $\lbrack j,i\rbrack$ for $V^*$.
By Proposition
\ref{prop:dualdec} the decomposition
 $\lbrack i,j\rbrack$ for $V$ is dual to
the decomposition
 $\lbrack j,i\rbrack$ for $V^*$.
Therefore 
$(\eta_i,\eta^*_j)\not=0$.
\\
\noindent (ii) Use Proposition
\ref{prop:flagdual} with $n=0$. Recall that
$V$ is nontrivial so it has diameter $d\geq 1$.
\hfill $\Box$ \\

\begin{lemma}
\label{lem:normexist}
Let $V$ denote an evaluation module for $\boxtimes_q$ that has
diameter $d$. 
Given mutually distinct
$i,j,k,\ell \in \Z_4$ there exists a unique
basis $\lbrace u_n \rbrace_{n=0}^d$
for $V$ such that {\rm (i), (ii)} hold below:
\begin{enumerate}
\item[\rm (i)] for $0 \leq n \leq d$ the vector $u_n$ is contained in
component $n$ of the decomposition $\lbrack k,\ell\rbrack $ of
$V$;
\item[\rm (ii)] $\eta_j = \sum_{n=0}^d u_n$.
\end{enumerate}
We denote this basis by
$\lbrack i,j,k,\ell\rbrack$.
\end{lemma}
\noindent {\it Proof:}  
We first show that the basis
$\lbrack i,j,k,\ell\rbrack$ exists. By Lemma
\ref{lem:24exist} there exists an
$\lbrack i,j,k,\ell\rbrack$-basis for $V$.
Denote this basis by $\lbrace v_n \rbrace_{n=0}^d$.
Define $v=\sum_{n=0}^d v_n$ and note that $v\not=0$.
By Definition
\ref{def:basis} $v$ is contained in component $0$
of the flag $\lbrack j \rbrack $ for $V$.
By construction $\eta_j$ is a basis for component $0$ of the
flag $\lbrack j \rbrack $ for $V$.
Therefore there exists $0 \not=\alpha \in \F$ such that 
$v = \alpha \eta_j$. Define $u_n = \alpha^{-1}v_n$  for $0 \leq n \leq d$.
Then $\lbrace u_n \rbrace_{n=0}^d$ is a basis for $V$ that
satisfies the requirements (i), (ii).
One checks using Lemma
\ref{lem:basisprime}
that the basis
$\lbrack i,j,k,\ell\rbrack$ is unique.
\hfill $\Box$ \\

\begin{lemma} 
\label{lem:invnorm}
Let $V$ denote an evaluation module for $\boxtimes_q$, and
pick mutually distinct $i,j,k,\ell$ in $\Z_4$. Then the
bases $\lbrack i,j,k,\ell\rbrack$  and
$\lbrack i,j,\ell,k\rbrack$  for $V$ are the inversion
of each other.
\end{lemma}
\noindent {\it Proof:}   By Lemma
\ref{lem:normexist} and the meaning of inversion.
\hfill $\Box$ \\

\begin{proposition}
\label{prop:comp0d}
Let $V$ denote an evaluation module for $\boxtimes_q$
that has diameter $d$.
Pick mutually distinct $i,j,k,\ell$ in
$\Z_4$ and consider the basis
$\lbrack i,j,k,\ell\rbrack$ of $V$. For this basis the components
0 and $d$ are given below.

\medskip
\centerline{
\begin{tabular}[t]{ccc}
{\rm component $0$} &  &
   {\rm component $d$} 
   \\ \hline 
 $\frac{(\eta_j, \eta^*_\ell)}{(\eta_k,\eta^*_\ell)} \eta_k$
 &
 &
 $ \frac{(\eta_j, \eta^*_k)}{(\eta_\ell,\eta^*_k)} \eta_\ell$
     \end{tabular}}
     \medskip
\end{proposition}
\noindent {\it Proof:}  
Denote the basis by $\lbrace u_n \rbrace_{n=0}^d$.
Recall that for $0 \leq n \leq d$ the vector
$u_n$ is contained in component $n$ of the decomposition
$\lbrack k,\ell\rbrack$ of $V$.
By construction $\eta_k$ (resp. $\eta_\ell$) is a basis
for component $0$ (resp. $d$)  of the decomposition
$\lbrack k, \ell \rbrack$ of $V$.
Similarly 
 $\eta^*_k$ (resp. $\eta^*_\ell$) is a basis
for component $0$ (resp. $d$)  of the decomposition
$\lbrack k, \ell \rbrack$ of $V^*$.
By Proposition
\ref{prop:dualdec} the
decomposition
$\lbrack k,\ell\rbrack $ of $V$ is dual to
decomposition
$\lbrack \ell,k\rbrack $ of $V^*$.
Therefore $(u_n,\eta^*_\ell)=0$ for $1 \leq n\leq d$ and
$(u_n,\eta^*_k)=0$ for $0 \leq n\leq d-1$.
Moreover there exist $\alpha, \beta \in \F$ such that
$u_0 = \alpha\eta_k$ and $u_d = \beta\eta_\ell$.
Now using Lemma
\ref{lem:normexist}(ii),
\begin{eqnarray*}
(\eta_j,\eta^*_\ell) 
= 
\sum_{n=0}^d (u_n, \eta^*_\ell)
= 
(u_0,\eta^*_\ell)
= 
\alpha (\eta_k,\eta^*_\ell)
\end{eqnarray*}
so
$\alpha = (\eta_j,\eta^*_\ell)/(\eta_k,\eta^*_\ell)$.
Similarly
\begin{eqnarray*}
(\eta_j,\eta^*_k) 
= 
\sum_{n=0}^d (u_n, \eta^*_k)
= 
(u_d,\eta^*_k)
= 
\beta (\eta_\ell,\eta^*_k)
\end{eqnarray*}
so
$\beta = (\eta_j,\eta^*_k)/(\eta_\ell,\eta^*_k)$.
The result follows.
\hfill $\Box$ \\

\section{Transition matrices between the 24 bases}

\noindent
Let $V$ denote an evaluation module for 
$\boxtimes_q$ that has
diameter $d$.
Recall the 24 bases for $V$ from
Lemma
\ref{lem:normexist}.
In this section we compute the transition matrices
between certain pairs of bases among the 24. First we clarify
 a
few terms.

\medskip
\noindent 
Suppose we are given two bases for $V$, denoted
$\lbrace u_n\rbrace_{n=0}^d$
and 
$\lbrace v_n\rbrace_{n=0}^d$. By the {\it transition matrix from
$\lbrace u_n\rbrace_{n=0}^d$
to
$\lbrace v_n\rbrace_{n=0}^d$} we mean the
matrix $S \in {\rm Mat}_{d+1}(\F)$ such that
$v_n = \sum_{r=0}^d S_{rn} u_r$ for
$0 \leq n \leq d$.
Let $S$ denote the transition matrix from 
$\lbrace u_n\rbrace_{n=0}^d$
to
$\lbrace v_n\rbrace_{n=0}^d$.
Then $S^{-1}$ exists and equals the transition matrix
from
$\lbrace v_n\rbrace_{n=0}^d$
to
$\lbrace u_n\rbrace_{n=0}^d$.

\medskip
\noindent 
Let $\lbrace w_n\rbrace_{n=0}^d$ denote a basis for
$V$ and let $T$ denote the transition matrix from
$\lbrace v_n\rbrace_{n=0}^d$
to
$\lbrace w_n\rbrace_{n=0}^d$.
Then $ST$ is the transition matrix from
$\lbrace u_n\rbrace_{n=0}^d$
to
$\lbrace w_n\rbrace_{n=0}^d$.

\medskip
\noindent 
Let ${ A \in \rm End}(V)$ 
and let 
$ M$ denote the matrix in ${\rm Mat}_{d+1}(\F)$
that represents $A$ with respect to
$\lbrace u_n\rbrace_{n=0}^d$. Then
the matrix $S^{-1}MS$  represents $A$
with respect to 
$\lbrace v_n\rbrace_{n=0}^d$.
\medskip

\noindent The matrix
$Z$ 
is defined in Appendix I.
Let $\lbrace v_n\rbrace_{n=0}^d$ denote a basis
for $V$ and consider the inverted basis
$\lbrace v_{d-n}\rbrace_{n=0}^d$.
 Then $Z$ is the transition matrix from
 $\lbrace v_n\rbrace_{n=0}^d$ to
$\lbrace v_{d-n}\rbrace_{n=0}^d$.

\medskip

\noindent
Pick mutually distinct $i,j,k,\ell$ in 
$\Z_4$ and consider the basis
$\lbrack i,j,k,\ell \rbrack $ of $V$.
 We will display the transition matrix 
from this  basis
to each of the bases 
\begin{eqnarray*}
\lbrack j,i,k,\ell \rbrack,
\qquad 
\lbrack i,k,j,\ell \rbrack,
\qquad 
\lbrack i,j,\ell,k \rbrack.
\end{eqnarray*}
As we will see, the first transition matrix
is diagonal, the second is lower triangular,
and the third one is equal to $Z$.

\medskip
\noindent We now consider the transitions
of type
$\lbrack i,j,k,\ell \rbrack \to
\lbrack j,i,k,\ell \rbrack$. We will be discussing
the matrices $D_q(t)$ and $\mathcal D_q(t)$ defined in
Appendix I.

\begin{theorem} 
\label{thm:diagtrans}
Let $V$ denote an evaluation module for $\boxtimes_q$,
with diameter $d$ and evaluation parameter $t$.
In the table below we display some transition
matrices between the 24 bases for $V$ from
Lemma
\ref{lem:normexist}.
 Each transition matrix
is diagonal. 
Pick $r \in \Z_4$,
and first assume that $r$ is even.

\medskip
\centerline{
\begin{tabular}[t]{c|c}
{\rm transition}  & 
   {\rm transition matrix}
   \\ \hline
   \hline
\\
$ \lbrack r,r+1,r+2,r+3\rbrack \to
 \lbrack r+1,r,r+2,r+3\rbrack $ &
$
D_q(t)
\frac{
(\eta_r, \eta^*_{r+3})}
{
(\eta_{r+1}, \eta^*_{r+3})}
$
\\
\\
$ \lbrack r+1,r,r+2,r+3\rbrack \to
 \lbrack r,r+1,r+2,r+3\rbrack $ &
$
\bigl(D_q(t)\bigr)^{-1}
\frac{
(\eta_{r+1}, \eta^*_{r+3})}
{(\eta_r, \eta^*_{r+3})}
$
\\
\\
$ \lbrack r,r+1,r+3,r+2\rbrack \to
 \lbrack r+1,r,r+3,r+2\rbrack $ &
$
\bigl(D_{q^{-1}}(t)\bigr)^{-1}
\frac{
(\eta_{r}, \eta^*_{r+2})}
{(\eta_{r+1}, \eta^*_{r+2})}
$
\\
\\
$ \lbrack r+1,r,r+3,r+2\rbrack \to
 \lbrack r,r+1,r+3,r+2\rbrack $ &
$
D_{q^{-1}}(t)
\frac{
(\eta_{r+1}, \eta^*_{r+2})}
{
(\eta_{r}, \eta^*_{r+2})}
$
\\
\\
\hline
\\
$ \lbrack r,r+2,r+1,r+3\rbrack \to
 \lbrack r+2,r,r+1,r+3\rbrack $ &
$
\mathcal D_q(t)
\frac{(\eta_r,\eta^*_{r+3})}{(\eta_{r+2},\eta^*_{r+3})}
$
\\
\\
$ \lbrack r+2,r,r+1,r+3\rbrack \to
 \lbrack r,r+2,r+1,r+3\rbrack $ &
$
\bigl(\mathcal D_q(t)\bigr)^{-1}
\frac{
(\eta_{r+2},\eta^*_{r+3})
}{
(\eta_r,\eta^*_{r+3})
}
$
     \end{tabular}}
     \medskip
\noindent Next assume that $r$ is odd. Then in the above
table replace $t$ by $t^{-1}$.
\end{theorem}
\noindent {\it Proof:}  
By Lemma
\ref{lem:rhoeffect2} and the construction,
we may assume that 
$r=0$.
We consider six cases. In the following discussion all bases
mentioned are for $V$.
\\
\noindent $\lbrack 0,1,2,3\rbrack \to 
 \lbrack 1,0,2,3\rbrack$.
Let
$\lbrace u_n \rbrace_{n=0}^d$ and
$\lbrace v_n \rbrace_{n=0}^d$
denote the bases
$\lbrack 0,1,2,3\rbrack$
and
 $\lbrack 1,0,2,3\rbrack$ respectively.
Let $D \in 
{\rm Mat}_{d+1}(\F)$
 denote the transition matrix
from
$\lbrace u_n \rbrace_{n=0}^d$ to
$\lbrace v_n \rbrace_{n=0}^d$.
For $0 \leq n \leq d$ the vectors $u_n$, $v_n$ are contained
in component $n$ of the decomposition $\lbrack 2,3\rbrack$.
Therefore $D$ is diagonal. By Theorem
\ref{thm:xijaction} the matrix representing
$x_{12}$ with respect to $\lbrace u_n \rbrace_{n=0}^d$
is equal to $E_q$, and the 
 matrix representing
$x_{12}$ with respect to $\lbrace v_n \rbrace_{n=0}^d$
is equal to $G_q(t)$.
Therefore  $E_q D = D G_q(t)$. Comparing this with
(\ref{eq:GE})  we find that 
there exists $0 \not=\alpha \in \F$ such that
 $D=\alpha D_q(t)$.
We now find $\alpha$. 
The
$(0,0)$-entry of $D_q(t)$ is $1$, so
the $(0,0)$-entry of $D$ is $\alpha$.
Therefore $v_0 = \alpha u_0$.
By Proposition
\ref{prop:comp0d}, both
\begin{eqnarray*}
u_0 = \frac{(\eta_{1},\eta^*_{3})}{(\eta_{2},\eta^*_{3})}
\eta_{2},
\qquad \qquad
v_0 = \frac{(\eta_0,\eta^*_{3})}{(\eta_{2},\eta^*_{3})}
\eta_{2}.
\end{eqnarray*}
Therefore $\alpha = 
(\eta_0,\eta^*_{3})/(\eta_{1},\eta^*_{3})$.
\\
\noindent $\lbrack 1,0,2,3\rbrack \to 
 \lbrack 0,1,2,3\rbrack$.
This transition matrix is the inverse of the transition matrix
for
$\lbrack 0,1,2,3\rbrack \to 
 \lbrack 1,0,2,3\rbrack$.
\\
\noindent $\lbrack 0,1,3,2\rbrack \to 
 \lbrack 1,0,3,2\rbrack$.
Let
$\lbrace u_n \rbrace_{n=0}^d$ and
$\lbrace v_n \rbrace_{n=0}^d$
denote the bases
$\lbrack 0,1,3,2\rbrack$
and
 $\lbrack 1,0,3,2\rbrack$ respectively.
Let $D \in 
{\rm Mat}_{d+1}(\F)$
 denote the transition matrix
from
$\lbrace u_n \rbrace_{n=0}^d$ to
$\lbrace v_n \rbrace_{n=0}^d$.
For $0 \leq n \leq d$ the vectors $u_n$, $v_n$ are contained
in component $n$ of the decomposition $\lbrack 3,2\rbrack$.
Therefore $D$ is diagonal. By Theorem
\ref{thm:xijaction} the matrix representing
$x_{30}$ with respect to $\lbrace u_n \rbrace_{n=0}^d$
is equal to $G_{q^{-1}}(t)$, and the 
 matrix representing
$x_{30}$ with respect to $\lbrace v_n \rbrace_{n=0}^d$
is equal to $E_{q^{-1}}$.
Therefore  $G_{q^{-1}}(t) D = D E_{q^{-1}}$. Comparing this with
(\ref{eq:GE})  we find that
there exists $0 \not=\alpha \in \F$ such that
$D=\alpha (D_{q^{-1}}(t))^{-1}$.
We now find $\alpha$.
The matrix  
 $D_{q^{-1}}(t)$ is diagonal with $(0,0)$-entry $1$, so
$D$ has $(0,0)$-entry $\alpha$.
Therefore $v_0 = \alpha u_0$.
By Proposition
\ref{prop:comp0d}, both
\begin{eqnarray*}
u_0 = \frac{(\eta_{1},\eta^*_{2})}{(\eta_{3},\eta^*_{2})}
\eta_{3},
\qquad \qquad
v_0 = \frac{(\eta_0,\eta^*_{2})}{(\eta_{3},\eta^*_{2})}
\eta_{3}.
\end{eqnarray*}
Therefore $\alpha = 
(\eta_0,\eta^*_{2})/(\eta_{1},\eta^*_{2})$.
\\
\noindent $\lbrack 1,0,3,2\rbrack \to 
 \lbrack 0,1,3,2\rbrack$.
This transition matrix is the inverse of the transition matrix
for
$\lbrack 0,1,3,2\rbrack \to 
 \lbrack 1,0,3,2\rbrack$.
\\
\noindent $\lbrack 0,2,1,3\rbrack \to 
 \lbrack 2,0,1,3\rbrack$.
Let
$\lbrace u_n \rbrace_{n=0}^d$ and
$\lbrace v_n \rbrace_{n=0}^d$
denote the bases
$\lbrack 0,2,1,3\rbrack$
and
 $\lbrack 2,0,1,3\rbrack$ respectively.
Let $D \in 
{\rm Mat}_{d+1}(\F)$
 denote the transition matrix
from
$\lbrace u_n \rbrace_{n=0}^d$ to
$\lbrace v_n \rbrace_{n=0}^d$.
For $0 \leq n \leq d$ the vectors $u_n$, $v_n$ are contained
in component $n$ of the decomposition $\lbrack 1,3\rbrack$.
Therefore $D$ is diagonal. By Theorem
\ref{thm:xijaction} the matrix representing
$x_{01}$ with respect to $\lbrace u_n \rbrace_{n=0}^d$
is equal to $F_q(t^{-1})$, and the 
 matrix representing
$x_{01}$ with respect to $\lbrace v_n \rbrace_{n=0}^d$
is equal to $E_q$.
Therefore  $F_q(t^{-1}) D = D E_q$. 
Comparing this with
(\ref{eq:FE})
we find that
there exists  $0 \not=\alpha \in \F$ such that
$D=\alpha \mathcal D_q(t)$.
We now 
find $\alpha$. 
The matrix 
$\mathcal D_q(t)$ has $(0,0)$-entry 1,
so 
$D$ has $(0,0)$-entry $\alpha$. 
Therefore $v_0 = \alpha u_0$.
By Proposition
\ref{prop:comp0d}, both
\begin{eqnarray*}
u_0 = \frac{(\eta_{2},\eta^*_{3})}{(\eta_{1},\eta^*_{3})}
\eta_{1},
\qquad \qquad
v_0 = \frac{(\eta_0,\eta^*_{3})}{(\eta_{1},\eta^*_{3})}
\eta_{1}.
\end{eqnarray*}
Therefore $\alpha = 
(\eta_0,\eta^*_{3})/(\eta_{2},\eta^*_{3})$.
\\
\noindent $\lbrack 2,0,1,3\rbrack \to 
 \lbrack 0,2,1,3\rbrack$.
This transition matrix is the inverse of the transition matrix
for
$\lbrack 0,2,1,3\rbrack \to 
 \lbrack 2,0,1,3\rbrack$.
\\
\hfill $\Box$ \\

\noindent We now consider the transitions
of type $\lbrack i,j,k,\ell \rbrack \to 
\lbrack i,k,j,\ell \rbrack$.
We will be discussing the 
matrix $T_q$ defined in
Appendix I. 

\begin{theorem} 
\label{thm:lttrans}
Let $V$ denote an evaluation module for $\boxtimes_q$.
In the table below we display some transition
matrices between the 24 bases for $V$ from
Lemma
\ref{lem:normexist}.
 Each transition matrix
is lower triangular. 
 Pick $r \in \Z_4$.

\medskip
\centerline{
\begin{tabular}[t]{c|c}
{\rm transition} & 
   {\rm transition matrix }
   \\ \hline 
   \hline
\\
$ \lbrack r,r+1,r+2,r+3\rbrack \to
 \lbrack r,r+2,r+1,r+3\rbrack $ &
$
T_q
\frac{(\eta_{r+2},\eta^*_{r+3})}{(\eta_{r+1},\eta^*_{r+3})}
$
\\
\\
$ \lbrack r+1,r,r+2,r+3\rbrack \to
 \lbrack r+1,r+2,r,r+3\rbrack $ &
$
T_q
\frac{(\eta_{r+2},\eta^*_{r+3})}{(\eta_{r},\eta^*_{r+3})}
$
\\
\\
$ \lbrack r,r+1,r+3,r+2\rbrack \to
 \lbrack r,r+3,r+1,r+2\rbrack $ &
$
T_{q^{-1}}
\frac{(\eta_{r+3},\eta^*_{r+2})}{(\eta_{r+1},\eta^*_{r+2})}
$
\\
\\
$ \lbrack r+1,r,r+3,r+2\rbrack \to
 \lbrack r+1,r+3,r,r+2\rbrack $ &
$
T_{q^{-1}}
\frac{(\eta_{r+3},\eta^*_{r+2})}{(\eta_{r},\eta^*_{r+2})}
$
\\
\\
\hline
\\
$ \lbrack r,r+2,r+1,r+3\rbrack \to
 \lbrack r,r+1,r+2,r+3\rbrack $ &
$
T_{q^{-1}}
\frac{(\eta_{r+1},\eta^*_{r+3})}{(\eta_{r+2},\eta^*_{r+3})}
$
\\
\\
$ \lbrack r+2,r,r+1,r+3\rbrack \to
 \lbrack r+2,r+1,r,r+3\rbrack $ &
$
T_q
\frac{(\eta_{r+1},\eta^*_{r+3})}{(\eta_{r},\eta^*_{r+3})}
$
     \end{tabular}}
     \medskip
\end{theorem}
\noindent {\it Proof:}  
These transition matrices were found in
\cite[Theorem~15.4]{fduqe}.
In that article the notation is a bit different
from what we are presently using. To translate between
the notations use Lemma
\ref{lem:comparebasis}
and
\cite[Definition~13.4]{fduqe}.
\hfill $\Box$ \\


\noindent We now consider the transitions of
type
$\lbrack i,j,k,\ell\rbrack  \to 
\lbrack i,j,\ell,k\rbrack$.

\begin{lemma}
Let $V$ denote an evaluation module for $\boxtimes_q$,
and pick
mutually distinct $i,j,k,\ell$ in
$\Z_4$. Then the transition matrix from the basis
$\lbrack i,j,k,\ell\rbrack$ to the basis
$\lbrack i,j,\ell,k\rbrack$ is equal to $Z$.
\end{lemma}
\noindent {\it Proof:}  
By Lemma
\ref{lem:invnorm}.
\hfill $\Box$ \\


\section{Comments on the bilinear form}

\noindent
Throughout this section $V$ denotes
an evaluation module for $\boxtimes_q$, with
diameter $d$ and evaluation parameter $t$.
Recall 
from Definition
\ref{def:cornervector} the vectors
$\lbrace \eta_i \rbrace_{i \in \Z_4}$ in $V$,
and the vectors $\lbrace \eta^*_i\rbrace_{i \in \Z_4}$ in
$V^*$.
We now consider how the scalars
\begin{eqnarray}
(\eta_i,\eta^*_j) \qquad \qquad i,j \in \Z_4, \qquad i\not=j
\label{eq:list}
\end{eqnarray}
are related.

\begin{proposition} 
\label{lem:rel}
With the above notation,
\begin{eqnarray*}
\frac{(\eta_0, \eta^*_1)(\eta_2,\eta^*_3)}
{(\eta_2, \eta^*_1)(\eta_0,\eta^*_3)}
= t^d q^{d(d-1)},
\qquad \qquad 
\frac{(\eta_1, \eta^*_2)(\eta_3,\eta^*_0)}
{(\eta_3, \eta^*_2)(\eta_1,\eta^*_0)}
= t^{-d} q^{d(d-1)}.
\end{eqnarray*}
\end{proposition}
\noindent {\it Proof:}  
Throughout this proof all bases mentioned are for $V$.
Pick $r \in \Z_4$.
Let $\lbrace u_n \rbrace_{n=0}^d$ and 
$\lbrace v_n \rbrace_{n=0}^d$ denote the bases
$\lbrack r,r+2,r+1, r+3\rbrack$ and
$\lbrack r+2,r,r+1, r+3\rbrack$, respectively.
We relate $u_d$, $v_d$  in two ways.
By the fifth row of the table in 
Theorem
\ref{thm:diagtrans},
\begin{eqnarray}
v_d = t^{de} q^{d(d-1)} \frac{(\eta_r, \eta^*_{r+3})}{(\eta_{r+2},\eta^*_{r+3})} u_d,
\label{eq:vd3}
\end{eqnarray}
where $e=(-1)^r$. Observe that $\lbrace u_{d-n}\rbrace_{n=0}^d$
is the basis
$\lbrack r,r+2,r+3, r+1\rbrack$ and
 $\lbrace v_{d-n}\rbrace_{n=0}^d$
is the basis
$\lbrack r+2,r,r+3, r+1\rbrack$.
By the sixth row of the table in Theorem
\ref{thm:diagtrans} (with $r$ replaced by $r+2$),
\begin{eqnarray}
v_d = \frac{(\eta_r, \eta^*_{r+1})}{(\eta_{r+2},\eta^*_{r+1})} u_d.
\label{eq:vd4}
\end{eqnarray}
Comparing 
(\ref{eq:vd3}), 
(\ref{eq:vd4}) we obtain
\begin{eqnarray*}
\frac{(\eta_r, \eta^*_{r+1})}{(\eta_{r+2},\eta^*_{r+1})} 
\,
\frac{(\eta_{r+2}, \eta^*_{r+3})}{(\eta_{r},\eta^*_{r+3})} = 
t^{de} q^{d(d-1)}.
\end{eqnarray*}
The result folllows.
\hfill $\Box$ \\

\begin{proposition}
\label{prop:relate1}
With the above notation,
\begin{eqnarray*}
&&
\frac{(\eta_0, \eta^*_2)(\eta_1,\eta^*_3)}
{(\eta_1, \eta^*_2)(\eta_0,\eta^*_3)}
= (1-tq^{d-1})(1-tq^{d-3})\cdots (1-tq^{1-d}),
\\
&&
\frac{(\eta_1, \eta^*_3)(\eta_2,\eta^*_0)}
{(\eta_2, \eta^*_3)(\eta_1,\eta^*_0)}
= (1-t^{-1}q^{d-1})(1-t^{-1}q^{d-3})\cdots (1-t^{-1}q^{1-d}),
\\
&&
\frac{(\eta_2, \eta^*_0)(\eta_3,\eta^*_1)}
{(\eta_3, \eta^*_0)(\eta_2,\eta^*_1)}
= (1-tq^{d-1})(1-tq^{d-3})\cdots (1-tq^{1-d}),
\\
&&
\frac{(\eta_3, \eta^*_1)(\eta_0,\eta^*_2)}
{(\eta_0, \eta^*_1)(\eta_3,\eta^*_2)}
= (1-t^{-1}q^{d-1})(1-t^{-1}q^{d-3})\cdots (1-t^{-1}q^{1-d}).
\end{eqnarray*}
\end{proposition}
\noindent {\it Proof:}  
Throughout this proof all bases mentioned are for $V$.
Pick $r \in \Z_4$. Let
$\lbrace u_n \rbrace_{n=0}^d$ and 
$\lbrace v_n \rbrace_{n=0}^d$ denote the bases
$\lbrack r,r+1,r+2, r+3\rbrack$ and
$\lbrack r+1,r,r+2, r+3\rbrack$, respectively.
We relate $u_d$, $v_d$  in two ways.
By the first row of the table in 
Theorem
\ref{thm:diagtrans},
\begin{eqnarray}
v_d = (t^e q^{1-d}; q^2)_d \frac{(\eta_r, \eta^*_{r+3})}{(\eta_{r+1},\eta^*_{r+3})} u_d,
\label{eq:vd1}
\end{eqnarray}
where $e=(-1)^r$. Observe that $\lbrace u_{d-n}\rbrace_{n=0}^d$
is the basis
$\lbrack r,r+1,r+3, r+2\rbrack$ and
 $\lbrace v_{d-n}\rbrace_{n=0}^d$
is the basis
$\lbrack r+1,r,r+3, r+2\rbrack$.
By the third row of the table in Theorem
\ref{thm:diagtrans},
\begin{eqnarray}
v_d = \frac{(\eta_r, \eta^*_{r+2})}{(\eta_{r+1},\eta^*_{r+2})} u_d.
\label{eq:vd2}
\end{eqnarray}
Comparing 
(\ref{eq:vd1}), 
(\ref{eq:vd2}) we obtain
\begin{eqnarray*}
\frac{(\eta_r, \eta^*_{r+2})}{(\eta_{r+1},\eta^*_{r+2})} 
\,
\frac{(\eta_{r+1}, \eta^*_{r+3})}{(\eta_{r},\eta^*_{r+3})} = 
(t^e q^{1-d}; q^2)_d.
\end{eqnarray*}
The result folllows.
\hfill $\Box$ \\

\noindent 
We view the following result as a $\boxtimes_q$-analog
of 
\cite[Proposition~13.11]{fduqe}.

\begin{corollary}\label{cor:noT}
With the above notation,
\begin{eqnarray*}
&&
\frac{
(\eta_0, \eta^*_1)
(\eta_1, \eta^*_2)
(\eta_2, \eta^*_0)}
{
(\eta_1, \eta^*_0)
(\eta_2, \eta^*_1)
(\eta_0, \eta^*_2)}
= (-1)^d q^{d(d-1)},
\qquad \quad
\frac{
(\eta_1, \eta^*_2)
(\eta_2, \eta^*_3)
(\eta_3, \eta^*_1)}
{
(\eta_2, \eta^*_1)
(\eta_3, \eta^*_2)
(\eta_1, \eta^*_3)}
= (-1)^d q^{d(d-1)},
\\
&&
\frac{
(\eta_2, \eta^*_3)
(\eta_3, \eta^*_0)
(\eta_0, \eta^*_2)}
{
(\eta_3, \eta^*_2)
(\eta_0, \eta^*_3)
(\eta_2, \eta^*_0)}
= (-1)^d q^{d(d-1)},
\qquad \quad 
\frac{
(\eta_3, \eta^*_0)
(\eta_0, \eta^*_1)
(\eta_1, \eta^*_3)}
{
(\eta_0, \eta^*_3)
(\eta_1, \eta^*_0)
(\eta_3, \eta^*_1)}
= (-1)^d q^{d(d-1)}.
\end{eqnarray*}
\end{corollary}
\noindent {\it Proof:}  We verify the first equation;
the others are similarly verified.
Consider the first two equations in Proposition
\ref{prop:relate1}. The right-hand side of the 
first one is equal to $(-1)^d t^d $ times the right-hand side
of the second one. Therefore, the left-hand side
of the first one is
equal to $(-1)^d t^d $ times the left-hand side
of the second one.
In this equality eliminate $t^d$ using the
first equation in Proposition
\ref{lem:rel}.
\hfill $\Box$ \\

\begin{note}\rm
By Propositions
\ref{lem:rel},
\ref{prop:relate1} 
the scalars (\ref{eq:list})
are determined by the sequence
\begin{eqnarray}
(\eta_0, \eta^*_1), \quad 
(\eta_0, \eta^*_2), \quad 
(\eta_0, \eta^*_3), \quad 
(\eta_1, \eta^*_2),\quad 
(\eta_2, \eta^*_1), \quad 
(\eta_3, \eta^*_0),\quad
(\eta_3, \eta^*_1).
\label{eq:free2}
\end{eqnarray}
The scalars 
(\ref{eq:free2})
are ``free'' in the following sense.
Given a sequence  $\theta$ 
of seven nonzero scalars in $\F$,
there exist vectors $\eta_i, \eta^*_i$ 
$(i \in \Z_4)$
 as in Definition
\ref{def:cornervector}
such that the sequence 
(\ref{eq:free2})
is equal to $\theta$.
\end{note}


\section{Exchangers}

\noindent
Throughout this section $V$ denotes an evaluation module
for $\boxtimes_q$, with diameter $d$ and evaluation parameter $t$.
We investigate a type of map
in ${\rm{End}}(V)$ called an exchanger.

\begin{lemma}
\label{lem:preex}
For each standard generator $x_{ij}$ of 
$\boxtimes_q$ and 
each  $X \in {\rm{End}}(V)$, 
the following are equivalent:
\begin{enumerate}
\item[\rm (i)] $X$ is invertible and
 $X x_{ij} X^{-1} = x_{i+2,j+2}$ holds on $V$;
\item[\rm (ii)] $X$ sends the decomposition
$\lbrack i,j\rbrack $ of $V$ to the decomposition
$\lbrack i+2,j+2\rbrack $ of $V$.
\end{enumerate}
\end{lemma}
\noindent {\it Proof:} 
Recall from Section 6 that for $0 \leq n \leq d$, component
$n$ of the decomposition $\lbrack i,j\rbrack$ (resp.
$\lbrack i+2,j+2\rbrack$) of $V$
is the eigenspace of $x_{ij}$
(resp. $x_{i+2,j+2}$) with eigenvalue
$q^{d-2n}$. The result follows from this and linear algebra.
\hfill $\Box$ \\

\begin{definition}
\label{def:exch}
\rm
By an {\it exchanger} for $V$ we mean an invertible
$X \in 
{\rm End}(V)$ such that on $V$ the equation
\begin{eqnarray*}
X x_{ij} X^{-1} = x_{i+2,j+2}
\end{eqnarray*}
holds for each standard generator $x_{ij}$ of $\boxtimes_q$.
\end{definition}

\noindent We now consider the existence and uniqueness of
the exchangers. We start with existence.

\begin{lemma} 
\label{lem:exchmeaning}
For 
$X \in 
{\rm End}(V)$
the following are equivalent:
\begin{enumerate}
\item[\rm (i)] $X$ is an exchanger for $V$;
\item[\rm (ii)] $X$ is an isomorphism of $\boxtimes_q$-modules
from $V$ to $V$ twisted via $\rho^2$.
\end{enumerate}
\end{lemma}
\noindent {\it Proof:}  
Recall that $\rho^2$ sends each standard generator 
$x_{ij}\mapsto x_{i+2,j+2}$.
\hfill $\Box$ \\

\begin{lemma} 
\label{lem:exexist}
There exists an exchanger for $V$.
\end{lemma}
\noindent {\it Proof:}  
By Corollary
\ref{cor:ttinv} and
Lemma
\ref{lem:exchmeaning}.
\hfill $\Box$ \\

\noindent We now consider the uniqueness of the exchangers.

\begin{lemma} 
\label{lem:exchunique}
Let $\Psi$ denote an
 exchanger for $V$.
Then for all $X \in 
{\rm End}(V)$ the following are equivalent:
\begin{enumerate}
\item[\rm (i)] $X$ is an exchanger for $V$;
\item[\rm (ii)] there exists $0 \not=\alpha \in \F$
such that 
$X = \alpha \Psi $.
\end{enumerate}
\end{lemma}
\noindent {\it Proof:}  
${\rm (i)}\Rightarrow {\rm (ii)}$
 The composition $G=X\Psi^{-1}$ commutes with each
 standard generator of $\boxtimes_q$, and therefore
 everything in $\boxtimes_q$.
Since $\F$ is algebraically closed and 
$V$ has finite positive dimension,
there exists an eigenspace 
$W \subseteq V$ for $G$. 
 Let $\alpha \in \F$
 denote the corresponding eigenvalue. Then $\alpha\not=0$ since
 $G$ is invertible.
Since $G$ commutes with everything in $\boxtimes_q$,
 we see that $W$ is a $\boxtimes_q$-submodule of $V$. The $\boxtimes_q$-module
 $V$ is irreducible so $W=V$. Therefore $G=\alpha I$ so 
 $X = \alpha \Psi$.
\\
\noindent 
${\rm (ii)}\Rightarrow {\rm (i)}$ Routine.
\hfill $\Box$ \\

\noindent We now characterize the exchangers
in various ways.

\begin{lemma} 
\label{lem:various}
For $X \in {\rm End}(V)$ the following are equivalent:
\begin{enumerate}
\item[\rm (i)]
$X$ is an exchanger for $V$;
\item[\rm (ii)]
for all $i \in \Z_4$,
$X$ sends the flag $\lbrack i\rbrack $ for 
$V$ to the
flag
$\lbrack i+2\rbrack $ for $V$;
\item[\rm (iii)]
for all distinct $i,j \in \Z_4$,
$X$ sends the decomposition $\lbrack i,j\rbrack $ of
$V$ to the
decomposition $\lbrack i+2,j+2\rbrack $ of $V$.
\end{enumerate}
\end{lemma}
\noindent {\it Proof:}  
${\rm (i)}\Rightarrow {\rm (ii)}$ 
Let $i \in \Z_4$ be given. By
Lemma
\ref{lem:exchmeaning}, $X$ sends the flag 
$\lbrack i\rbrack $ for 
$V$ to the
flag 
$\lbrack i\rbrack $ 
for $V$ twisted via $\rho^2$.
By 
Lemma \ref{lem:twistflag},
the flag
$\lbrack i\rbrack $ 
for $V$ twisted via $\rho^2$ is the same thing
as the flag $\lbrack i+2\rbrack$ for $V$.
The result follows.
\\
\noindent 
${\rm (ii)}\Rightarrow {\rm (iii)}$ 
For the flags and decompositions under discussion,
their relationship is described near the end of
 Section 6.
\\
\noindent
${\rm (iii)}\Rightarrow {\rm (i)}$ 
By Lemma
\ref{lem:preex} and
Definition
\ref{def:exch}.
\hfill $\Box$ \\

\begin{lemma}
\label{lem:exchar2}
For $X \in {\rm End}(V)$ the following are equivalent:
\begin{enumerate}
\item[\rm (i)]
$X$ is an exchanger for $V$;
\item[\rm (ii)]
for all mutually distinct $i,j,k,\ell$ in $\Z_4$,
$X$ sends each $\lbrack i,j,k,\ell\rbrack$-basis
for $V$ to an 
$\lbrack i+2,j+2,k+2,\ell+2\rbrack$-basis for $V$;
\item[\rm (iii)]
there exist mutually distinct $i,j,k,\ell$ in $\Z_4$
such that $X$ sends the basis $\lbrack i,j,k,\ell\rbrack$
for $V$ to an 
$\lbrack i+2,j+2,k+2,\ell+2\rbrack$-basis for $V$.
\end{enumerate}
\end{lemma}
\noindent {\it Proof:}  
${\rm (i)}\Rightarrow {\rm (ii)}$ 
Let 
 mutually distinct $i,j,k,\ell$ in $\Z_4$ be given.
By 
Lemma
\ref{lem:exchmeaning}, $X$ sends
 each $\lbrack i,j,k,\ell\rbrack$-basis
for $V$ to an 
$\lbrack i,j,k,\ell\rbrack$-basis for 
$V$ twisted via $\rho^2$.
By
Lemma \ref{lem:rhoeffect}, 
an $\lbrack i,j,k,\ell\rbrack$-basis for 
$V$ twisted via $\rho^2$ is the same thing as an
$\lbrack i+2,j+2,k+2,\ell+2\rbrack$-basis for 
$V$.
The result follows.
\\
\noindent
${\rm (ii)}\Rightarrow {\rm (iii)}$ 
Clear.
\\
\noindent
${\rm (iii)}\Rightarrow {\rm (i)}$ 
By Lemma
\ref{lem:exexist} there exists an exchanger
$\Psi$ for $V$.
By the implication 
${\rm (i)}\Rightarrow {\rm (ii)}$ above,
$\Psi$ sends the basis
 $\lbrack i,j,k,\ell\rbrack$
of $V$ to an 
$\lbrack i+2,j+2,k+2,\ell+2\rbrack$-basis of $V$.
By assumption, 
$X$ also sends the basis
 $\lbrack i,j,k,\ell\rbrack$
of $V$ to an 
$\lbrack i+2,j+2,k+2,\ell+2\rbrack$-basis of $V$.
Now by Lemma
\ref{lem:basisprime},
there exists $0 \not=\alpha \in \F$ such
that $X=\alpha \Psi$.
Now  by Lemma
\ref{lem:exchunique},
$X$ is an exchanger for $V$.
\hfill $\Box$ \\

\begin{lemma} For 
 $X \in {\rm End}(V)$ the following are equivalent:
\begin{enumerate}
\item[\rm (i)]
$X$ is an exchanger for $V$;
\item[\rm (ii)]
$X$ is an exchanger for 
${}^\rho V$.
\end{enumerate}
\end{lemma}
\noindent {\it Proof:}  
Use
Lemma
\ref{lem:twistflag}
and Lemma
\ref{lem:various}(i),(ii).
\hfill $\Box$ \\

\begin{lemma} Let $X$ denote an exchanger for $V$. Then
$X^{-1}$ is an exchanger for $V$.
\end{lemma} 
\noindent {\it Proof:}  
Use Lemma
\ref{lem:various}(i),(ii).
\hfill $\Box$ \\

\noindent We will return to exchangers shortly.

\begin{lemma} 
\label{lem:15pt5}
Consider the
 24 bases for $V$ from
Lemma
\ref{lem:normexist}.
In the table below we display some transition
matrices between 
these bases.
For each transition matrix 
we give the $(i,j)$-entry for $0 \leq i,j \leq d$.
Pick $r \in \Z_4$,
and first assume that $r$ is even.

\medskip
\centerline{
\begin{tabular}[t]{c|c}
{\rm transition matrix} & 
   {\rm $(i,j)$-entry}
   \\ \hline
\\
$ \lbrack r,r+2,r+1,r+3\rbrack \to
 \lbrack r+2,r,r+3,r+1\rbrack $ &
$
\delta_{i+j,d}t^iq^{i(d-1)} 
\frac{(\eta_r,\eta^*_{r+3})}{(\eta_{r+2},\eta^*_{r+3})}
$
\\
\\
$ \lbrack r+2,r,r+1,r+3\rbrack \to
 \lbrack r,r+2,r+3,r+1\rbrack $ &
$
\delta_{i+j,d}t^{-i}q^{i(1-d)}
\frac{
(\eta_{r+2},\eta^*_{r+3})
}{
(\eta_r,\eta^*_{r+3})
}
$
     \end{tabular}}
     \medskip
\noindent Next assume that $r$ is odd. Then in the above
table replace $t$ by $t^{-1}$.
\end{lemma}
\noindent {\it Proof:}  
To find the transition matrix 
$\lbrack r,r+2,r+1,r+3\rbrack\to
\lbrack r+2,r,r+3,r+1\rbrack$,
compute the product of transition matrices
\begin{eqnarray*}
\lbrack r,r+2,r+1,r+3\rbrack \to
\lbrack r+2,r,r+1,r+3\rbrack \to
\lbrack r+2,r,r+3,r+1\rbrack.
\end{eqnarray*}
In this product 
the first transition matrix is from Theorem
\ref{thm:diagtrans} and the second one is equal to $Z$.
To find the transition matrix
$\lbrack r+2,r,r+1,r+3\rbrack\to
\lbrack r,r+2,r+3,r+1\rbrack$,
compute the product of transition matrices
\begin{eqnarray*}
\lbrack r+2,r,r+1,r+3\rbrack \to
\lbrack r,r+2,r+1,r+3\rbrack \to
\lbrack r,r+2,r+3,r+1\rbrack.
\end{eqnarray*}
In this product 
the first transition matrix is from Theorem
\ref{thm:diagtrans} and the second one is equal to $Z$.
\hfill $\Box$ \\

\begin{theorem}
\label{thm:exchrep}
There exists an exchanger $ \mathcal X$
for 
$V$ that is described as follows.
In the table below,
each row contains a basis for $V$ from
Lemma
\ref{lem:normexist},
 and the entries of 
a matrix in 
${\rm Mat}_{d+1}(\F)$.
The matrix represents $\mathcal X$ with respect to the basis.

\medskip
\centerline{
\begin{tabular}[t]{c|c}
{\rm basis} & 
   {\rm $(i,j)$-entry for $0 \leq i,j\leq d$}  
   \\ \hline \hline
 $\lbrack 0,2,1,3 \rbrack $   
& $\delta_{i+j,d}t^{i} q^{i(d-1)-\binom{d}{2}}$
   \\ 
    $\lbrack 0,2,3,1 \rbrack $   
& $\delta_{i+j,d}t^{d-i} q^{\binom{d}{2}-i(d-1)}$
\\
    $\lbrack 2,0,3,1 \rbrack $   
& $\delta_{i+j,d}t^{i} q^{i(d-1)-\binom{d}{2}}$
   \\ 
    $\lbrack 2,0,1,3 \rbrack $   
& $\delta_{i+j,d}t^{d-i} q^{\binom{d}{2}-i(d-1)}$
\\
\hline
    $\lbrack 1,3,2,0 \rbrack $   
   & $\delta_{i+j,d}(-1)^d t^{d-i} q^{i(d-1)-\binom{d}{2}}$
   \\
    $\lbrack 1,3,0,2 \rbrack $   
  & $\delta_{i+j,d}(-1)^d t^i q^{\binom{d}{2}-i(d-1)}$
   \\
    $\lbrack 3,1,0,2 \rbrack $   
   & $\delta_{i+j,d}(-1)^d t^{d-i} q^{i(d-1)-\binom{d}{2}}$
   \\
   $\lbrack 3,1,2,0 \rbrack $   
  & $\delta_{i+j,d}(-1)^d t^i q^{\binom{d}{2}-i(d-1)}$
\end{tabular}}
     \medskip
\end{theorem}
\noindent {\it Proof:}  
The above table has 8 rows. For 
$1 \leq h \leq 8$, 
row $h$ of the table 
contains a basis and the entries of a matrix
which we call $X_h$.
Define $\mathcal X_h\in
 {\rm End}(V)$ 
such that $\mathcal X_h$ is represented by $X_h$ with respect to
the basis in row $h$. 
We claim that $\mathcal X_h$ is an exchanger for $V$.
To prove the claim, first assume that $h=1$.
Compare
 row $1$ of the above table
with row $1$ of the table in Lemma
\ref{lem:15pt5} (with $r=0$). 
The comparison shows that 
the matrix $X_1$ is a scalar multiple of
the transition matrix 
from the basis
$\lbrack 0,2,1,3\rbrack$ to
the basis $\lbrack
2,0,3,1\rbrack$.
 Therefore
 $\mathcal X_1$ sends the basis
$\lbrack 0,2,1,3\rbrack$ to a
$\lbrack 2,0,3,1\rbrack$-basis.
Now 
 $\mathcal X_1$ is an exchanger for $V$, in view of
 Lemma
\ref{lem:exchar2}(i),(iii).
The claim is proven for $h=1$, and for $2 \leq h \leq 8$
the argument is similar. 
We now show that $\mathcal X_h$ is
independent of $h$ for $1 \leq h \leq 8$.
\\
\noindent
We show $\mathcal X_1 = \mathcal X_2$.
By construction the matrix
$X_1$
represents $\mathcal X_1$
with respect to
$\lbrack 0,2,1,3\rbrack $.
The transition matrix from
$\lbrack 0,2,1,3\rbrack $ to
$ \lbrack 0,2,3,1\rbrack$ is equal to $Z$.
Therefore the matrix $Z X_1 Z$ represents 
$\mathcal X_1$ with respect to 
 $\lbrack 0,2,3,1\rbrack$.
By construction $X_2$ represents
$\mathcal X_2$ with respect to
 $\lbrack 0,2,3,1\rbrack$.
One checks $ X_1 Z = Z X_2$ so
$ ZX_1 Z = X_2$.
Therefore
$\mathcal X_1 = \mathcal X_2$.
By a similar argument
\begin{eqnarray*}
\mathcal X_3 = \mathcal X_4,
\qquad 
\mathcal X_5 = \mathcal X_6,
\qquad 
\mathcal X_7 = \mathcal X_8.
\end{eqnarray*}
We show $\mathcal X_2 = \mathcal X_3$.
By Theorem
\ref{thm:diagtrans},
the transition matrix from
$\lbrack 2,0,3,1\rbrack$ to
$\lbrack 0,2,3,1\rbrack$ is a scalar multiple of
$\mathcal D_q(t)$. One checks
$X_3 \mathcal D_q(t) = \mathcal D_q(t)X_2$.
Therefore
$\mathcal X_2 = \mathcal X_3$.
We show $\mathcal X_5 = \mathcal X_8$.
By Theorem
\ref{thm:diagtrans},
the transition matrix 
from $\lbrack 1,3,2,0\rbrack$ to
$\lbrack 3,1,2,0\rbrack$ is a scalar multiple of
$\mathcal D_q(t^{-1})$. One checks
$X_5 \mathcal D_q(t^{-1}) = \mathcal D_q(t^{-1})X_8$.
Therefore 
$\mathcal X_5 = \mathcal X_8$.
We  show 
$\mathcal X_1 = \mathcal X_6$.
Since each of
$\mathcal X_1$, $\mathcal X_6$ is an exchanger, 
By Lemma
\ref{lem:exchunique} there exists 
$0 \not=\alpha \in \F$ such that
$\mathcal X_1 = \alpha \mathcal X_6$.
We show $\alpha = 1$. Let $T$ denote the transition
matrix from
$\lbrack 0,2,1,3\rbrack$ to
$\lbrack 1,3,0,2\rbrack$. By construction
$X_1 T = \alpha T X_6$.  We now find $T$.
To this end, compute the product of transition matrices
\begin{eqnarray*}
\lbrack 0,2,1,3\rbrack
\to
\lbrack 0,1,2,3\rbrack
\to
\lbrack 1,0,2,3\rbrack
\to
\lbrack 1,0,3,2\rbrack
\to
\lbrack 1,3,0,2\rbrack.
\end{eqnarray*}
In this product each factor is given in
Section 13.
The computation shows that $T$ is a scalar multiple of
$T_{q^{-1}} D_q(t) Z T_{q^{-1}}$. Therefore
\begin{eqnarray}
X_1 
T_{q^{-1}} D_q(t) Z T_{q^{-1}}
=
\alpha
T_{q^{-1}} D_q(t) Z T_{q^{-1}} X_6.
\label{eq:x15}
\end{eqnarray}
Using the fact that $T_{q^{-1}}$ is lower triangular
and $D_q(t)$ is diagonal, we routinely compute the
$(0,d)$-entry for each side of
(\ref{eq:x15}). Comparing these entries we find
$\alpha = 1$. Therefore
$\mathcal X_1 = \mathcal X_6$.
We have shown that $\mathcal X_h$ is independent of
$h$ for $1 \leq h \leq 8$. The result follows.
\hfill $\Box$ \\

\begin{definition}\rm
By the {\it standard exchanger} for $V$,
we mean the 
exchanger $\mathcal X$ 
from Theorem
\ref{thm:exchrep}.
\end{definition}

\begin{theorem}
\label{thm:exchinv}
For the standard exchanger $\mathcal X$ of $V$,
\begin{eqnarray*}
\mathcal X^2 = t^d I.
\end{eqnarray*}
\end{theorem}
\noindent {\it Proof:} 
Let $M$ denote the matrix 
that represents $\mathcal X$ with respect to
the basis $\lbrack 0,2,1,3\rbrack$. The matrix
$M$ is given in the first row of the table in Theorem
\ref{thm:exchrep}. By matrix multiplication
$M^2 = t^d I.$ The result follows.
\hfill $\Box$ \\

\begin{proposition} 
\label{prop:stexchinv}
The following coincide:
\begin{enumerate}
\item[\rm (i)]
the standard exchanger for 
${}^\rho V$;
\item[\rm (ii)]
$(-1)^d$ times the inverse of the standard exchanger for $V$.
\end{enumerate}
\end{proposition}
\noindent {\it Proof:}
Referring to Theorem
\ref{thm:exchrep},
compare the matrices that represent
$\mathcal X$ with respect to the bases
$\lbrack 0,2,1,3\rbrack$ and
$\lbrack 1,3,2,0\rbrack$.
Recall from Lemma
\ref{lem:rhoeffect2}
that $t^{-1}$ is the evaluation parameter for
${}^\rho V$. Also, by Theorem 
\ref{thm:exchinv} $\mathcal X^{-1} = t^{-d}\mathcal X$.
\hfill $\Box$ \\

\begin{corollary}
The following coincide:
\begin{enumerate}
\item[\rm (i)] the standard exchanger for $V$ twisted via $\rho^2$;
\item[\rm (ii)] the standard exchanger for $V$.
\end{enumerate}
\end{corollary}
\noindent {\it Proof:} 
Apply Proposition \ref{prop:stexchinv} twice.
\hfill $\Box$ \\

\noindent 
Let $\mathcal X$ denote the standard exchanger for $V$.
We now describe what $\mathcal X$ does to the
vectors $\lbrace \eta_i\rbrace_{i \in \Z_4}$
from Definition
\ref{def:cornervector}.

\begin{proposition}
\label{prop:etaMove}
Let $\mathcal X$ denote the standard exchanger for $V$. Then
\begin{eqnarray*}
&&
\mathcal X \eta_0 =  q^{-\binom{d}{2}} 
\frac{(\eta_0,\eta^*_1)}{(\eta_2, \eta^*_1)}\eta_2,
\qquad \qquad 
\mathcal X \eta_1 =  (-1)^d q^{\binom{d}{2}} 
\frac{(\eta_1,\eta^*_0)}{(\eta_3, \eta^*_0)}\eta_3,
\\
&&
\mathcal X \eta_2 =  q^{-\binom{d}{2}} 
\frac{(\eta_2,\eta^*_3)}{(\eta_0, \eta^*_3)}\eta_0,
\qquad \qquad 
\mathcal X \eta_3 =  (-1)^d q^{\binom{d}{2}} 
\frac{(\eta_3,\eta^*_2)}{(\eta_1, \eta^*_2)}\eta_1.
\end{eqnarray*}
\end{proposition}
\noindent {\it Proof:} 
We verify the first equation.
Let $\lbrace u_n\rbrace_{n=0}^d$ 
and $\lbrace v_n\rbrace_{n=0}^d$ 
denote the bases $\lbrack 2,0,3,1\rbrack$
and
 $\lbrack 0,2,1,3\rbrack$ for $V$, respectively.
Let $T$
denote 
the transition matrix from 
$\lbrace u_n\rbrace_{n=0}^d$ 
to $\lbrace v_n\rbrace_{n=0}^d$;  this matrix is given in Lemma
\ref{lem:15pt5}.
Let $M$ denote the matrix that 
represents $\mathcal X$ with respect to
$\lbrace u_n\rbrace_{n=0}^d$; this matrix is given in
Theorem
\ref{thm:exchrep}. Comparing $M$ and $T$ we obtain
\begin{eqnarray*}
M = \alpha T, \qquad \qquad 
\alpha = q^{-\binom{d}{2}} \frac{(\eta_0,\eta^*_1)}{(\eta_2,\eta^*_1)}.
\end{eqnarray*}
Therefore $\mathcal X u_n = \alpha v_n$ for $0 \leq n \leq d$.
By
Lemma 
\ref{lem:normexist}(ii),
$\eta_0 = \sum_{n=0}^d u_n$ and
$\eta_2 = \sum_{n=0}^d v_n$. Consequently
$\mathcal X \eta_0 = \alpha \eta_2$. The first equation is verified.
The remaining equations are similarly verified.
\hfill $\Box$ \\

\noindent 
Let $\mathcal X$ denote the standard exchanger for $V$.
We now describe what $\mathcal X$ does to the 24 bases
for $V$ from Lemma
\ref{lem:normexist}.

\begin{theorem} 
Let $\mathcal X$ denote the standard exchanger
for $V$.
For mutually distinct $i,j,k,\ell$ in $\Z_4$,
consider the bases 
$\lbrack i,j,k,\ell\rbrack $ 
and $\lbrack i+2,j+2,k+2,\ell+2\rbrack $ of $V$.
The map $\mathcal X$  sends
$\lbrack i,j,k,\ell\rbrack $ 
to a scalar multiple of
$\lbrack i+2,j+2,k+2,\ell+2\rbrack $.
The scalar is given in the tables below.

\bigskip

\centerline{
\begin{tabular}[t]{cccc|c}
$i$ & $j$ & $k$ & $\ell$ & 
   {\rm scalar }  
   \\ \hline \hline
 &&&&\\
 $0$ & $2$ & $1$ & $3$
& 
 $
 q^{-\binom{d}{2}}
 \frac{(\eta_2, \eta^*_3)}{(\eta_0, \eta^*_3)}
 $   
   \\ 
  &&&& \\
 $0$ & $2$ & $3$ &$1$
& 
 $
 q^{-\binom{d}{2}}
 \frac{(\eta_2, \eta^*_3)}{(\eta_0, \eta^*_3)}
    $   
\\
&&&&
\\
$2$ & $0$ & $3$ &$1$
& 
    $
 q^{-\binom{d}{2}}
 \frac{(\eta_0, \eta^*_1)}{(\eta_2, \eta^*_1)}
    $   
   \\ 
&&&&
\\
$2$ & $0$ & $1$ &$3$
& 
 $
 q^{-\binom{d}{2}}
 \frac{(\eta_0, \eta^*_1)}{(\eta_2, \eta^*_1)} $
\\
\end{tabular}
\qquad \qquad 
\begin{tabular}[t]{cccc|c}
$i$ & $j$ & $k$ & $\ell$ & 
   {\rm scalar }  
   \\ \hline \hline
&&&&\\
$1$ & $3$ & $2$ &$0$
& 
    $
 (-1)^d q^{\binom{d}{2}}\frac{(\eta_3, \eta^*_2)}{(\eta_1, \eta^*_2)}
    $   
\\
&&&&
\\
$1$ & $3$ & $0$ &$2$
& 
    $
 (-1)^d q^{\binom{d}{2}}\frac{(\eta_3, \eta^*_2)}{(\eta_1, \eta^*_2)}
$
\\
&&&&\\
$3$ & $1$ & $0$ &$2$
& 
$
 (-1)^d q^{\binom{d}{2}}\frac{(\eta_1, \eta^*_0)}{(\eta_3, \eta^*_0)}
$
\\
&&&&\\
$3$ & $1$ & $2$ &$0$
& 
    $
 (-1)^d q^{\binom{d}{2}}\frac{(\eta_1, \eta^*_0)}{(\eta_3, \eta^*_0)}
    $ 
\\
\end{tabular}
}

\bigskip

\centerline{
\begin{tabular}[t]{cccc|c}
$i$ & $j$ & $k$ & $\ell$ & 
   {\rm scalar }  
   \\ \hline \hline
&&&&\\
$0$ & $1$ & $2$ & $3$
& 
 $
(-1)^d q^{\binom{d}{2}}
 \frac{(\eta_1, \eta^*_0)}{
 (\eta_3, \eta^*_0)}
$
   \\ 
&&&&\\
$0$ & $1$ & $3$ & $2$
& 
 $
(-1)^d q^{\binom{d}{2}}
 \frac{(\eta_1, \eta^*_0)}{
 (\eta_3, \eta^*_0)}
$
   \\ 
&&&&\\
$1$ & $0$ & $3$ & $2$
& 
 $
q^{-\binom{d}{2}}
 \frac{(\eta_0, \eta^*_1)}{
 (\eta_2, \eta^*_1)}
$
   \\ 
&&&&\\
$1$ & $0$ & $2$ & $3$
& 
 $
q^{-\binom{d}{2}}
 \frac{(\eta_0, \eta^*_1)}{
 (\eta_2, \eta^*_1)}
$
   \\ 
\end{tabular}
\qquad \qquad 
\begin{tabular}[t]{cccc|c}
$i$ & $j$ & $k$ & $\ell$ & 
   {\rm scalar }  
   \\ \hline \hline
&&&&\\
$1$ & $2$ & $3$ & $0$
& 
 $
q^{-\binom{d}{2}}
 \frac{(\eta_2, \eta^*_3)}{
 (\eta_0, \eta^*_3)}
$
   \\ 
&&&&\\
$1$ & $2$ & $0$ & $3$
& 
 $
q^{-\binom{d}{2}}
 \frac{(\eta_2, \eta^*_3)}{
 (\eta_0, \eta^*_3)}
$
   \\ 
&&&&\\
$2$ & $1$ & $0$ & $3$
& 
 $
(-1)^d q^{\binom{d}{2}}
 \frac{(\eta_1, \eta^*_0)}{
 (\eta_3, \eta^*_0)}
$
   \\ 
&&&&\\
$2$ & $1$ & $3$ & $0$
& 
 $
(-1)^d q^{\binom{d}{2}}
 \frac{(\eta_1, \eta^*_0)}{
 (\eta_3, \eta^*_0)}
$
   \\ 
\end{tabular}
}

\bigskip

\centerline{
\begin{tabular}[t]{cccc|c}
$i$ & $j$ & $k$ & $\ell$ & 
   {\rm scalar }  
   \\ \hline \hline
&&&&\\
$2$ & $3$ & $0$ & $1$
& 
 $
(-1)^d q^{\binom{d}{2}}
 \frac{(\eta_3, \eta^*_2)}{
 (\eta_1, \eta^*_2)}
$
   \\ 
&&&&\\
$2$ & $3$ & $1$ & $0$
& 
 $
(-1)^d q^{\binom{d}{2}}
 \frac{(\eta_3, \eta^*_2)}{
 (\eta_1, \eta^*_2)}
$
   \\ 
&&&&\\
$3$ & $2$ & $1$ & $0$
& 
 $
q^{-\binom{d}{2}}
 \frac{(\eta_2, \eta^*_3)}{
 (\eta_0, \eta^*_3)}
$
   \\ 
&&&&\\
$3$ & $2$ & $0$ & $1$
& 
 $
q^{-\binom{d}{2}}
 \frac{(\eta_2, \eta^*_3)}{
 (\eta_0, \eta^*_3)}
$
   \\ 
\end{tabular}
\qquad \qquad 
\begin{tabular}[t]{cccc|c}
$i$ & $j$ & $k$ & $\ell$ & 
   {\rm scalar }  
   \\ \hline \hline
&&&&\\
$3$ & $0$ & $1$ & $2$
& 
 $
q^{-\binom{d}{2}}
 \frac{(\eta_0, \eta^*_1)}{
 (\eta_2, \eta^*_1)}
$
   \\ 
&&&&\\
$3$ & $0$ & $2$ & $1$
& 
 $
q^{-\binom{d}{2}}
 \frac{(\eta_0, \eta^*_1)}{
 (\eta_2, \eta^*_1)}
$
   \\ 
&&&&\\
$0$ & $3$ & $2$ & $1$
& 
 $
(-1)^d q^{\binom{d}{2}}
 \frac{(\eta_3, \eta^*_2)}{
 (\eta_1, \eta^*_2)}
$
   \\ 
&&&&\\
 $0$ & $3$ & $1$ & $2$
& 
 $
(-1)^d q^{\binom{d}{2}}
 \frac{(\eta_3, \eta^*_2)}{
 (\eta_1, \eta^*_2)}
$
\\
\end{tabular}
}

\bigskip

\end{theorem}
\noindent {\it Proof:} 
Let $\lbrace u_n\rbrace_{n=0}^d$ 
and
$\lbrace v_n\rbrace_{n=0}^d$ 
denote the bases
$\lbrack i,j,k,\ell \rbrack$ and
$\lbrack i+2,j+2,k+2,\ell+2 \rbrack$, respectively.
By Lemma
\ref{lem:basisprime}
and Lemma \ref{lem:exchar2}, 
there exists $0 \not=\alpha \in \F$ such that
$\mathcal X u_n = \alpha v_n$ for $0 \leq n \leq d$.
By Lemma
\ref{lem:normexist}(ii), 
$\eta_j = \sum_{n=0}^d u_n $ and
$\eta_{j+2} = \sum_{n=0}^d v_n $.
Therefore 
$\mathcal X \eta_j = \alpha \eta_{j+2}$.
Evaluating this using Proposition
\ref{prop:etaMove} we find that $\alpha$ is as shown
 in the tables.
\hfill $\Box$ \\


\section{Leonard pairs of $q$-Racah type}

\noindent We now turn our attention 
to Leonard pairs. There is a general family 
of Leonard pairs said to have $q$-Racah type
\cite[Section~5]{hwhuang},
\cite[Example~5.3]{terPA}.
In this section, we show that
 for any Leonard pair of $q$-Racah type,
 the underlying vector space becomes an evaluation
 module for 
 $\boxtimes_q$ in a natural way.
\medskip

\noindent 
Let $a,b,c $ denote nonzero scalars in $\F$. Recall the
equitable generators $x,y,z$ for 
$U_q(\mathfrak{sl}_2)$.

\begin{definition}
\label{def:ABC}
\rm
Using $a,b,c$ 
we define
some elements in
$U_q(\mathfrak{sl}_2)$:
\begin{eqnarray*}
&&{\bf A} = ax+ a^{-1} y + bc^{-1}\frac{xy-yx}{q-q^{-1}},
\\
&&{\bf B} = by+ b^{-1} z + ca^{-1}\frac{yz-zy}{q-q^{-1}},
\\
&&{\bf C} = cz+ c^{-1} x + ab^{-1}\frac{zx-xz}{q-q^{-1}}.
\end{eqnarray*}
The sequence $\bf A, B, C$ is called the {\it Askey-Wilson
triple} for $a,b,c$.
\end{definition}

\begin{lemma} 
{\rm \cite[Proposition~1.1]{uawe}}.
\label{lem:awZ3}
We have
\begin{eqnarray*}
&&
{\bf A} + \frac{q{\bf BC}-q^{-1}{\bf CB}}{q^2-q^{-2}} = \frac{(a+a^{-1})\Lambda
+(b+b^{-1})(c+c^{-1})}{q+q^{-1}},
\\
&&
{\bf B} + \frac{q{\bf CA}-q^{-1}{\bf AC}}{q^2-q^{-2}} = \frac{(b+b^{-1}) \Lambda
+(c+c^{-1})(a+a^{-1})}{q+q^{-1}},
\\
&&
{\bf C} + \frac{q{\bf AB}-q^{-1}{\bf BA}}{q^2-q^{-2}} = \frac{(c+c^{-1}) \Lambda
+(a+a^{-1})(b+b^{-1})}{q+q^{-1}}.
\end{eqnarray*}
\end{lemma}

\begin{note}\rm The equations in
Lemma
\ref{lem:awZ3} are a variation on the $\Z_3$-symmetric
Askey-Wilson relations
\cite[Theorem~10.1]{hwhuang},
\cite[Section 1]{daha}.
\end{note}

\begin{definition}
\label{def:ABCp}
\rm
Using $a,b,c$ we define some more elements in
$U_q(\mathfrak{sl}_2)$:
\begin{eqnarray*}
&&{\bf A'} = ay+ a^{-1} z + cb^{-1}\frac{yz-zy}{q-q^{-1}},
\\
&&{\bf B'} = bx+ b^{-1} y + ac^{-1}\frac{xy-yx}{q-q^{-1}},
\\
&&{\bf C'} = cz+ c^{-1} x +ba^{-1}\frac{zx-xz}{q-q^{-1}}.
\end{eqnarray*}
The sequence $\bf A',B', C'$ is called the
{\it dual Askey-Wilson triple}
for $a,b,c$.
\end{definition}

\begin{lemma} 
\label{lem:ABCvsABCp}
Let $\bf A',B',C'$ denote the
dual Askey-Wilson triple for
$a,b,c$. Then $\bf B',A', C'$ is the Askey-Wilson
triple for $b,a,c$.
\end{lemma}
\noindent {\it Proof:}  Compare
the equations in
Definition
\ref{def:ABC} and
Definition
\ref{def:ABCp}.
\hfill $\Box$ \\

\begin{lemma} 
\label{lem:awp}
Let 
$\bf A',B',C'$ denote the dual Askey-Wilson triple 
for
$a,b,c$.
Then
\begin{eqnarray*}
&&
{\bf A'} + \frac{q{\bf C'B'}-q^{-1}{\bf B'C'}}{q^2-q^{-2}} = \frac{(a+a^{-1})\Lambda
+(b+b^{-1})(c+c^{-1})}{q+q^{-1}},
\\
&&
{\bf B'} + \frac{q{\bf A'C'}-q^{-1}{\bf C'A'}}{q^2-q^{-2}} = \frac{(b+b^{-1}) \Lambda
+(c+c^{-1})(a+a^{-1})}{q+q^{-1}},
\\
&&
{\bf C'} + \frac{q{\bf B'A'}-q^{-1}{\bf A'B'
}}{q^2-q^{-2}} = \frac{(c+c^{-1}) \Lambda
+(a+a^{-1})(b+b^{-1})}{q+q^{-1}}.
\end{eqnarray*}
\end{lemma}
\noindent {\it Proof:} 
By 
Lemma \ref{lem:ABCvsABCp} the triple 
$\bf B',A',C'$ is the Askey-Wilson triple 
for
$b,a,c$.
Apply 
Lemma \ref{lem:awZ3} to this triple.
\hfill $\Box$ \\

\noindent Throughout this section and the next,
fix an integer $d\geq 1$.
Using $a,b,c,d$ we now define some parameters in $\F$.

\begin{definition}
\label{def:pa}
\rm Define
\begin{eqnarray*}
\theta_n = aq^{2n-d} + a^{-1}q^{d-2n},
\qquad
\qquad 
\theta^*_n = bq^{2n-d} + b^{-1}q^{d-2n}
\end{eqnarray*}
for $0 \leq n \leq d$, 
and
\begin{eqnarray*}
&&\varphi_n = a^{-1}b^{-1}q^{d+1}(q^n-q^{-n})(q^{n-d-1}-q^{d-n+1})
(q^{-n}-abcq^{n-d-1})(q^{-n}-abc^{-1}q^{n-d-1}),
\\
&&\phi_n = ab^{-1}q^{d+1}(q^n-q^{-n})(q^{n-d-1}-q^{d-n+1})
(q^{-n}-a^{-1}bcq^{n-d-1})(q^{-n}-a^{-1}bc^{-1}q^{n-d-1})
\end{eqnarray*}
for $1 \leq n \leq d$.
\end{definition}

\noindent In order to avoid degenerate situations in
Definition
\ref{def:pa}, we sometimes impose
restrictions on how $a,b,c,d$ are related.
We now describe these restrictions.

\begin{definition} 
\label{def:feas}
\rm
The sequence $(a,b,c,d)$ is called {\it feasible} 
whenever
the following {\rm (i), (ii)} hold:
\begin{enumerate}
\item[\rm (i)] neither of $a^2$, $b^2$ is among
$q^{2d-2}, q^{2d-4}, \ldots, q^{2-2d}$;
\item[\rm (ii)] none of $abc$, $a^{-1}bc$, $ab^{-1}c$,
$abc^{-1}$ is among $q^{d-1}, q^{d-3}, \ldots, q^{1-d}$.
\end{enumerate}
\end{definition}

\begin{lemma} The following are equivalent:
\begin{enumerate}
\item[\rm (i)] the sequence $(a,b,c,d)$ is feasible;
\item[\rm (ii)] the
$\lbrace \theta_n \rbrace_{n=0}^d$ are mutually distinct, the
$\lbrace \theta^*_n \rbrace_{n=0}^d$ are mutually distinct, and the
$\lbrace \varphi_n \rbrace_{n=1}^d$,
$\lbrace \phi_n \rbrace_{n=1}^d$
are all nonzero.
\end{enumerate}
\end{lemma}
\noindent {\it Proof:}  
This is routinely checked.
\hfill $\Box$ \\

\noindent The literature on Leonard pairs contains the
notion of a parameter array
\cite{ter:LBUB}, \cite{terPA}. For our purpose we do not need the 
full definition; just the following feature.

\begin{lemma}
{\rm 
\cite[Lemma~7.3]{hwhuang}}.
The following are equivalent:
\begin{enumerate}
\item[\rm (i)] the sequence $(a,b,c,d)$ is feasible;
\item[\rm (iii)] the sequence 
$(\lbrace 
\theta_n \rbrace_{n=0}^d;
\lbrace \theta^*_n \rbrace_{n=0}^d;
\lbrace \varphi_n \rbrace_{n=1}^d;
\lbrace \phi_n \rbrace_{n=1}^d
)$
is a parameter array.
\end{enumerate}
\end{lemma}

\begin{lemma} 
\label{lem:fs}
Assume that $(a,b,c,d)$ is feasible. Then each of
the following sequences is feasible:
\begin{eqnarray}
(a^{-1},b,c,d); \quad \quad
(a,b^{-1},c,d); \quad \quad
(a,b,c^{-1},d); \quad \quad
(b,a,c,d).
\label{eq:moves}
\end{eqnarray}
\end{lemma}
\noindent {\it Proof:}  
By Definition
\ref{def:feas}.
\hfill $\Box$ \\

\noindent  In 
Lemma
\ref{lem:fs} we gave some feasible sequences.
We now consider how their parameter arrays are related.

\begin{lemma} 
\label{lem:patable}
Assume that $(a,b,c,d)$ is feasible.
In the table below, each row contains a feasible sequence
and the corresponding parameter array.

\medskip
\centerline{
\begin{tabular}[t]{c|c}
 {\rm feasible sequence} 
  & 
 {\rm corresp. parameter array} 
\\
\hline
$(a,b,c,d)$ & 
$
(\lbrace \theta_n \rbrace_{n=0}^d;
\lbrace \theta^*_n \rbrace_{n=0}^d;
\lbrace \varphi_n \rbrace_{n=1}^d;
\lbrace \phi_n \rbrace_{n=1}^d
)$
\\
$(a^{-1},b,c,d)$ & 
$(\lbrace\theta_{d-n} \rbrace_{n=0}^d;
\lbrace \theta^*_n \rbrace_{n=0}^d;
\lbrace \phi_n \rbrace_{n=1}^d;
\lbrace \varphi_n \rbrace_{n=1}^d
)$
\\
$(a,b^{-1},c,d)$ & 
$(\lbrace \theta_n \rbrace_{n=0}^d;
\lbrace \theta^*_{d-n} \rbrace_{n=0}^d;
\lbrace \phi_{d-n+1} \rbrace_{n=1}^d;
\lbrace \varphi_{d-n+1} \rbrace_{n=1}^d
)$
\\
$(a,b,c^{-1},d)$ & 
$(\lbrace \theta_n \rbrace_{n=0}^d;
\lbrace \theta^*_n \rbrace_{n=0}^d;
\lbrace \varphi_n \rbrace_{n=1}^d;
\lbrace \phi_n \rbrace_{n=1}^d
)$
\\
$(b,a,c,d)$ & 
$(\lbrace \theta^*_n \rbrace_{n=0}^d;
\lbrace \theta_n \rbrace_{n=0}^d;
\lbrace \varphi_n \rbrace_{n=1}^d;
\lbrace \phi_{d-n+1} \rbrace_{n=1}^d
)$
     \end{tabular}}
\end{lemma}
\noindent {\it Proof:}  
Use Definition
\ref{def:pa}.
\hfill $\Box$ \\

%

\begin{lemma}
{\rm \cite[Theorem~1.9]{terLS}}.
\label{lem:abclp}
Assume 
that $(a,b,c,d)$ is feasible. Then there exists a 
Leonard pair over $\F$
that is described as follows. In one basis the pair
is represented by
    \begin{eqnarray}
      \left(
   \begin{array}{cccccc}
  \theta_0 &    & & & &{\bf 0}  \\
  1 &  \theta_1&  & & &\\
  &   1 &  \theta_2 & &  &\\
  & & \cdot & \cdot & &\\
  & & &\cdot & \cdot & \\
 {\bf 0} & &  && 1 & \theta_d 
 \end{array}
 \right),
\qquad \qquad
      \left(
   \begin{array}{cccccc}
  \theta^*_0 & \varphi_1   & & & & {\bf 0}  \\
   &  \theta^*_1& \varphi_2 & & &\\
  &    &  \theta^*_2 & \cdot & &\\
  & &  & \cdot & \cdot &\\
  & &  & & \cdot & \varphi_d \\
 {\bf 0} & & & & & \theta^*_d 
 \end{array}
 \right).
\label{eq:aas}
\end{eqnarray}
In another basis the 
 pair is represented by
    \begin{eqnarray}
      \left(
   \begin{array}{cccccc}
  \theta_d &    & & & &{\bf 0}  \\
  1 &  \cdot &  & & &\\
   &  \cdot & \cdot & & &\\
  &&   \cdot &  \theta_2 & & \\
  && &    1 & \theta_1 & \\
 {\bf 0} && & & 1 & \theta_0 
 \end{array}
 \right),
\qquad \qquad
      \left(
   \begin{array}{cccccc}
  \theta^*_0 & \phi_1 &  & & & {\bf 0}  \\
   &  \theta^*_1& \phi_2 & & &\\
  &    &  \theta^*_2 & \cdot & &\\
  &  &  &  \cdot & \cdot & \\
  & & & & \cdot & \phi_d \\
 {\bf 0} && & &  & \theta^*_d 
 \end{array}
 \right).
\label{eq:aas2}
\end{eqnarray}
Up to isomorphism
the above Leonard pair is uniquely determined by
$(a,b,c,d)$.
\end{lemma}

\begin{definition}\rm
Assume that $(a,b,c,d)$ is feasible.
The Leonard pair from
Lemma \ref{lem:abclp} is said to {\it correspond} to
 $(a,b,c,d)$.
\end{definition}

\noindent Not every Leonard pair arises from the
construction of 
Lemma \ref{lem:abclp}. The ones that do are said to
have {\it $q$-Racah type}
\cite[Section~5]{hwhuang},
\cite[Example~5.3]{terPA}.
\medskip

\noindent A Leonard pair of $q$-Racah type corresponds
to more than one feasible sequence. This is explained in the
next result.

\begin{lemma}
\label{lem:same}
Let $A, B$ denote a Leonard pair of
$q$-Racah type, with feasible sequence
$(a,b,c,d)$.
Then each of 
\begin{eqnarray}
 \label{eq:obv1}
&&
(a,b,c,d); \quad \quad \;
(a^{-1},b,c,d); \quad \quad\;
(a,b^{-1},c,d); \quad \quad \;
(a,b,c^{-1},d);
\\&&
 \label{eq:obv2}
(a,b^{-1},c^{-1},d); \quad 
(a^{-1},b,c^{-1},d); \quad 
(a^{-1},b^{-1},c,d); \quad 
(a^{-1},b^{-1},c^{-1},d)
\end{eqnarray}
is a feasible sequence for $A, B$.
The Leonard pair $A,B$  has no other feasible sequence.
\end{lemma}
\noindent {\it Proof:}  
This follows from 
\cite[Lemma~12.2]{ter:LBUB}
and
Lemma
\ref{lem:patable}.
\hfill $\Box$ \\

\begin{lemma} 
\label{lem:BA}
Let $A, B$ denote a Leonard pair
of $q$-Racah type, with feasible sequence $(a,b,c,d)$.
Then the Leonard pair $B, A$ is of $q$-Racah type,
with feasible sequence
$(b,a,c,d)$.
\end{lemma}
\noindent {\it Proof:}  
By \cite[Theorem~1.11]{terLS}
and
Lemma
\ref{lem:patable}.
\hfill $\Box$ \\

\begin{proposition} 
\label{prop:bringinC}
{\rm \cite[Theorem~10.1]{hwhuang}.}
Let $A,B$ denote a Leonard pair over $\F$
that has $q$-Racah type.
%
Let $V$ denote the underlying vector space.
Then there exists a unique  $C \in
{\rm End}(V)$ such that
\begin{eqnarray*}
&&
{A} + \frac{q{BC}-q^{-1}{CB}}{q^2-q^{-2}} = \frac{(a+a^{-1})
(q^{d+1}+q^{-d-1})
+(b+b^{-1})(c+c^{-1})}{q+q^{-1}},
\\
&&
{ B} + \frac{q{CA}-q^{-1}{AC}}{q^2-q^{-2}} = \frac{(b+b^{-1})
(q^{d+1}+q^{-d-1})
+(c+c^{-1})(a+a^{-1})}{q+q^{-1}},
\\
&&
{ C} + \frac{q{AB}-q^{-1}{BA}}{q^2-q^{-2}} = \frac{(c+c^{-1}) 
(q^{d+1}+q^{-d-1})
+(a+a^{-1})(b+b^{-1})}{q+q^{-1}}.
\end{eqnarray*}
Here $(a,b,c,d)$ denotes a feasible sequence for the Leonard pair
$A,B$.
\end{proposition}

\begin{definition}\rm Referring to Proposition
\ref{prop:bringinC}, we call $C$ the {\it $\Z_3$-symmetric
completion} of the Leonard pair $A$, $B$.
\end{definition}

\begin{definition}
\label{def:dualcompletion}
\rm 
Let $A,B$ denote a Leonard pair of $q$-Racah type.
By the {\it dual
$\Z_3$-symmetric
completion} of $A,B$ we mean the
$\Z_3$-symmetric
completion of the Leonard pair $B,A$.
\end{definition}

\begin{proposition} 
\label{prop:bringinCp}
Let $A,B$ denote a Leonard pair 
of $q$-Racah type, with
dual $\Z_3$-symmetric completion $C'$.
 Then
\begin{eqnarray*}
&&
{A} + \frac{q{C'B}-q^{-1}{BC'}}{q^2-q^{-2}} = \frac{(a+a^{-1})
(q^{d+1}+q^{-d-1})
+(b+b^{-1})(c+c^{-1})}{q+q^{-1}},
\\
&&
{ B} + \frac{q{AC'}-q^{-1}{C'A}}{q^2-q^{-2}} = \frac{(b+b^{-1})
(q^{d+1}+q^{-d-1})
+(c+c^{-1})(a+a^{-1})}{q+q^{-1}},
\\
&&
{ C'} + \frac{q{BA}-q^{-1}{AB}}{q^2-q^{-2}} = \frac{(c+c^{-1}) 
(q^{d+1}+q^{-d-1})
+(a+a^{-1})(b+b^{-1})}{q+q^{-1}}.
\end{eqnarray*}
Here $(a,b,c,d)$ denotes a feasible sequence for the Leonard pair
$A,B$. 
\end{proposition}
\noindent {\it Proof:} 
By Definition
\ref{def:dualcompletion}, $C'$ is the
$\Z_3$-symmetric completion of the Leonard pair
$B,A$.
By Lemma
\ref{lem:BA}, the Leonard pair $B,A$ has a feasible sequence $(b,a,c,d)$.
Using these comments,
apply Proposition
\ref{prop:bringinC} to the Leonard pair $B,A$.
\hfill $\Box$ \\

\noindent Let $A,B$ denote a Leonard pair that has
$q$-Racah type. 
We now relate  its
$\Z_3$-symmetric completion  and
dual $\Z_3$-symmetric completion.

\begin{lemma}
\label{lem:diff}
Let $A,B$ denote a Leonard pair 
of $q$-Racah type, with
$\Z_3$-symmetric completion $C$ and
dual $\Z_3$-symmetric completion $C'$.
Then
\begin{eqnarray*}
C'-C = \frac{AB-BA}{q-q^{-1}}.
\end{eqnarray*}
\end{lemma}
\noindent {\it Proof:}  
Subtract the last equation in
Proposition 
\ref{prop:bringinC}
from the last equation in
Proposition
\ref{prop:bringinCp}.
\hfill $\Box$ \\

\begin{lemma} 
{\rm \cite[Theorem~5.8]{tersplit}.}
\label{lem:antiaut}
Let $A,B$ denote a Leonard pair over $\F$.
Let $V$ denote the underlying vector space.
Then there exists a unique antiautomorphism
$\dagger$ of 
${\rm End}(V)$ that fixes each of $A$, $B$.
Moreover $\dagger^2=1$.
\end{lemma}

\begin{lemma} 
With reference to Lemma
\ref{lem:antiaut}, assume that
 $A,B$
has $q$-Racah type. Then
its $\Z_3$-symmetric completion $C$ and
 dual $\Z_3$-symmetric completion $C'$
are swapped by
$\dagger$.
\end{lemma}
\noindent {\it Proof:}  
In the last equation of
Proposition
\ref{prop:bringinC}, apply $\dagger$ to
each term. Compare the resulting equation with the
last equation of 
Proposition
\ref{prop:bringinCp}.
\hfill $\Box$ \\

\medskip

\noindent We now use
$U_q(\mathfrak{sl}_2)$ to construct 
 Leonard pairs of $q$-Racah type.

\begin{theorem}
\label{thm:ABCmain}
Assume that $(a,b,c,d)$ is feasible.
Let ${\bf A,B,C}$ denote the Askey-Wilson triple
for $a,b,c$.
Then the following {\rm (i)--(iii)} hold.
\begin{enumerate}
\item[\rm (i)] The pair ${\bf A,  B}$ acts on the
$U_q(\mathfrak{sl}_2)$-module $\mathbf V_d$ as a Leonard pair.
\item[\rm(ii)] This
Leonard pair 
corresponds to $(a,b,c,d)$.
\item[\rm(iii)] 
The element $\bf C$ 
acts on $\mathbf V_d$ as the 
$\Z_3$-symmetric completion of this Leonard pair.
\end{enumerate}
\end{theorem}
\noindent {\it Proof:} (i), (ii) 
The parameter array
$(
\lbrace \theta_n \rbrace_{n=0}^d;
\lbrace \theta^*_n \rbrace_{n=0}^d;
\lbrace \varphi_n \rbrace_{n=1}^d;
\lbrace \phi_n \rbrace_{n=1}^d)$
for $(a,b,c,d)$ is shown in Definition
\ref{def:pa}.
Let $\lbrace u_n \rbrace_{n=0}^d$ denote a
$\lbrack y \rbrack_{row}^{inv}$-basis for $\mathbf V_d$.
Consider the matrices in
${\rm Mat}_{d+1}(\F)$ 
that represent $x,y,z$ with respect to
$\lbrace u_n \rbrace_{n=0}^d$.
These matrices are given in Lemma 3.22.
The matrix representing $x$ is lower bidiagonal,
with $(n,n)$-entry $q^{d-2n}$ for
$0 \leq n \leq d$ and
$(n,n-1)$-entry $q^d-q^{d-2n}$ for
$1 \leq n \leq d$.
The matrix representing $y$ is diagonal,
with $(n,n)$-entry $q^{2n-d}$ for
$0 \leq n \leq d$.
The matrix representing $z$ is upper bidiagonal,
with $(n,n)$-entry $q^{d-2n}$ for
$0 \leq n \leq d$ and
$(n-1,n)$-entry $q^{-d}-q^{d-2n+2}$ for
$1 \leq n \leq d$.
Using this data and
Definition
\ref{def:ABC}
we compute the matrices in 
${\rm Mat}_{d+1}(\F)$ 
that represent $\bf A,B$ with respect to
$\lbrace u_n \rbrace_{n=0}^d$.
The matrix representing $\bf A$ is lower bidiagonal,
with $(n,n)$-entry $\theta_{d-n}$ for
$0 \leq n \leq d$ and
$(n,n-1)$-entry 
\begin{eqnarray}
\label{eq:Aentry}
aq^d(q^n-q^{-n})(q^{-n}-a^{-1}bc^{-1}q^{n-d-1})
\end{eqnarray}
for $1 \leq n \leq d$.
The matrix representing $\bf B$ is upper bidiagonal,
with $(n,n)$-entry $\theta^*_{n}$ for
$0 \leq n \leq d$ and
$(n-1,n)$-entry 
\begin{eqnarray}
\label{eq:Bentry}
b^{-1}q(q^{n-d-1}-q^{d-n+1})(q^{-n}-a^{-1}bcq^{n-d-1})
\end{eqnarray}
for $1 \leq n \leq d$.
We now adjust the basis
$\lbrace u_n \rbrace_{n=0}^d$.
For $1 \leq n \leq d$ let
$\alpha_n$ (resp. $\beta_i$) denote the scalar
(\ref{eq:Aentry}) 
(resp. (\ref{eq:Bentry})).
We have $\alpha_n \beta_n = \phi_n$, so each of
 $\alpha_n,  \beta_n$ is nonzero.
Define $u'_n = \alpha_1 \alpha_2 \cdots \alpha_n u_n$
for $0 \leq n \leq d$.
By construction $\lbrace u'_n \rbrace_{n=0}^d$
is a basis for $\mathbf V_d$. With respect to this
basis the matrices representing $\bf A, B$
are the ones shown in
(\ref{eq:aas2}).
Therefore the pair $\bf A,B$ acts on $\mathbf V_d$
as a Leonard pair that corresponds to
$(a,b,c,d)$. 
\\
\noindent (iii)
Recall that $\Lambda$ acts on $\mathbf V_d$ as
$q^{d+1}+q^{-d-1}$ times the identity.
Using this comment, compare the
last equation in Lemma
\ref{lem:awZ3}
with the last equation in
Proposition 
\ref{prop:bringinC}.
\hfill $\Box$ \\

\begin{lemma}
\label{lcor:abcuq}
Let $A,B$ denote a Leonard pair over $\F$
of $q$-Racah type, 
 with feasible sequence $(a,b,c,d)$.
Let ${\bf A,B,C}$ denote the Askey-Wilson triple
for $a,b,c$.
Assume that there exists a
 $U_q(\mathfrak{sl}_2)$-module structure on the underlying
vector space $V$ such that $A=\bf A$ and $B=\bf B$ on
$V$. Then the
 $U_q(\mathfrak{sl}_2)$-module $V$ is irreducible.
\end{lemma}
\noindent {\it Proof:} By
Lemma \ref{lem:LPQR}.
\hfill $\Box$ \\


\begin{corollary}
\label{cor:abcuq}
Let $A,B$ denote a Leonard pair over $\F$
of $q$-Racah type, 
 with feasible sequence $(a,b,c,d)$.
Let ${\bf A,B,C}$ denote the Askey-Wilson triple
for $a,b,c$.
Then  the following {\rm (i), (ii)}  hold.
\begin{enumerate}
\item[\rm (i)]
There exists a
unique type 1 $U_q(\mathfrak{sl}_2)$-module structure on the underlying
vector space $V$ such that $A=\bf A$ and $B=\bf B$ on
$V$.
   \item[\rm (ii)] The element $\bf C$ acts on $V$
    as the $\Z_3$-symmetric completion of $A,B$.
\end{enumerate}
\end{corollary}
\noindent {\it Proof:}  
The Leonard pair $A,B$ corresponds to
$(a,b,c,d)$.
The Leonard pair in Theorem
\ref{thm:ABCmain} also corresponds to
$(a,b,c,d)$. Therefore these Leonard pairs
are isomorphic.
Let $\zeta:V\to \mathbf V_d$ denote an isomorphism
of Leonard pairs 
from $A,B$ to the Leonard pair
in 
 Theorem
\ref{thm:ABCmain}.
Via $\zeta$ we transport the 
$U_q(\mathfrak{sl}_2)$-module structure from $\mathbf V_d$
to 
 $V$. 
This turns $V$ into a 
$U_q(\mathfrak{sl}_2)$-module that is isomorphic to $\mathbf V_d$.
By construction $A=\bf A$ and $B=\bf B$ on $V$.
Also by Theorem
\ref{thm:ABCmain}(iii),
$\bf C$ acts on $V$ as the $\Z_3$-symmetric completion
of $A,B$.
The uniqueness assertion in (i)
follows from 
Lemma \ref{lem:LPQR}.
\hfill $\Box$ \\


\begin{theorem}
\label{thm:ABCpmain}
Assume that $(a,b,c,d)$ is feasible.
Let ${\bf A',B',C'}$ denote the dual Askey-Wilson triple
for $a,b,c$.
Then the following {\rm (i)--(iii)} hold.
\begin{enumerate}
\item[\rm (i)] The pair ${\bf A', B'}$ acts on the
$U_q(\mathfrak{sl}_2)$-module $\mathbf V_d$ as a Leonard pair.
\item[\rm (ii)]
This Leonard pair 
corresponds to $(a,b,c,d)$.
\item[\rm (iii)]
The element $\bf C'$ acts on $\mathbf V_d$ as the 
dual $\Z_3$-symmetric completion of this Leonard pair.
\end{enumerate}
\end{theorem}
\noindent {\it Proof:} Conceptually this proof is
similar to the proof of Theorem
\ref{thm:ABCmain}, but as the details are different
they will be displayed.
\\
\noindent
(i), (ii) 
The parameter array
$(
\lbrace \theta_n \rbrace_{n=0}^d;
\lbrace \theta^*_n \rbrace_{n=0}^d;
\lbrace \varphi_n\rbrace_{n=1}^d;
\lbrace \phi_n \rbrace_{n=1}^d)$
for $(a,b,c,d)$ is shown in Definition
\ref{def:pa}.
Let $\lbrace v_n \rbrace_{n=0}^d$ denote a
$\lbrack y \rbrack_{row}$-basis for $\mathbf V_d$.
Consider the matrices in
${\rm Mat}_{d+1}(\F)$ 
that represent $x,y,z$ with respect to
$\lbrace v_n \rbrace_{n=0}^d$.
These matrices are given in Lemma 3.9.
The matrix representing $x$ is upper bidiagonal,
with $(n,n)$-entry $q^{2n-d}$ for
$0 \leq n \leq d$ and
$(n-1,n)$-entry $q^d-q^{2n-d-2}$ for
$1 \leq n \leq d$.
The matrix representing $y$ is diagonal,
with $(n,n)$-entry $q^{d-2n}$ for
$0 \leq n \leq d$.
The matrix representing $z$ is lower bidiagonal,
with $(n,n)$-entry $q^{2n-d}$ for
$0 \leq n \leq d$ and
$(n,n-1)$-entry $q^{-d}-q^{2n-d}$ for
$1 \leq n \leq d$.
Using this data and
Definition
\ref{def:ABCp}
we compute the matrices in 
${\rm Mat}_{d+1}(\F)$ 
that represent $\bf A',B'$ with respect to
$\lbrace v_n \rbrace_{n=0}^d$.
The matrix representing $\bf A'$ is lower bidiagonal,
with $(n,n)$-entry $\theta_{d-n}$ for
$0 \leq n \leq d$ and
$(n,n-1)$-entry 
\begin{eqnarray}
\label{eq:Apentry}
b^{-1}cq(q^n-q^{-n})(q^{-n}-a^{-1}bc^{-1}q^{n-d-1})
\end{eqnarray}
for $1 \leq n \leq d$.
The matrix representing $\bf B'$ is upper bidiagonal,
with $(n,n)$-entry $\theta^*_{n}$ for
$0 \leq n \leq d$ and
$(n-1,n)$-entry 
\begin{eqnarray}
\label{eq:Bpentry}
ac^{-1}q^d(q^{n-d-1}-q^{d-n+1})(q^{-n}-a^{-1}bcq^{n-d-1})
\end{eqnarray}
for $1 \leq n \leq d$.
We now adjust the basis
$\lbrace v_n \rbrace_{n=0}^d$.
For $1 \leq n \leq d$ let
$\alpha_n$ (resp. $\beta_n$) denote the scalar
(\ref{eq:Apentry}) 
(resp. (\ref{eq:Bpentry})).
We have $\alpha_n \beta_n = \phi_n$, so each of
 $\alpha_n,  \beta_n$ is nonzero.
Define $v'_n = \alpha_1 \alpha_2 \cdots \alpha_n v_n$
for $0 \leq n \leq d$.
By construction $\lbrace v'_n \rbrace_{n=0}^d$
is a basis for $\mathbf V_d$. With respect to this
basis the matrices representing $\bf A', B'$
are the ones shown in
(\ref{eq:aas2}).
Therefore the pair $\bf A', B'$ acts on $\mathbf V_d$
as a Leonard pair that corresponds to
$(a,b,c,d)$. 
\\
\noindent (iii)
Recall that $\Lambda$ acts on $\mathbf V_d$ as
$q^{d+1}+q^{-d-1}$ times the identity.
Using this comment, compare the
last equation in Lemma
\ref{lem:awp}
with the last equation in
Proposition 
\ref{prop:bringinCp}.
\hfill $\Box$ \\

\begin{lemma}
\label{lcor:dualabcuq}
Let $A,B$ denote a Leonard pair over $\F$
of $q$-Racah type, 
 with feasible sequence $(a,b,c,d)$.
Let ${\bf A',B',C'}$ denote the dual Askey-Wilson triple
for $a,b,c$.
Assume that there exists a
 $U_q(\mathfrak{sl}_2)$-module structure on the underlying
vector space $V$ such that $A=\bf A'$ and $B=\bf B'$ on
$V$. Then the
 $U_q(\mathfrak{sl}_2)$-module $V$ is irreducible.
\end{lemma}
\noindent {\it Proof:} By
Lemma \ref{lem:LPQR}.
\hfill $\Box$ \\


\begin{corollary}
\label{cor:abcpuq}
Let $A,B$ denote a Leonard pair over $\F$
of $q$-Racah type, 
 with feasible sequence $(a,b,c,d)$.
Let $\bf A',
 B',
 C'$ denote the dual
Askey-Wilson triple for
$a,b,c$.
Then the following {\rm (i), (ii)} hold.
\begin{enumerate}
\item[\rm (i)]
There exists a
unique type 1 $U_q(\mathfrak{sl}_2)$-module structure on the underlying
vector space $V$ such that $A=\bf A'$ and $B=\bf B'$ on
$V$.
 \item[\rm (ii)] The element $\bf C'$ acts on $V$
 as the dual $\Z_3$-symmetric completion of $A,B$.
\end{enumerate}
\end{corollary}
\noindent {\it Proof:}  
Similar to the proof of 
Corollary
\ref{cor:abcuq}.
\hfill $\Box$ \\

\noindent We now bring in $\boxtimes_q$.
Recall the injections
$\kappa_i: 
U_q(\mathfrak{sl}_2) \to \boxtimes_q$ from
Lemma 4.4.

\begin{definition}
\label{def:newABC}
\rm Referring to Definition
\ref{def:ABC}, we identify 
$\bf A,
B, C$ with their images 
under the injection 
$\kappa_2:
U_q(\mathfrak{sl}_2) \to \boxtimes_q$. 
 Thus
\begin{eqnarray*}
{\bf A} &=& ax_{01}+a^{-1}x_{13} +
bc^{-1}\frac{\lbrack x_{01}, x_{13}
\rbrack}{q-q^{-1}},
\\
{\bf B} &=& bx_{13}+b^{-1}x_{30} +
ca^{-1}\frac{\lbrack x_{13}, x_{30}
\rbrack}{q-q^{-1}},
\\
{\bf C} &=& cx_{30}+c^{-1}x_{01} +
ab^{-1}\frac{\lbrack x_{30}, x_{01}
\rbrack}{q-q^{-1}}.
\end{eqnarray*}
\end{definition}

\begin{definition}
\label{def:newABCp}
\rm 
Referring to Definition
\ref{def:ABCp}, we identify 
$\bf A',
B', C'$ with their images under the
injection
$\kappa_0:
U_q(\mathfrak{sl}_2) \to \boxtimes_q$. 
 Thus
\begin{eqnarray*}
{\bf A'} &=& ax_{31}+a^{-1}x_{12} +
cb^{-1}\frac{\lbrack x_{31}, x_{12}
\rbrack}{q-q^{-1}},
\\
{\bf B'} &=& bx_{23}+b^{-1}x_{31} +
ac^{-1}\frac{\lbrack x_{23}, x_{31}
\rbrack}{q-q^{-1}},
\\
{\bf C'} &=& cx_{12}+c^{-1}x_{23} +
ba^{-1}\frac{\lbrack x_{12}, x_{23}
\rbrack}{q-q^{-1}}.
\end{eqnarray*}
\end{definition}

\begin{lemma}
\label{lem:ABexposed}
Assume that $(a,b,c,d)$ is feasible,
and define
$t=abc^{-1}$. On the $\boxtimes_q$-module $\mathbf V_d(t)$,
\begin{eqnarray}
&&{\bf A} =  
ax_{01} + a^{-1}x_{12} 
=
{\bf A'},
\label{eq:unpack1}
\\
&&{\bf B} =  
bx_{23} + b^{-1}x_{30}
=
{\bf B'}.
\label{eq:unpack2}
\end{eqnarray}
\end{lemma}
\noindent {\it Proof:}  
Use Lemma 9.5 (with $t=abc^{-1}$).
\hfill $\Box$ \\



\begin{theorem} 
\label{thm:tetqLP}
Assume that $(a,b,c,d)$ is feasible, and define
$t=abc^{-1}$. 
Then the following {\rm (i)--(iv)} hold.
\begin{enumerate}
\item[\rm (i)]
The pair 
\begin{eqnarray}
 a x_{01} + a^{-1}x_{12},
\qquad \qquad 
 b x_{23} + b^{-1}x_{30}
\label{eq:ABdisguise}
\end{eqnarray}
acts on the $\boxtimes_q$-module $\mathbf V_d(t)$ as a Leonard pair
of $q$-Racah type.
\item[\rm (ii)]
This Leonard pair corresponds to
$(a,b,c,d)$.
\item[\rm (iii)]
The element $\bf C$ from
Definition \ref{def:newABC}
acts on $\mathbf V_d(t)$ as the
$\Z_3$-symmetric completion of this Leonard pair.
\item[\rm (iv)]
The element $\bf C'$ 
from Definition \ref{def:newABCp}
acts on $\mathbf V_d(t)$ as the
dual $\Z_3$-symmetric completion of this
Leonard pair.
\end{enumerate}
\end{theorem}
\noindent {\it Proof:}  
We first obtain
(i)--(iii).
Using the homomorphism
$\kappa_2 :
U_q(\mathfrak{sl}_2) \to \boxtimes_q$ we turn
the $\boxtimes_q$-module $\mathbf V_d(t)$ into
a 
$U_q(\mathfrak{sl}_2)$-module isomorphic to
$\mathbf V_d$. Apply
Theorem
\ref{thm:ABCmain} to this
$U_q(\mathfrak{sl}_2)$-module,
and use the equations on the left in
(\ref{eq:unpack1}),
(\ref{eq:unpack2}).
This yields (i)--(iii).
Next we obtain (iv).
Using the homomorphism
$\kappa_0 :
U_q(\mathfrak{sl}_2) \to \boxtimes_q$ we turn
the $\boxtimes_q$-module $\mathbf V_d(t)$ into
a 
$U_q(\mathfrak{sl}_2)$-module isomorphic to
$\mathbf V_d$. Apply
Theorem
\ref{thm:ABCpmain} to this
$U_q(\mathfrak{sl}_2)$-module,
 and 
use the equations on the right in
(\ref{eq:unpack1}),
(\ref{eq:unpack2}). This yields (iv), along with a second proof of
(i), (ii).
\hfill $\Box$ \\

\begin{corollary}
\label{thm:tetqAB}
Let $A,B$ denote a Leonard pair over $\F$ of $q$-Racah
type, with feasible sequence $(a,b,c,d)$.
Define $t=abc^{-1}$. Then 
the following {\rm (i)--(iii)} hold.
\begin{enumerate}
\item[\rm (i)]
The underlying
vector space $V$ supports
a unique $t$-evaluation module for $\boxtimes_q$ such that
on $V$,
\begin{eqnarray}
A = a x_{01} + a^{-1}x_{12},
\qquad \qquad 
B = b x_{23} + b^{-1}x_{30}.
\label{eq:ABneed}
\end{eqnarray}
\item[\rm (ii)] The element ${\bf C}$ 
from Definition \ref{def:newABC}
acts on
$V$ as the $\Z_3$-symmetric completion of 
$A,B$.
\item[\rm (iii)] The element ${\bf C'}$
from Definition \ref{def:newABCp}
acts on
$V$ as the dual $\Z_3$-symmetric completion of 
$A,B$.
\end{enumerate}
\end{corollary}
\noindent {\it Proof:}  
The Leonard pair $A,B$ corresponds to $(a,b,c,d)$.
The Leonard pair in
Theorem
\ref{thm:tetqLP} also corresponds to $(a,b,c,d)$.
Therefore these Leonard pairs are isomorphic.
Let $\partial : V \to \mathbf V_d(t)$ denote
an isomorphism of Leonard pairs from $A,B$ to the
Leonard pair in Theorem
\ref{thm:tetqLP}.
Via $\partial$ we transport the
$\boxtimes_q$-module structure from $\mathbf V_d(t)$ to $V$.
This turns $V$ into a $\boxtimes_q$-module isomorphic to
$\mathbf V_d(t)$. By construction and Theorem
\ref{thm:tetqLP}, the equations 
(\ref{eq:ABneed}) hold on $V$, and the assertions
(ii), (iii) are valid.
The uniqueness assertion in (i) follows from
Lemma
\ref{lem:LPQR}.
\hfill $\Box$ \\

\begin{note}\rm 
Let $A,B$ denote a Leonard pair of $q$-Racah type.
Using
Lemma
\ref{lem:same}
and Corollary 
\ref{thm:tetqAB}, we get eight $\boxtimes_q$-module
structures on the underlying vector space.
\end{note}

\section{The compact basis}

\noindent
Let $A,B$ denote a Leonard pair over $\F$ that has  
$q$-Racah type.
In this section we introduce a certain 
basis for the underlying vector space,
with respect to 
which the matrices representing $A$ and $B$ are
 tridiagonal with attractive entries. We call this
 basis the compact basis.
 We show
  how the compact basis
 is related to the $\boxtimes_q$-module structure
 discussed in Corollary
\ref{thm:tetqAB}. 

\begin{proposition}
\label{lem:2TD}
Let $A,B$ denote a Leonard pair over $\F$ of $q$-Racah type,
with feasible sequence $(a,b,c,d)$.
Then there exists a basis for the underlying vector space $V$,
with respect to which $A$ and $B$ are represented
by the following tridiagonal matrices in
${\rm Mat}_{d+1}(\F)$:

\medskip
\centerline{
\begin{tabular}[t]{c|ccc}
{\rm element} & 
   {\rm $(n,n-1)$-entry}  &  
   {\rm $(n,n)$-entry}  &  
   {\rm $(n-1,n)$-entry} 
   \\ \hline 
    $A$  
&
$c^{-1}(1-q^{-2n})$ 
&
$(a+a^{-1})q^{d-2n}$
&
$c(1-q^{2d-2n+2})$ 
\\
$B$
&
$q^{-d-1}(1-q^{2n})$ 
&
$(b+b^{-1})q^{2n-d}$ 
&
$q^{d+1}(1-q^{2n-2d-2})$ 
     \end{tabular}}
     \medskip
\noindent
\end{proposition}
\noindent {\it Proof:}  
Recall the $\boxtimes_q$-module
structure on $V$ from
Corollary
\ref{thm:tetqAB}. 
This is an evaluation module, with
evaluation parameter $t=abc^{-1}$.
Let $\lbrace u_n\rbrace_{n=0}^d $ denote a $\lbrack 1,3,0,2\rbrack$-basis
for $V$. 
Consider the
matrices in
${\rm Mat}_{d+1}(\F)$  that represent
$x_{01}$, 
$x_{12}$, 
$x_{23}$, 
$x_{30}$ 
with respect to 
 $\lbrace u_n\rbrace_{n=0}^d$.
These matrices are given in the last row of 
the first table in Theorem 11.1 (with $r=3$).
Their entries are described as follows.
The matrix representing $x_{01}$ is upper bidiagonal, 
with $(n,n)$-entry $q^{d-2n}$ for $0 \leq n \leq d$
and $(n-1,n)$-entry $q^{-1}t^{-1}(1-q^{2d-2n+2})$
for $1 \leq n \leq d$.
The matrix representing $x_{12}$ is lower bidiagonal, 
with $(n,n)$-entry $q^{d-2n}$ for $0 \leq  n \leq d$
and $(n,n-1)$-entry $qt(1-q^{-2n})$
for $1 \leq n \leq d$.
The matrix representing $x_{23}$ is lower bidiagonal, 
with $(n,n)$-entry $q^{2n-d}$ for $0 \leq n \leq d$
and $(n,n-1)$-entry $q^{-d}(1-q^{2n})$ 
for $1 \leq n \leq d$.
The matrix representing $x_{30}$ is upper bidiagonal, 
with $(n,n)$-entry $q^{2n-d}$ for $0 \leq n \leq d$
and $(n-1,n)$-entry $q^d(1-q^{2n-2d-2})$
for $1 \leq n \leq d$.
Now consider the matrices in
${\rm Mat}_{d+1}(\F)$  that represent
 $A$ and $B$
with respect to 
$\lbrace u_n\rbrace_{n=0}^d $.
Their entries are found using (\ref{eq:ABneed}) and
the above comments.
The matrix representing $A$ is tridiagonal,
with $(n,n-1)$-entry $qbc^{-1}(1-q^{-2n})$ for
$1 \leq n \leq d$,
$(n,n)$-entry $(a+a^{-1})q^{d-2n}$ for $0 \leq n \leq d$,
and 
 $(n-1,n)$-entry $q^{-1}b^{-1}c(1-q^{2d-2n+2})$ for
$1 \leq n \leq d$.
The matrix representing $B$ is  tridiagonal,
with $(n,n-1)$-entry $bq^{-d}(1-q^{2n})$ for
$1 \leq n \leq d$,
$(n,n)$-entry $(b+b^{-1})q^{2n-d}$ for $0 \leq n \leq d$,
and 
 $(n-1,n)$-entry $b^{-1}q^{d}(1-q^{2n-2d-2})$ for
$1 \leq n \leq d$.
We now adjust the basis
$\lbrace u_n \rbrace_{n=0}^d$. Define $v_n = q^n b^n u_n$ for
$0 \leq n \leq d$. Then 
$\lbrace v_n \rbrace_{n=0}^d$ is a basis for $V$.
With respect to this basis the matrices representing $A$ and
$B$ are as shown in the theorem statement. 
\hfill $\Box$ \\

\begin{example} \rm The matrices from
Proposition
\ref{lem:2TD} look as follows for
$d=3$. The matrix representing $A$ is
\begin{eqnarray*}
\left(
\begin{array}{ cccc}
(a+a^{-1})q^3    &  c(1-q^6)  & 0 &  0
\\
c^{-1}(1-q^{-2})    & (a+a^{-1})q   & c(1-q^4)  &  0
\\
0    & c^{-1}(1-q^{-4})   & (a+a^{-1})q^{-1}  & c(1-q^2) 
\\
0    &  0   & c^{-1}(1-q^{-6}) & (a+a^{-1})q^{-3} 
\end{array}
\right)
\end{eqnarray*}
The matrix representing $B$ is
\begin{eqnarray*}
\left(
\begin{array}{ cccc}
(b+b^{-1})q^{-3}    &  q^4(1-q^{-6})  & 0 &  0
\\
q^{-4}(1-q^{2})    & (b+b^{-1})q^{-1}   & q^4(1-q^{-4})  &  0
\\
0    & q^{-4}(1-q^{4})   & (b+b^{-1})q  & q^4(1-q^{-2}) 
\\
0    &  0   & q^{-4}(1-q^{6}) & (b+b^{-1})q^{3} 
\end{array}
\right)
\end{eqnarray*}
\end{example}

\begin{definition}
\label{def:compact}
\rm The basis
 discussed in Proposition
\ref{lem:2TD} is said to be {\it compact}.
\end{definition}

\begin{note}\rm Our motivation for Definition
\ref{def:compact} is that the entries shown
in
the table of Proposition
\ref{lem:2TD} are rather concise.
\end{note}

\begin{note}\rm  The existence of the compact basis was
hinted at in
\cite{GYZlinear},
\cite[Section~6]{rosen1},
\cite[Section~2.3]{rosen}.
\end{note}

\noindent We comment on the uniqueness of the compact basis.

\begin{lemma}
\label{lem:2TDunique}
Let $A,B$ denote a Leonard pair over $\F$ of $q$-Racah type,
with feasible sequence $(a,b,c,d)$. Let 
$\lbrace v_n \rbrace_{n=0}^d$ denote a basis for the underlying
vector space $V$ that meets the
requirements of
Proposition
\ref{lem:2TD}. Let
$\lbrace v'_n \rbrace_{n=0}^d$ denote any basis for
 $V$.
 Then
the following are equivalent:
\begin{enumerate}
\item[\rm (i)] 
the basis 
$\lbrace v'_n \rbrace_{n=0}^d$ meets 
the requirements of
Proposition
\ref{lem:2TD};
\item[\rm (ii)] 
there exists a nonzero $\alpha \in \F$ such that
$v'_n = \alpha v_n$ for $0 \leq n \leq d$.
\end{enumerate}
\end{lemma}
\noindent {\it Proof:}  
${\rm (i)}\Rightarrow {\rm (ii)}$  
Consider the map
$\psi \in {\rm End}(V)$ that sends
$v_n \mapsto v'_n$ for $0 \leq n \leq d$.
The matrix representing $A$ with respect to
$\lbrace v_n\rbrace_{n=0}^d$  is equal to
the matrix representing $A$ with respect to 
$\lbrace v'_n\rbrace_{n=0}^d$.
Therefore $\psi$ commutes with $A$.
By a similar argument
 $\psi$ commutes with $B$.
Now by Lemma \ref{lem:LPQR},
$\psi$ commutes with everything in 
${\rm End}(V)$. Consequently there exists 
$\alpha \in \F$ such that $\psi = \alpha I$.
We have
$\alpha \not=0$ since
$\psi \not=0$.
By construction
 $v'_n = \alpha v_n $ for $0 \leq n \leq d$.
\\
\noindent ${\rm (ii)}\Rightarrow {\rm (i)}$ 
Clear.
\hfill $\Box$ \\

\begin{theorem}
Let $A,B$ denote a Leonard pair over $\F$ of $q$-Racah type,
with feasible sequence $(a,b,c,d)$. Let 
$\lbrace v_n \rbrace_{n=0}^d$ denote a basis for the underlying
vector space $V$. Then
the following are equivalent:
\begin{enumerate}
\item[\rm (i)] 
the basis 
$\lbrace v_n \rbrace_{n=0}^d$ meets 
the requirements of
Proposition
\ref{lem:2TD};
\item[\rm (ii)] 
the basis $\lbrace q^{-n}b^{-n}v_n \rbrace_{n=0}^d$ is
a $\lbrack 1,3,0,2\rbrack$-basis of $V$, for
 the $\boxtimes_q$-module structure in
Corollary
\ref{thm:tetqAB}.
\end{enumerate}
\end{theorem}
\noindent {\it Proof:}  
This follows from 
the proof of
Proposition
\ref{lem:2TD}, together with
Lemma
\ref{lem:2TDunique}.
\hfill $\Box$ \\

\begin{proposition} 
\label{prop:CCp}
Let $A,B$ denote a Leonard pair over $\F$ of $q$-Racah type,
with $\Z_3$-symmetric completion $C$ and
 dual $\Z_3$-symmetric completion $C'$.
Let $(a,b,c,d)$ denote a feasible sequence for $A,B$.
Consider the matrices in 
${\rm Mat}_{d+1}(\F)$ 
that
represent 
$C$ and $C'$ with respect to 
a basis from
Proposition
\ref{lem:2TD}.
The matrix representing $C$ is upper triangular, with 
$(n,n)$-entry
\begin{eqnarray*}
cq^{2n-d}+c^{-1}q^{d-2n}
\end{eqnarray*}
for $0 \leq n \leq d$,
$(n-1,n)$-entry
\begin{eqnarray*}
(q^{d-n+1}-q^{n-d-1})\bigl((b+b^{-1})cq^{n}-(a+a^{-1})q^{d-n+1}\bigr)
\end{eqnarray*}
for $1 \leq n \leq d$, 
$(n-2,n)$-entry
\begin{eqnarray*}
cq^{d+1}(q^{d-n+1}-q^{n-d-1})(q^{d-n+2}-q^{n-d-2})
\end{eqnarray*}
for $2 \leq n \leq d$, and all other entries 0.
The matrix representing $C'$ is lower triangular, with
$(n,n)$-entry
\begin{eqnarray*}
cq^{d-2n}+c^{-1}q^{2n-d}
\end{eqnarray*}
for $0 \leq n \leq d$,
$(n,n-1)$-entry
\begin{eqnarray*}
(q^n-q^{-n})\bigl((a+a^{-1})q^{-n}-(b+b^{-1})c^{-1}q^{n-d-1}\bigr)
\end{eqnarray*}
for $1 \leq n \leq d$, 
$(n,n-2)$-entry
\begin{eqnarray*}
c^{-1}q^{-d-1}(q^{n}-q^{-n})(q^{n-1}-q^{1-n})
\end{eqnarray*}
for $2 \leq n \leq d$, and all other entries 0.
\end{proposition}
\noindent {\it Proof:}  
The matrix representing $C$ is obtained
using Proposition
\ref{lem:2TD} and the last equation in
Proposition
\ref{prop:bringinC}. 
The matrix representing $C'$ is similarly
obtained
using the last equation in
Proposition
\ref{prop:bringinCp}. 
\hfill $\Box$ \\

\noindent We now summarize Proposition 
\ref{lem:2TD} and Proposition
\ref{prop:CCp}. 

\begin{corollary} 
Let $A,B$ denote a Leonard pair over $\F$ of $q$-Racah type.
Consider a compact basis for this pair.
In the table below, 
for each of the displayed maps
we describe the matrix that represents it
with respect to the basis.

\medskip
\centerline{
\begin{tabular}[t]{c|c}
{\rm map} & {\rm representing matrix}  
   \\ \hline 
    $A$ & {\rm irred. tridiagonal}
    \\
    $B$ & {\rm irred. tridiagonal}
\\
$ qAB-q^{-1}BA$ & {\rm upper triangular}
\\
$ qBA-q^{-1}AB$ & {\rm lower triangular}
     \end{tabular}}
     \medskip
\end{corollary}



\section{Appendix I: Some matrix definitions}

\noindent In this appendix we define and discuss the matrices
that were used earlier in the paper. 

\medskip

\noindent Fix an integer $d\geq 1$ and a nonzero
$t \in \F$ that is not among
$\lbrace q^{d-2n+1}\rbrace_{n=1}^d$.
Let $I$ denote the identity matrix in
${\rm Mat}_{d+1}(\F)$.
\medskip

\begin{definition}
\label{def:z}
\rm
Let $Z$ 
denote the matrix in $ 
{\rm Mat}_{d+1}(\F)$ with 
$(i,j)$-entry 
$\delta_{i+j,d}$ for
 $0 \leq i,j\leq d$.
Note that $Z^2=I$.
\end{definition}

\begin{example} \rm
For  $d=3$,
\begin{eqnarray*}
&&
Z = \quad 
\left(
\begin{array}{ cccc}
0&0&0&1  \\
  0 &0&1&0\\
   0& 1 &0&0\\
 1&0  &0   &0
\end{array}
\right).
\end{eqnarray*}
\end{example}

\noindent For each lemma in this appendix, the
proof is routine and left to the reader.

\begin{lemma}
\label{Alem:ZBZ}
For $B \in 
{\rm Mat}_{d+1}(\F)$  and 
$0 \leq i,j\leq d$ the following coincide:
\begin{enumerate}
\item[\rm (i)] the $(i,j)$-entry of $ZBZ$;
\item[\rm (i)] the $(d-i,d-j)$-entry of $B$.
\end{enumerate}
\end{lemma}

\noindent We now consider the matrices used in Section 13.
They are
\begin{eqnarray}
\label{eq:three}
D_q(t), \quad \mathcal D_q(t), \quad T_q.
\end{eqnarray}

\medskip
\noindent Recall the notation
\begin{eqnarray*}
(a;q)_n = (1-a)(1-aq) \cdots (1-aq^{n-1})
\qquad \qquad n=0,1,2,\ldots
\end{eqnarray*}
We interpret $(a;q)_0=1$.

\begin{definition}
\label{Adef:dqt}
\rm
Let $D_q(t)$ denote the diagonal matrix in 
${\rm Mat}_{d+1}(\F)$
with $(i,i)$-entry $(tq^{d-2i+1};q^2)_i$
for $0 \leq i \leq d$.
\end{definition}

\begin{example}\rm
For $d=3$,
\begin{eqnarray*}
D_q(t) = {\rm diag} \bigl(1, 1-tq^2, (1-tq^2)(1-t), 
(1-tq^2)(1-t)(1-tq^{-2}) \bigr).
\end{eqnarray*}
\end{example}

\begin{lemma} We have
\begin{eqnarray*}
\bigl(D_q(t)\bigr)^{-1}  = \frac{ZD_{q^{-1}}(t) Z}{(tq^{1-d};q^2)_d}.
\end{eqnarray*}
\end{lemma}

\begin{definition}
\label{Adef:caldqt}
\rm
Let $\mathcal D_q(t)$ denote the diagonal matrix in 
${\rm Mat}_{d+1}(\F)$
with $(i,i)$-entry $t^iq^{i(d-1)}$
for $0 \leq i \leq d$.
\end{definition}

\begin{example}\rm
For $d=3$,
\begin{eqnarray*}
{\mathcal D}_q(t) = {\rm diag} (1, tq^2, t^2q^4, 
t^3q^6).
\end{eqnarray*}
\end{example}

\begin{lemma} We have
\begin{eqnarray*}
\bigl(\mathcal D_q(t)\bigr)^{-1}  = \mathcal D_{q^{-1}}(t^{-1})
=
t^{-d}q^{d(1-d)}Z\mathcal D_q(t)Z.
\end{eqnarray*}
\end{lemma}

\medskip
\noindent We recall some notation. 
For integers $n\geq i \geq 0$ define
\begin{eqnarray*}
\left[\begin{array}{c} n  \\
i \end{array} \right]_q &=& 
\frac{\lbrack n \rbrack^!_q}{\lbrack i \rbrack^!_q \lbrack n-i\rbrack^!_q }.
\end{eqnarray*}

\begin{definition}\rm Let $T_q$ denote the lower triangular
matrix in 
${\rm Mat}_{d+1}(\F)$ with $(i,j)$-entry
\begin{eqnarray*}
(-1)^j q^{j(1-i)}
\left[\begin{array}{c} i  \\
j \end{array} \right]_q
\end{eqnarray*}
for $0 \leq j\leq i \leq d$.
\end{definition}

\begin{example}\rm For $d=3$,
\begin{eqnarray*}
&&
T_q=
\left(
\begin{array}{ cccc}
1 & 0 &0&0  \\
1 &  -1& 0 &0\\
1 & -q^{-1}\lbrack 2 \rbrack_q &  q^{-2} & 0\\
1 &-q^{-2}\lbrack 3 \rbrack_q &q^{-4} \lbrack 3 \rbrack_q &  -q^{-6} 
\end{array}
\right).
\end{eqnarray*}
\end{example}

\noindent By \cite[Theorem~15.4]{fduqe}
we have 
$(T_q)^{-1} = T_{q^{-1}}$.

\medskip

\noindent We are done discussing the matrices
(\ref{eq:three}).
We now consider the matrices
used in Section 11. They are
\begin{eqnarray}
\label{eq:seven}
K_q, \quad
E_q, \quad
F_q(t), \quad
G_q(t), \quad
L_q(t), \quad
S_q(t), \quad
M_q(t).
\end{eqnarray}

\begin{definition}
\label{Adef:Kq}
\rm Let $K_q$ denote the diagonal matrix
in 
${\rm Mat}_{d+1}(\F)$ with
$(i,i)$-entry $q^{d-2i}$ for 
$0 \leq i \leq d$.
\end{definition}

\begin{example} \rm For
$d=3$,
\begin{eqnarray*}
K_q = \mbox{\rm diag}( q^3, q, q^{-1}, q^{-3}).
\end{eqnarray*}
\end{example}

\begin{lemma}
\label{lem:ki}
We have
\begin{eqnarray*}
(K_q)^{-1} = K_{q^{-1}} = Z K_q Z.
\end{eqnarray*}
\end{lemma}


\begin{definition}
\label{Adef:eq}
\rm
Let 
$E_q$  denote the upper bidiagonal matrix
in 
${\rm Mat}_{d+1}(\F)$ 
with $(i,i)$-entry
$q^{2i-d}$ for $0 \leq i \leq d$ and
$(i-1,i)$-entry $q^d-q^{2i-2-d}$
for $1 \leq i \leq d$.
Note that $E_q$ has constant row sum $q^d$.
\end{definition}

\begin{example}\rm For $d=3$,
\begin{eqnarray*}
&&
E_q=
\left(
\begin{array}{ cccc}
q^{-3} & q^3-q^{-3}  &0&0  \\
0 &  q^{-1}& q^3-q^{-1} &0\\
0 & 0 &  q &  q^3-q \\
 0&0&0 &  q^3 
\end{array}
\right).
\end{eqnarray*}
\end{example}

\begin{lemma} 
\label{lem:eqform1}
We have
$E_q = X^{-1}K_q X$ where $X=T_q Z$. 
\end{lemma}

\noindent In Section 11 we refer to $(E_q)^{-1}$. 
\begin{lemma} 
The matrix $(E_q)^{-1}$ is upper triangular with
$(i,j)$-entry
\begin{eqnarray*}
(q^{2(d-j+1)};q^2)_{j-i} q^{d-2j}
\end{eqnarray*}
for $0 \leq i \leq j \leq d$.
\end{lemma}

\begin{example}\rm
For $d=3$,
\begin{eqnarray*}
&&
(E_q)^{-1}=
\left(
\begin{array}{ cccc}
q^{3} & (1-q^{6})q  &(1-q^4)(1-q^6)q^{-1}&(1-q^2)(1-q^4)(1-q^6)q^{-3}  \\
0 &  q& (1-q^4)q^{-1} &(1-q^2)(1-q^4)q^{-3}\\
0 & 0 &  q^{-1} &  (1-q^2)q^{-3} \\
 0&0&0 &  q^{-3} 
\end{array}
\right).
\end{eqnarray*}
\end{example}

\begin{definition}
\label{Adef:fq}
\rm
Let 
$F_q(t)$  denote the upper bidiagonal matrix
in 
${\rm Mat}_{d+1}(\F)$ 
with $(i,i)$-entry
$q^{2i-d}$ for $0 \leq i \leq d$ and
$(i-1,i)$-entry $(q^d-q^{2i-2-d})q^{1-d}t$
for $1 \leq i \leq d$.
\end{definition}

\begin{example}\rm For $d=3$,
\begin{eqnarray*}
&&
F_q(t)=
\left(
\begin{array}{ cccc}
q^{-3} & (q^3-q^{-3})q^{-2}t  &0&0  \\
0 &  q^{-1}& (q^3-q^{-1})q^{-2}t &0\\
0 & 0 &  q &  (q^3-q)q^{-2}t \\
 0&0&0 &  q^3 
\end{array}
\right).
\end{eqnarray*}
\end{example}

\begin{lemma} We have
\begin{eqnarray}
F_q(t) =
\mathcal D_{q}(t^{-1})
E_q 
\bigl(\mathcal D_{q}(t^{-1})\bigr)^{-1}
\label{eq:FE}
\end{eqnarray}
and also
\begin{eqnarray*}
F_q(t) = K_{q^{-1}} - t
\frac{\lbrack E_{q^{-1}}, K_{q^{-1}}\rbrack}{q-q^{-1}}.
\end{eqnarray*}
\end{lemma}

\begin{definition}
\label{Adef:Gqt}
\rm 
Let 
$G_q(t)$  denote the upper bidiagonal matrix
in 
${\rm Mat}_{d+1}(\F)$ 
with $(i,i)$-entry
$q^{2i-d}$ for $0 \leq i \leq d$ and
$(i-1,i)$-entry $(q^d-q^{2i-2-d})(1-t q^{d-2i+1})$
for $1 \leq i \leq d$.
\end{definition}

\begin{example}\rm For $d=3$,
\begin{eqnarray*}
&&
G_q(t)=
\left(
\begin{array}{ cccc}
q^{-3} & (q^3-q^{-3})(1-tq^{2})  &0&0  \\
0 &  q^{-1}& (q^3-q^{-1})(1-t) &0\\
0 & 0 &  q &  (q^3-q)(1-tq^{-2}) \\
 0&0&0 &  q^3 
\end{array}
\right).
\end{eqnarray*}
\end{example}

\begin{lemma} We have
\begin{eqnarray}
\label{eq:GE}
G_q(t) = 
\bigl(D_{q}(t)\bigr)^{-1}
E_q 
 D_q (t)
\end{eqnarray}
and also
\begin{eqnarray*}
G_q(t) = E_q - K_q + F_{q^{-1}}(t).
\end{eqnarray*}
\end{lemma}


\begin{definition}
\label{Adef:Lqt}
\rm 
Let 
$L_q(t)$  denote the upper bidiagonal matrix
in 
${\rm Mat}_{d+1}(\F)$ 
with $(i,i)$-entry
$q^{2i-d}$ for $0 \leq i \leq d$ and
$(i-1,i)$-entry $(q^d-q^{2i-2-d})(1-tq^{d-2i+1})^{-1}$
for $1 \leq i \leq d$.
\end{definition}

\begin{example}\rm For $d=3$,
\begin{eqnarray*}
&&
L_q(t)=
\left(
\begin{array}{ cccc}
q^{-3} & \frac{q^3-q^{-3}}{1-tq^{2}}  &0&0  \\
0 &  q^{-1}& \frac{q^3-q^{-1}}{1-t} &0\\
0 & 0 &  q &  \frac{q^3-q}{1-tq^{-2}} \\
 0&0&0 &  q^3 
\end{array}
\right).
\end{eqnarray*}
\end{example}

\begin{lemma}
\label{lem:lqi}
We have
\begin{eqnarray*}
L_q(t) = D_q(t) E_q \bigl(D_q(t)\bigr)^{-1}.
\end{eqnarray*}
\end{lemma}

\noindent In Section 11 we refer to $\bigl(L_q(t)\bigr)^{-1}$.
\begin{lemma} The matrix 
 $\bigl(L_q(t)\bigr)^{-1}$ is upper triangular with
 $(i,j)$-entry
\begin{eqnarray*}
\frac{(q^{2(d-j+1)};q^2)_{j-i}}{(tq^{d-2j+1};q^2)_{j-i}} q^{d-2j}
\end{eqnarray*}
for $0 \leq i \leq j\leq d$.
\end{lemma}

\begin{example} \rm 
For $d=3$,
\begin{eqnarray*}
 \bigl(L_q(t)\bigr)^{-1}
 = 
\left(
\begin{array}{ cccc}
q^3 &
\frac{(1-q^6)q}{1-tq^{2}} 
&
\frac{(1-q^4)(1-q^6)q^{-1}}{(1-t)(1-tq^2)} 
&
\frac{(1-q^{2})(1-q^{4})(1-q^{6})q^{-3}}{(1-tq^{-2})(1-t)(1-tq^2)} 
\\
0
& 
q
&
\frac{(1-q^4)q^{-1}}{1-t} 
&
\frac{(1-q^{2})(1-q^{4})q^{-3}}{(1-tq^{-2})(1-t)} 
\\
0 &
0
&
q^{-1}
& 
\frac{(1-q^{2})q^{-3}}{1-tq^{-2}} 
\\
0& 0 & 
0
&
q^{-3}
\end{array}
\right)
\end{eqnarray*}
\end{example}




\begin{definition}
\label{Adef:Sqt}
\rm
Let $S_q(t)$ denote the tridiagonal matrix in
${\rm Mat}_{d+1}(\F)$ 
with $(i-1,i)$-entry
$(q^d-q^{2i-d-2})(1-tq^{2i-d-1})$
for $1 \leq i \leq d$,
$(i,i-1)$-entry
$(q^{2i-d}-q^{-d})q^{2i-d-1}t$
for $1 \leq i \leq d$, and
$(i,i)$-entry
\begin{eqnarray*}
q^d - 
(q^{2i-d}-q^{-d})q^{2i-d-1}t
-
(q^d-q^{2i-d})(1-tq^{2i-d+1})
\end{eqnarray*}
for $0 \leq i \leq d$.
Note that $S_q(t)$ has constant row sum $q^d$.
\end{definition}

\begin{example}\rm
For $d=3$,

{\small
\begin{eqnarray*}
&&S_q(t)=
\\
&&
E_q+ t\left(
\begin{array}{ cccc}
q-q^{-5}  & q^{-5}-q  &0&0  \\
q^{-3}-q^{-5} & (q-q^{-1})(q^2+1-q^{-4})
& q^{-1}-q^3 &0\\
0 & q-q^{-3} &  
(q-q^{-1})(q^4-1-q^{-2})
& q^3-q^5  \\
 0&0& q^5-q^{-1}  &  
q^{-1}-q^5
\end{array}
\right)
\end{eqnarray*}
}
\end{example}

\begin{lemma} We have
\begin{eqnarray*}
S_q(t) = T_q G_{q^{-1}}(t) T_{q^{-1}}
\end{eqnarray*}
\noindent and also
\begin{eqnarray*}
S_q(t) = E_q - t\,\frac{\lbrack E_q, Z E_{q^{-1}}Z\rbrack}{q-q^{-1}}.
\end{eqnarray*}
\end{lemma}



\begin{definition}
\label{Adef:M}
\rm We define
$M_q(t) \in
{\rm Mat}_{d+1}(\F)$ as follows.
For $0 \leq i,j\leq d$ the 
   $(i,j)$-entry of
$M_q(t)$ is given in the table below. 

\medskip
\centerline{
\begin{tabular}[t]{c|c}
{\rm case } & 
   {\rm $(i,j)$-entry of $M_q(t)$} 
   \\ \hline 
   \\
    $i-j>1$  
&
   $0$
\\
\\
$i-j=1$ 
&
  $\frac{q^{-1}-q^{2i-1}}{t^{-1}-q^{2i-d-1}}$
  \\
\\
$j\geq i, \quad i\not=0, \quad j\not=d$ 
&
 $
 \frac{1-tq^{d+1}}{1-tq^{2j-d+1}}
\, 
\frac{1-tq^{-d-1}}{1-tq^{2i-d-1}}
\,
\frac{(q^{2i-2d};q^2)_{j-i}}{(tq^{2i-d+1};q^2)_{j-i}}
\, 
 \frac{1}{q^{d-2j}} $
  \\
\\
 $i=0, \quad j\not=d$ 
&
 $
 \frac{1-tq^{d+1}}{1-tq^{2j-d+1}}
\, 
\frac{(q^{-2d};q^2)_{j}}{(tq^{1-d};q^2)_{j}}
\, 
 \frac{1}{q^{d-2j}} $
  \\
\\
$i\not=0, \quad j=d$ 
&
 $
\frac{1-tq^{-d-1}}{1-tq^{2i-d-1}}
\,
\frac{(q^{2i-2d};q^2)_{d-i}}{(tq^{2i-d+1};q^2)_{d-i}}
\, 
 {q^d} $
  \\
\\
 $i=0,  \quad j=d$ 
&
 $
\frac{(q^{-2d};q^2)_{d}}{(tq^{1-d};q^2)_{d}}
\, 
 {q^d} $
     \end{tabular}}
     \medskip

\end{definition}

\begin{example} \rm
For $d=3$,
\begin{eqnarray*}
M_q(t) = 
\left(
\begin{array}{ cccc}
\frac{1-tq^4}{q^3(1-tq^{-2})} &
\frac{(1-tq^4)(1-q^{-6})}{q(1-tq^{-2})(1-t)} 
&
\frac{q(1-tq^4)(1-q^{-4})(1-q^{-6})}{(1-tq^{-2})(1-t)(1-tq^2)} 
&
\frac{q^3(1-q^{-2})(1-q^{-4})(1-q^{-6})}{(1-tq^{-2})(1-t)(1-tq^2)} 
\\
\frac{q^{-1}-q}{t^{-1}-q^{-2}}
& 
\frac{(1-tq^4)(1-tq^{-4})}{q(1-tq^{-2})(1-t)} 
&
\frac{q(1-tq^4)(1-tq^{-4})(1-q^{-4})}{(1-tq^{-2})(1-t)(1-tq^2)} 
&
\frac{q^3(1-tq^{-4})(1-q^{-2})(1-q^{-4})}{(1-tq^{-2})(1-t)(1-tq^2)} 
\\
0 &
\frac{q^{-1}-q^3}{t^{-1}-1}
&
\frac{q(1-tq^4)(1-tq^{-4})}{(1-t)(1-tq^2)} 
& 
\frac{q^3(1-tq^{-4})(1-q^{-2})}{(1-t)(1-tq^2)} 
\\
0& 0 & 
\frac{q^{-1}-q^5}{t^{-1}-q^2}
&
\frac{q^3(1-tq^{-4})}{1-tq^2}
\end{array}
\right)
\end{eqnarray*}
\end{example}

\begin{lemma} We have
\begin{eqnarray*}
M_q(t) =  T_{q^{-1}} \bigl(L_q(t)\bigr)^{-1} T_q.
\end{eqnarray*}
\end{lemma}

%
%

\section{Appendix II: Some matrix identities}

\noindent In this section we 
record some miscellaneous facts about the matrices
listed in (\ref{eq:three}),
(\ref{eq:seven}).

\begin{lemma}\label{lem:9pt4}
We have
\begin{eqnarray*}
&& 
t\bigl(ZE_qZ-F_q(t^{-1})\bigr) =
\frac{\lbrack E_{q^{-1}},ZF_{q^{-1}}(t)Z\rbrack}{q-q^{-1}},
\\
&&
t\bigl(G_q(t^{-1})-Z E_{q^{-1}}Z\bigr) =
\frac{\lbrack S_q(t),K_{q}\rbrack}{q-q^{-1}},
\\
&&
t\bigl(K_q-S_q(t^{-1})\bigr) =
\frac{\lbrack G_q(t),ZE_{q^{-1}}Z\rbrack}{q-q^{-1}}.
\end{eqnarray*}
\end{lemma}
\noindent {\it Proof:} 
These equations express how 
the relations from 
Lemma \ref{lem:eval} look in 
the 24 bases for
$\mathbf V_d(t)$
from Lemma
\ref{lem:normexist}.
\hfill $\Box$ \\

\begin{lemma} We have
\begin{eqnarray*}
&& 
t\bigl(ZE_qZ-M_q(t)\bigr) =
\frac{\lbrack E_{q^{-1}},M_{q}(t)\rbrack}{q-q^{-1}},
\qquad 
t^{-1}\bigl(F_{q^{-1}}(t)-K_q\bigr) =
\frac{\lbrack E_q,K_{q}\rbrack}{q-q^{-1}},
\\
&&
t\bigl(E_{q^{-1}}-K_q\bigr) =
\frac{\lbrack F_q(t),K_{q}\rbrack}{q-q^{-1}},
\quad 
 t^{-1}\bigl(ZG_{q^{-1}}(t)Z-\bigl(L_q(t)\bigr)^{-1} \bigr) =
\frac{\lbrack \bigl(L_q(t)\bigr)^{-1},ZS_{q^{-1}}(t^{-1})Z\rbrack}{q-q^{-1}},
\\
&&
t\bigl(K_q-(E_q)^{-1}\bigr) =
\frac{\lbrack G_q(t),(E_q)^{-1}\rbrack}{q-q^{-1}},
\qquad
t^{-1}\bigl(S_q(t)-E_q \bigr) =
\frac{\lbrack Z E_{q^{-1}} Z,E_q \rbrack }{q-q^{-1}},
\\
&&
t\bigl(S_q(t^{-1})-ZL_{q^{-1}}(t)Z \bigr) =
\frac{\lbrack Z L_{q^{-1}}(t)Z,G_q(t)\rbrack}{q-q^{-1}},
\qquad
t\bigl(G_q(t^{-1})-E_q \bigr) =
\frac{\lbrack E_q,K_q \rbrack}{q-q^{-1}},
\\
&&
t\bigl(ZE_{q^{-1}}Z-(E_q)^{-1} \bigr) =
\frac{\lbrack (E_q)^{-1},S_q(t)\rbrack}{q-q^{-1}},
\qquad
t\bigl(K_q-\bigl(L_q(t)\bigr)^{-1} \bigr) =
\frac{\lbrack E_q,\bigl(L_q(t)\bigr)^{-1}\rbrack }{q-q^{-1}},
\\
&&
t^{-1}\bigl(E_q-L_q(t) \bigr) =
\frac{\lbrack L_q(t),K_q\rbrack }{q-q^{-1}},
\qquad
 t^{-1}\bigl(ZF_{q^{-1}}(t)Z-M_q(t)\bigr) =
 \frac{\lbrack M_q(t),F_{q}(t^{-1})\rbrack}{q-q^{-1}}.
\end{eqnarray*}
\end{lemma}
\noindent {\it Proof:} 
These equations express how 
the relations from 
Lemma
\ref{lem:teq2}
 look in 
the 24 bases for
$\mathbf V_d(t)$
from Lemma
\ref{lem:normexist}.
\hfill $\Box$ \\

\begin{lemma}
\label{lem:eqformcomb}
We have
\begin{eqnarray*}
&&
T_q Z T_q Z T_q Z = (-1)^d q^{-d(d-1)}I,
\\
&&T_q \mathcal D_q(t) 
T_q D_{q^{-1}}(t^{-1}) T_{q^{-1}} 
\bigl(D_q(t)\bigr)^{-1}
= I.
\end{eqnarray*}
\end{lemma}
\noindent {\it Proof:} 
For $n \in \N$, let
$B_0,B_1,\ldots,B_n$ denote a sequence of bases for $\mathbf V_d(t)$
such that $B_0=B_n$. For $1 \leq i \leq n$ let
$T_i$ denote the transition matrix from $B_{i-1}$ to $B_i$.
Then $T_1T_2\cdots T_n=I$. Both equations in the lemma
statement are obtained in this way, using an appropriate
sequence of bases from the 24 given in Lemma
\ref{lem:normexist}.
\hfill $\Box$ \\

\begin{lemma} We have
\begin{eqnarray*}
&&
Z T_q Z G_q(t) = F_q(t) Z T_q Z,
\quad
\qquad 
\qquad 
\qquad
E_{q^{-1}} = 
T_{q^{-1}}E_q T_q,
\\
&&
L_q(t)ZT_{q^{-1}}Z L_{q^{-1}}(t^{-1}) Z T_q Z = I,
\qquad 
\qquad 
S_q(t) = Z T_{q^{-1}} F_{q^{-1}}(t) T_q Z,
\\
&&
D_q(t^{-1})S_q(t) \bigl(D_q(t^{-1})\bigr)^{-1} = Z S_{q^{-1}}(t)Z,
\qquad 
M_q(t) = Z T_q L_{q^{-1}}(t^{-1}) T_{q^{-1}}Z,
\\
&&
M_{q^{-1}}(t^{-1}) = \bigl(\mathcal D_q(t)\bigr)^{-1} M_q(t) 
\mathcal D_q(t),
\qquad 
\qquad 
M_q(t) Z M_{q^{-1}}(t^{-1}) Z = I. 
\end{eqnarray*}
\end{lemma}
\noindent {\it Proof:} 
These are all applications of the linear
algebra principle from the fourth paragraph of
Section 13.
\hfill $\Box$ \\


\begin{lemma} In the table below we list some 3-tuples $u,v,w$
of matrices in
${\rm Mat}_{d+1}(\F)$.
For each case
\begin{eqnarray*}
\frac{quv-q^{-1}vu}{q-q^{-1}}=I,
\qquad
\frac{qvw-q^{-1}wv}{q-q^{-1}}=I,
\qquad
\frac{qwu-q^{-1}uw}{q-q^{-1}}=I.
\end{eqnarray*}


\centerline{
\begin{tabular}[t]{ccc}
\\
$u$ & $v$ & $w$
\\
\hline
$E_q$ &  $K_q$ &  $ZE_{q^{-1}}Z$
\\
$L_q(t)$ & $K_q$ & $ZG_{q^{-1}}(t)Z$
\\
$S_{q^{-1}}(t)$ & 
$(E_{q^{-1}})^{-1}$ &
$G_{q^{-1}}(t^{-1})$
\\
$E_q$ 
&
$\bigl(L_q(t^{-1})\bigr)^{-1}$ &
$Z S_{q^{-1}}(t)Z$ 
\\
$Z G_q(t)Z$ &   $K_{q^{-1}}$ &  $L_{q^{-1}}(t)$
\\
$G_q(t^{-1})$ & $(E_q)^{-1}$ & $S_q(t)$
\\
$ZS_q(t)Z$ & 
 $\bigl(L_{q^{-1}}(t^{-1})\bigr)^{-1}$ &
$E_{q^{-1}}$
\\
$F_q(t^{-1})$ &
$E_{q^{-1}}$ &
$M_q(t)$
\\
$M_{q^{-1}}(t)$ &
$E_q$ &
$F_{q^{-1}}(t^{-1})$
\\
$F_q(t^{-1})$ &
$K_q$ &
$ZF_{q^{-1}}(t) Z$
     \end{tabular}}
\medskip

\end{lemma}
\noindent {\it Proof:} 
The relations 
(\ref{eq:tet2})
hold in $\boxtimes_q$.
The equations in the present
lemma
express how 
these relations
 look in 
the 24 bases for
$\mathbf V_d(t)$
from Lemma
\ref{lem:normexist}.
\hfill $\Box$ \\

\section{Acknowledgments}
The author thanks 
Hau-wen Huang and 
Kazumasa Nomura
for some insightful comments that lead to improvements
in the paper.

\bigskip

\noindent Tatsuro Ito \hfil\break
\noindent Division of Mathematical and Physical Sciences \hfil\break
\noindent Graduate School of Natural Science and Technology \hfil\break
\noindent Kanazawa University \hfil\break
\noindent Kakuma-machi, Kanazawa 920-1192, Japan \hfil\break
\noindent email: {\tt ito@se.kanazawa-u.ac.jp}\hfil\break
\bigskip

\noindent Hjalmar Rosengren \hfil\break
\noindent Department of Mathematical Sciences  \hfil\break
\noindent Chalmers University of Technology  \hfil\break
\noindent SE-412 96 G\"oteborg, Sweden \hfil\break
\noindent email: {\tt hjalmar@chalmers.se}\hfil\break 

\bigskip

\noindent Paul Terwilliger \hfil\break
\noindent Department of Mathematics \hfil\break
\noindent University of Wisconsin \hfil\break
\noindent 480 Lincoln Drive \hfil\break
\noindent Madison, WI 53706-1388 USA \hfil\break
\noindent email: {\tt terwilli@math.wisc.edu }\hfil\break


\end{document}